\documentclass[3p,times]{elsarticle}

\usepackage{graphicx}
\usepackage{psfrag}
\usepackage[dvipsnames]{xcolor}

\usepackage{amsmath}
\usepackage{amsthm}
\usepackage{amsfonts}
\usepackage{dsfont}
\usepackage{rotating}
\usepackage{multirow}
\usepackage{epstopdf}
\usepackage{newtxmath}
\usepackage{bm}
\usepackage{bbm}
\usepackage{comment}
\usepackage{booktabs}
\usepackage{accents}
\usepackage{empheq}
\newtheorem*{remark}{Remark}

\DeclareMathAlphabet\mathbfcal{OMS}{cmsy}{b}{n}  

\renewcommand{\vv}{\mathbf{v}}
\newcommand{\xx}{\mathbf{x}}

\newcommand{\q}{\mathbf{q}}
\newcommand{\bss}{\boldsymbol{\sigma}}

\newcommand{\dt}{{\Delta t}}

\newcommand{\Q}{\mathbf{U}}
\newcommand{\F}{\mathbf{F}}
\newcommand{\I}{\mathbf{I}}
\newcommand{\U}{\mathbf{U}}
\newcommand{\f}{\mathbf{f}}
\newcommand{\g}{\mathbf{g}}

\newcommand{\M}{\textnormal{Ma}}
\newcommand{\CFL}{\textnormal{CFL}}

\newcommand{\Pv}{\mathbf{\P}}
\newcommand{\uu}{\mathbf{u}}

\newcommand{\Rey}{\text{Re}}

\newcommand{\vbase}{\phi} 

\newcommand{\xv}{\mathbf{x}}

\newcommand{\xs}{x}
\newcommand{\ys}{y}
\newcommand{\vs}{\varv}

\newcommand{\hh}{h}

\renewcommand{\P} {{P}}            
\newcommand  {\E} {{e}}            
\newcommand{\xvP}{\xv_{\P}}        

\newcommand{\hP}{\hh_{\P}}

\newcommand{\Pset}{{P}}    

\newcommand{\NMB}{N}
\newcommand{\NP}{\NMB_{\Pset}}   

\newcommand{\abs}    [1]{|#1|}

\newcommand{\nor}  {\mathbf{n}}

\newcommand{\HONE}  {H^1}
\newcommand{\VS}[1] {V^{#1}}
\newcommand{\PS}[1]{\mathbbm{P}_{#1}}
\newcommand{\CS}[1] {C^{#1}}

\newcommand{\Vh}  {\VS{\hh}}

\newcommand{\vsh}{\vs}

\newcommand{\NDOF}{N^{\textrm{dof}}_{\P}}
\newcommand{\NDOFST}{N^{\textrm{dofst}}_{\P}}

\newcommand{\proj}{\Pi^\nabla_{\P,N}}
\newcommand{\projL}{\Pi^0_{\P,N}}

\begin{document}
	
	\begin{frontmatter}
		
		\journal{Elsevier}
		
		

		\title{Nonconforming Virtual Element basis functions for space-time Discontinuous Galerkin schemes on unstructured Voronoi meshes}
		
		\author[dmi]{Walter Boscheri$^*$}
		\ead{walter.boscheri@unife.it}
		\cortext[cor1]{Corresponding author}
        
        \author[disap]{Giulia Bertaglia}
        \ead{giulia.bertaglia@unife.it}
        
		\address[dmi]{Department of Mathematics and Computer Science, University of Ferrara, Via Niccol\`o Machiavelli 30, 44121 Ferrara, Italy}
		\address[disap]{Department of Environmental and Prevention Sciences, University of	Ferrara, Corso Ercole I d'Este 32, 44121 Ferrara, Italy}
%
\begin{abstract}
We introduce a new class of Discontinuous Galerkin (DG) methods for solving nonlinear conservation laws on unstructured Voronoi meshes that use a nonconforming Virtual Element basis defined within each polygonal control volume. The basis functions are actually evaluated as an $L_2$ projection of the virtual basis which remains unknown, along the lines of the Virtual Element Method (VEM). Contrarily to the VEM approach, the new basis functions lead to a nonconforming representation of the solution with discontinuous data across the element boundaries, as typically employed in DG discretizations. To improve the condition number of the resulting mass matrix, an orthogonalization of the full basis is proposed. The discretization in time is carried out following the ADER (Arbitrary order DERivative Riemann problem) methodology, which yields one-step fully discrete schemes that make use of a coupled space-time representation of the numerical solution. The space-time basis functions are constructed as a tensor product of the virtual basis in space and a one-dimensional Lagrange nodal basis in time. The resulting space-time stiffness matrix is stabilized by an extension of the dof--dof stabilization technique adopted in the VEM framework, hence allowing an element-local space-time Galerkin finite element predictor to be evaluated. The novel methods are referred to as VEM-DG schemes, and they are arbitrarily high order accurate in space and time.

The new VEM-DG algorithms are rigorously validated against a series of benchmarks in the context of compressible Euler and Navier--Stokes equations. Numerical results are verified with respect to literature reference solutions and compared in terms of accuracy and computational efficiency to those obtained using a standard modal DG scheme with Taylor basis functions. An analysis of the condition number of the mass and space-time stiffness matrix is also forwarded.
\end{abstract}
%
\end{frontmatter}


\section{Introduction} \label{sec.intro}
Discontinuous Galerkin (DG) methods have been first introduced in \cite{reed} for the solution of neutron transport equations, and subsequently applied to general nonlinear systems of hyperbolic conservation laws in \cite{cbs0,cbs1,cbs2,cbs3}. In the DG framework, the numerical solution is represented by piecewise polynomials within each mesh control volume, allowing jumps of the discrete solution across element boundaries. The main advantage of the DG approach is that it automatically provides high order of accuracy \emph{locally}, thus it does not require any reconstruction stencil like in finite volume solvers. Furthermore, DG schemes are typically more accurate compared to finite volume or finite difference methods, because the entire high order polynomial is evolved in time for each computational cell.

The discrete solution is represented in terms of an expansion that involves a set of basis functions and the associated degrees of freedom, also referred to as expansion coefficients. The basis can be either modal or nodal, provided that it achieves the formal order of accuracy of the method. In the nodal approach, the degrees of freedom correspond to the value of the numerical solution at the nodal points, while in the modal approach the expansion coefficients give the modes of the polynomial basis. 

Efficient DG schemes can be devised with \textit{nodal basis} if the nodes are carefully chosen, hence leading to nodal DG schemes with Gauss-Lobatto nodes \cite{GassnerDG_LES} or Gauss-Legendre nodes \cite{Exahype}. These schemes are typically referred to as spectral element methods (SEM) \cite{Kopriva2009}, which fit the Summation-By-Parts Simultaneous-Approximation-Term (SBP-SAT) framework \cite{CARPENTER199674,GassnerSIAM2013}. The nodal basis functions are defined on a reference element, to which the physical control volume is mapped. This means that the usage of nodal basis is restricted to either Cartesian meshes with quadrilaterals and hexahedra, or unstructured grids made of simplex control volumes,
namely triangles and tetrahedra. To overcome this limitation, in \cite{GassnerPoly} a nodal basis ansatz on more general unstructured meshes is proposed, and recently in \cite{ADERAFEDG} an agglomerated continuous finite element basis is devised at the sub-grid level of a Voronoi mesh. The nonlinear stability of nodal DG schemes is guaranteed by employing either artificial viscosity techniques \cite{PerssonAV,VegtAV,TavelliCNS,GassnerMunzAV} or sub-cell finite volume limiters \cite{DGLimiter1,DGsubcell_Gassner}.

Alternatively, \textit{modal basis} like Taylor basis can be easily employed on very general control volumes \cite{ArepoTN,DGBoltz}, although the associated computational cost increases since no reference element is available. Nevertheless, the adoption of a hierarchical modal approach is very convenient for designing slope and moment limiters in order to ensure the stability of DG schemes \cite{cbs4,Biswas_94,Burbeau_2001,Kri07,Kuzmin2013}.     

With the aim of dealing with general polygonal and polyhedral elements, the Virtual Element Method (VEM), which belongs to the class of continuous finite element methods, has recently emerged \cite{vem1,vem2,vem3,vem4,vem5,vem6}. In the VEM framework, the discrete solution is approximated by a set of basis functions that do not need to be explicitly determined, hence making them only \emph{virtually} defined. Indeed, the numerical solution is known using suitable operators that project the basis functions onto polynomial spaces of any degree, allowing for the discretization and appropriate approximation of the continuous linear functional and the bilinear form resulting from the variational formulation only through the knowledge of the degrees of freedom on the boundary of the mesh elements. Several applications in the field of solid mechanics, especially for elliptic problems of linear elasticity \cite{vem3,Gain2014}, linear elastodynamics \cite{Antonietti2021}, elastic problems with singularities and discontinuities \cite{BenvenutiChiozzi2019,BenvenutiChiozzi2022}, and fracture mechanics \cite{Hussein2019,Nguyen2018}, have largely proposed VEM strategies to solve the problem of interest. In more recent years, applications of VEM in the field of fluid dynamics have been forwarded, such as for flow problems in porous media \cite{Borio2022,Borio2021} or to solve the steady Navier--Stokes equations \cite{Beirao-Lovadina2018,Chernov2021,BeiraoStokes22,vem6,Antonietti22}. To account for discontinuous solution across the vertexes of the computational mesh, \emph{nonconforming} virtual element methods have been introduced in \cite{VEM_nc_org,vem_nc_2,vem_nc_3,vem_nc_4,vem_nc_5}. In this case, the numerical solution is always computed by solving a linear system on the entire computational mesh, but at the global level the conformity requirement between elements is relaxed, and the definition of the global discrete space is no longer continuous.

In this work, we further extend the application of VEM to solve nonlinear time-dependent Partial Differential Equations (PDEs), introducing an innovative class of high-order numerical methods which relate the Virtual Element Method to the compact framework of Discontinuous Galerkin Finite Element methods on unstructured Voronoi meshes, focusing on the solution of compressible viscous flows governed by the Navier--Stokes equations. The VEM projection operators are used to define the \emph{local} solution space within each control volume, hence using the virtual basis functions to approximate the discrete solution in the cell. Belonging the numerical scheme to the family of DG methods, the solution remains discontinuous across element boundaries, leading to a numerical method based on nonconforming virtual element basis functions. The new basis is a mixed nodal/modal basis, since it involves high order internal momenta starting from third order approximations. However, the interpolation property exhibited by nodal basis holds even for the degrees of freedom associated to the internal momenta, making the new nonconforming VEM basis potentially able to yield a quadrature-free DG scheme \cite{Exahype,ADERAFEDG}. The time discretization is performed relying on the ADER (Arbitrary order DERivative Riemann problem) approach \cite{Toro2006,Dumbser2008}, which, in contrast to Runge-Kutta schemes based on temporal sub-stages evaluations, permits to obtain high order in time with a one-step fully discrete method. This is achieved through the computation of a space-time predictor solution at the aid of an element-local finite element predictor. Therefore, we construct the space-time basis functions by performing a tensor product between the nonconforming VEM basis in space and a standard nodal Lagrange basis in time. As analyzed in \cite{Mascotto2018}, the mass and stiffness matrices arising from the VEM paradigm are highly ill-conditioned, thus requiring special care if they have to be inverted. Consequently, an \emph{ad hoc} stabilization technique is also introduced for the inversion of the space-time stiffness matrix needed in the ADER predictor strategy. Inspired by \cite{Berrone2017}, to further reduce the condition number of the VEM spatial mass matrix, an orthogonalization technique is proposed on the full basis.

The rest of the paper is structured as follows. In Section \ref{sec.pde} we introduce the governing equations. Section \ref{sec.numscheme} presents the numerical method, including the description of the novel nonconforming space-time VEM basis as well as the fully discrete DG scheme. The accuracy and the robustness of the proposed approach are demonstrated in Section \ref{sec.validation}, which contains a set of benchmark test problems in the field of compressible inviscid and viscous flows. Finally, Section \ref{sec.concl} provides some concluding remarks and gives an outlook to future work.

\section{Mathematical model} \label{sec.pde}
The mathematical model is given by a nonlinear system of conservation laws, reading in general form as
\begin{equation}
	\label{eqn.PDE}
	\frac{\partial \Q}{\partial t} + \nabla \cdot \F(\Q,\nabla \Q) = \mathbf{0}\,, \qquad 
	\xx \in \Omega \subset \mathds{R}^d\,, \quad 
	t \in \mathds{R}_0^+\,, \quad 
	\Q \in \Omega_{\Q} \subset \mathds{R}^{\gamma}\,,
\end{equation}
where $\Omega$ is a bounded domain in $d=2$ space dimensions, $\xx = (x,y)$ is the vector of spatial coordinates, and $t$ is the time. Here, $\Q$ denotes the vector of conserved variables defined in the space of admissible states $\Omega_{\Q}\subset \mathds{R}^{\gamma}$, while $\F(\Q,\nabla \Q) = (\f(\Q,\nabla \Q), \,\g(\Q,\nabla \Q))$ is the conservative non-linear flux tensor ($\f$ in $x$-direction, $\g$ in $y$-direction), which depends not only on the conserved state $\Q$ but also on its gradient $\nabla \Q$. In particular, we are interested in the compressible Navier--Stokes equations for a Newtonian fluid with heat conduction, which are based on the physical principle of conservation of mass, momentum and total energy. They can be written in the conservative form \eqref{eqn.PDE}, being
\begin{equation}
	\label{eqn.NS}
	\begin{aligned}
	\Q = 	\begin{pmatrix}
			\, \rho\\ 
			\rho \vv \\ 
			\rho E \,
			\end{pmatrix}\,,\qquad
	\F(\Q,\nabla \Q) = 	\begin{pmatrix}
							\, \rho \vv \\ 
							\rho (\vv \otimes \vv) + \bss(\Q, \nabla \Q) \\ 
							\vv \cdot \left(\rho E \I + \bss(\Q, \nabla \Q)\right) - \kappa \nabla T\,
							\end{pmatrix}\,.
	\end{aligned}
\end{equation}
In the above equations, $\rho(\xx, t)$ is the fluid density, $\vv(\xx, t) = (u,v)$ denotes the fluid velocity vector, $p(\xx,t)$ is the fluid pressure, $E(\xx,t)=\varepsilon+\frac1{2}|\vv|^2$ defines the specific total energy as the sum of the specific internal energy $\varepsilon(\xx,t)$ and the specific kinetic energy $k(\xx,t)=\frac1{2}|\vv|^2$, and $\I$ is the identity $d \times d$ matrix. The stress tensor $\bss(\Q, \nabla \Q)$ is defined under Stokes hypothesis as follows:
\begin{equation}
\label{stressT}
\bss(\Q, \nabla \Q) = \left( p + \frac{2}{3}\mu\, \nabla \cdot \vv\right) \I - \mu\left( \nabla\vv + \nabla\vv^T \right)\,,
\end{equation}
where $\mu$ is the dynamic viscosity, assumed to be constant. In the energy flux, $T(\xx,t)$ represents the fluid temperature, while $\kappa=\mu \gamma c_v \Pr^{-1}$ is the thermal conductivity coefficient, which depends on the viscosity $\mu$, the Prandtl number $\Pr$ and the specific heats $c_v$ and $c_p$, at constant volume and pressure, respectively, being $\gamma = c_p/c_v$ the adiabatic index. The specific heat at constant volume $c_v = R/(\gamma - 1)$ depends on the gas constant $R$. 


\subsection{Equations of state}
To close the system defined by \eqref{eqn.NS}, it is necessary to introduce two equations of state (EOS), a thermal one for the pressure $p = p(T,\rho)$ and a caloric one for the specific internal energy $\varepsilon=\varepsilon(T,\rho)$. In this work, we consider an ideal gas with the thermal and caloric EOS provided, respectively, by
\begin{equation}
\label{EOS}
\frac{p}{\rho}=RT\,, \qquad \varepsilon=c_v T\,.
\end{equation}
With this assumption, using the definition of the temperature obtained inverting the thermal EOS into the caloric EOS, it is possible to cancel out the temperature and to work with a single EOS of the form $\varepsilon(p,\rho)$, which defines the following linear relationship between pressure $p$ and internal energy $\varepsilon$:
\begin{equation}
\label{EOS1}
\varepsilon(p,\rho)=\frac{p}{\rho(\gamma-1)}\,.
\end{equation}

\subsection{Eigenvalues of the system}
The compressible Navier--Stokes equations \eqref{eqn.NS} involve a convective-acoustic sub-system, which corresponds to the compressible Euler equations, and a viscous sub-system. The eigenvalues can be analyzed separately for each sub-system \cite{ADERAFEDG,BosPar2021}, thus the convective-acoustic eigenvalues $\boldsymbol{\lambda}^c=\{\lambda_i^c\}_{i=1}^4$ are given by
\begin{equation}
\label{eig_c}
\lambda^c_1 = |\vv| - c \,, \qquad \lambda^c_2 = \lambda^c_3 = |\vv|\,, \qquad \lambda^c_4 = |\vv| + c\,,
\end{equation}
where $c = \sqrt{\gamma RT}$ is the sound speed, while the viscous eigenvalues $\boldsymbol{\lambda}^v=\{\lambda_i^v\}_{i=1}^4$, which depend on the viscosity and heat conduction properties of the fluid, read
\begin{equation}
\label{eig_v}
\lambda^v_1 = \lambda^v_2 = 0 \,, \qquad
\lambda^v_3 = \frac{4\mu}{3\rho} \,, \qquad \lambda^v_4 = \frac{\gamma \mu}{\rho \Pr}\,.
\end{equation}

\section{Numerical scheme} \label{sec.numscheme}
In this work, we propose to solve numerically the governing equations \eqref{eqn.PDE}-\eqref{eqn.NS} using a high order fully discrete one-step scheme which embeds the Virtual Element Method (VEM) \cite{vem2} in the definition of the basis functions of a Discontinuous Galerkin (DG) approach for the discretization in space, accounting for unstructured Voronoi meshes. To attain high accuracy in time, the ADER (Arbitrary order DERivative Riemann problem) methodology \cite{Dumbser2008} is properly modified to deal with the novel spatial basis functions.

\subsection{Discretization of the space-time computational domain} \label{ssec.domain}
In this section we present the details of the discretization of the computational domain in space and time.

\subsubsection{Space discretization}
The spatial domain $\Omega \subset \mathds{R}^2$ is discretized in $N_P$ non-overlapping polygonal control volumes $P_i$, $i = 1,...,N_P$, with boundary $\partial P_i$ and surface area $|P_i|$, defined by the following Voronoi-type tessellation:
\begin{equation*}
\mathcal{T}_{\Omega} = \bigcup \limits_{i=1}^{N_P}{P_i}\,. 
\end{equation*}
Each Voronoi element $P_i$ is characterized by $N_{e_i}$ vertexes $r$ (oriented counter-clockwise) with coordinates $\xv_{r}=(\xs_{r},\ys_{r})$, and $N_{e_i}$ straight edges $e$ that compose the boundary $\partial P_i$. For each edge $e$ we also define the outward pointing unit normal vector $\nor_{e}=(n_x,n_y)_e$.
The barycenter ${\xv_{P_i}}=(\xs_{P_i},\ys_{P_i})$ of each element is defined so that
\begin{equation}
	\label{eqn.xi}
	\xv_{P_i} = \frac1{N_{e_i}} \sum_{k=1}^{N_{e_i}}\xv_{r_k}\,. 
\end{equation}
Each Voronoi cell is associated to a characteristic length size $h_{P_i}$ computed as follows:
\begin{equation}
	\label{eqn.h}
	h_{P_i} = \frac{2 |P_i|}{\sum_{k=1}^{N_{e_i}} |\partial P_{i \, e_k}|}\,,
\end{equation}
where $|\partial P_{i \, e_k}|$ is the length of the edge $e_k$, $k=1,...,N_{e_i}$, of the element $P_i$. 

To construct the Voronoi grid we start by a primary Delaunay triangulation, whose vertexes provide the generators of the Voronoi cells. The interested reader is referred to \cite{ADERAFEDG} for all the details. Although we use a Voronoi mesh, this is not mandatory for the numerical scheme presented in this work. Indeed, we remark that any unstructured mesh can be adopted, provided that the following mesh regularity assumptions \cite{SIFVVEM,Mascotto2018} are fulfilled:
there exists a constant $\varrho > 0$ such that for every polygonal element $P_i \in\mathcal{T}_\Omega$ it holds that
\begin{description}
	\item[$(i)$] $P_i$ is star-shaped with respect to a disk with radius 
	$R_i \geq \varrho h_{P_i}$; 
	\item[$(ii)$] for every edge $e_k \in\partial P_i$ it holds that
	$|\partial P_{i \, e_k}| \geq \varrho h_{P_i}$.
\end{description}
This implies that the number of edges and vertexes is finite over the whole computational mesh, and that degenerate edges with zero length do not exist.

\subsubsection{Time discretization}
The time coordinate is bounded in the interval $[0,t_f]$, where $t_f \in \mathds{R}_0^+$ identifies the final time of the simulation. A sequence of $\dt = t^{n+1}-t^n$ time steps are used to discretize the time interval, such that $t \in [t^n;t^{n+1}]$:
\begin{equation}
\label{eqn.dt}
t = t^n + \tau \dt\,, \qquad \tau \in [0;1]\,,
\end{equation} 
where $\tau$ defines a reference time coordinate which maps the physical time interval $[t^n;t^{n+1}]$ into the range $[0;1]$.
The time step size is computed at each temporal iteration in order to respect a $\CFL$ (Courant--Friedrichs--Lewy) stability condition for explicit DG schemes \cite{ADERNSE,stedg2}, thus
\begin{equation}
\label{CFL}
\dt = \frac{\CFL}{2N+1}\, \frac{\min \limits_{\mathcal{T}_\Omega} h_{P}}{\max \limits_{\mathcal{T}_\Omega} \left( |\lambda^c_{\max}| + \frac{2(2N+1)}{h_i} \cdot |\lambda^v_{\max}| \right)_{P}} \,,
\end{equation}
where $\lambda^c_{\max}$ and $\lambda^v_{\max}$ represent the maximum eigenvalues defined in \eqref{eig_c}-\eqref{eig_v} of the element $P$, and $N$ is the degree of the chosen piecewise polynomial data representation, as further detailed.

\subsection{Data representation} \label{ssec.numsol}
For ease of reading, from now on we omit the subscript $i$ referring to polygon $P_i$, hence simply writing $P$, bearing in mind that we address local quantities within $P_i$. In each computational cell $P$, the discrete representation of the conserved quantity $\Q(\xx,t)$ at the current time $t^n$ is referred to as $\uu_h^n$, and it is given in terms of piecewise polynomials of degree $N \geq 0$ defined in the space $\PS{N}$ as
\begin{equation} 
\uu_h^n:=\Q\left(\xx |_{P},t^n \right) = \sum \limits_{\ell=1}^{\mathcal{N}} \phi_\ell(\xx) \, \hat{\uu}^{n}_{\ell} 
:= \phi_\ell(\xx) \, \hat{\uu}^{n}_{\ell} , \qquad \xx \in {P},
\label{eqn.uh}
\end{equation}
with $\phi_\ell(\xx)$ denoting the spatial basis functions used to span the space of piecewise polynomials $\PS{N}$ up to degree $N$, and $\mathcal{N}$ being the total number of degrees of freedom $\hat{\uu}^{n}_{\ell}$ for the cell $P$. Classical tensor index notation based on the Einstein summation convention is adopted, which implies summation over two equal indexes.

To compute the quantity $\uu_h^n$, an explicit definition of the basis functions $\phi_\ell(\xx)$ is needed. Due to the unstructured nature of the polygonal mesh, where, in principle, no reference element can be adopted (nor the unit triangle neither the unit square), a common solution consists in employing modal basis functions (like Taylor basis). Recently in \cite{ADERAFEDG}, an alternative ansatz has been forwarded for defining the basis functions, using local agglomerated finite element spaces. In the sequel, we present a novel strategy to approximate the discrete solution fitting the general form \eqref{eqn.uh}.

\subsection{The local Virtual Element space}
We define the \emph{local} Virtual Element space of degree $N$, $\Vh_N(P)$, on each Voronoi element $P \in \mathcal{T}_\Omega$ to be
\begin{equation}
\label{eq:Vh:def}
	\Vh_N(P) = \Big\{\,
	\vsh\in\HONE(P)\,:\, \Delta\vsh \in\PS{N-2}(P) \,, \quad
	\vsh|_{\partial P}\in\CS{0}(\partial P)\,, \quad
	\vsh|_{\E}\in\PS{N}(\E)\,\,\forall\E\in\partial P \, \Big\}\,,
\end{equation}
where $\vsh \in \Vh_N(P)$ is a function of the local Virtual Element space and $\PS{N}(P)$ denotes the space of polynomials of degree up to $N$ on $\P$. Hence, $\PS{N}(P)$ is a subspace of $\Vh_N(P)$ with dimension 
\begin{equation}
\label{eq:nk} 
	n_N := \textrm{dim} (\PS{N}(P)) = \frac{(N+1)(N+2)}{2}.
\end{equation} 

First, let us now define a basis for $\PS{N}(\P)$. We introduce the multi-index $\bm{\ell}=(\ell_1,\ell_2)$, thus, if $\bm{x}=(x_1,x_2)$, then $\bm{x}^{\bm{\ell}}=(x_1^{\ell_1},x_2^{\ell_2})$. We define the scaled monomials $m_{\bm{\ell}}$ as
\begin{align}\label{eq:m_alpha} 
	m_{\bm{\ell}} = \Big(\frac{\xv-\xvP}{\hP}\Big)^{\bm{\ell}}, \qquad 0 \leq \ell_1+\ell_2 \leq N.
\end{align}
A basis for $\PS{N}(\P)$ can then be provided by the set $\mathcal{M}_N(\P)$ of scaled monomials of degree less than or equal to $N$:
\begin{align}\label{eq:M_k} 
	\mathcal{M}_N(\P) := \{ m_{\bm{\ell}} : 0 \leq \ell_1+\ell_2 \leq N \}.
\end{align}
We also use the $\alpha$ subscript to denote the $\alpha$-th scaled monomial $m_{\alpha}$ of $\mathcal{M}_N(\P)$.

Next, let us focus on the definition of a canonical basis for the local virtual space. Each function $\vsh$ belonging to $\Vh_N(P)$ is uniquely defined by the following degrees of freedom (dof) \cite{vem2}: 
\begin{itemize}
	\item[i)] the $N_{e}$ point values of $\vsh$ at the vertexes of $P$;
	\item[ii)] the $(N-1)N_{e}$ point values of $\vsh$ on the $(N+1)$ Gauss-Lobatto quadrature points on each edge $e$;  
	\item[iii)] the scaled internal moments up to degree $(N-2)$ of $\vsh$ in $P$, defined as
	\begin{equation*}
		\label{eqn.dof_mom}
		\frac1{\abs{P}} \int_{P} \vsh \,m_{\alpha} \, d\xx, \qquad \alpha = 1,...,n_{N-2}\,,
	\end{equation*}
	where $n_{N-2}$ is the dimension of $\PS{N-2}(P)$ evaluated according to \eqref{eq:nk}.
\end{itemize}
The dimension of $\Vh_N(\P)$ and, therefore, also the number of degrees of freedom of each function $\vsh$ results then
\begin{equation}
\label{eq:dim_Vh} 
	\NDOF := \textrm{dim}(\Vh_N(\P)) =  N \cdot N_e + \frac{N(N-1)}{2}\,.
\end{equation}
We denote by $\hat{\vs}_{k}$ the value of the $k$-th degree of freedom of $\vsh$ and by $\{\varphi\}_{\ell=1}^{\NDOF}$ the local canonical basis for $\Vh_N(\P)$ such that for every $k=1, ..., \NDOF$ it holds
\begin{equation}
	\label{eqn.interp_prop}
	\{\hat{\varphi}_k\}_{\ell=1}^{\NDOF} = \delta_{k \ell},
\end{equation}
with $\delta_{k \ell}$ being the Kronecker delta function. Thus, we can represent each function $\vsh$ of the \emph{local} Virtual Element space in terms of its degrees of freedom using the following Lagrange interpolation:
\begin{equation}
\label{eq:vh_repr} 
	\vsh = \sum_{\ell=1}^{\NDOF} \hat{\vs}_{\ell} \,\varphi_\ell\,.
\end{equation}
Let us notice that the basis functions $\{\varphi_\ell\}_{\ell=1}^{\NDOF}$ are \emph{virtual}, meaning that they are not explicitly known but they are ensured to span the local Virtual Element space defined by \eqref{eq:Vh:def}. Only the degrees of freedom $\{\hat{\varphi}_k\}_{\ell=1}^{\NDOF}$ of the virtual basis are known.

\subsubsection{Elliptic projection operator}
Since the explicit definition of the virtual basis functions is not available, we need to design an appropriate projection operator $\proj: \Vh_N(\P) \rightarrow \PS{N}(\P)$, which maps functions from the Virtual Element space to the known polynomial space of degree $N$. Furthermore, such operator must be computable relying only on the \emph{known} degrees of freedom of the virtual basis. According to \cite{vem2}, the projector is defined from the orthogonality condition
\begin{equation}
\label{eq:vem_proj} 
	\int_{\P} \nabla p \cdot \nabla (\proj\vsh - \vsh) \, d\xx = 0\,, \qquad \forall p \in \PS{N}(\P)\,,\quad \forall \vsh \in \Vh_N(\P)\, ,
\end{equation}
where $p$ is a function of the polynomial space $\PS{N}(\P)$. The above relation is uniquely determined by prescribing an additional projection operator that acts onto constants, namely $P_0: \Vh_N(\P) \rightarrow \PS{0}(\P)$, fulfilling the following relation:
\begin{equation}
\label{eq:vem_proj_const} 
	P_0(\proj\vsh - \vsh) = 0\,.
\end{equation}
Following \cite{vem2}, we choose $P_0$ so that
\begin{align}
\label{eq:vem_proj_0} 
	P_0 \vsh = \begin{cases}
				\frac1{N_e} \sum_{k=1}^{N_e} \vsh(\xv_{r_k}) \quad &\textrm{for} \enskip N = 1\,, \\
				\frac1{\abs{\P}}\int_{\P} \vsh \, d\xx \quad &$\textrm{otherwise}$\,,
				\end{cases}
\end{align}
with $\xv_{r_k}$ being the coordinates of all vertexes of cell $\P$. The basis $\mathcal{M}_N(\P)$ given by \eqref{eq:M_k} is a basis of $\PS{N}(P)$ expressed in terms of the scaled monomials \eqref{eq:m_alpha}, thus we can rewrite the orthogonality condition \eqref{eq:vem_proj} as
\begin{equation}
\label{eq:vem_proj2} 
	\int_{\P} \nabla m_\alpha \cdot \nabla (\proj\vsh - \vsh) \, d\xx = 0\,, \qquad \alpha = 1,...,n_N\,.
\end{equation}
Moreover, being also $\proj\vsh$ an element of $\PS{N}(\P)$, we can use the basis $m_\alpha$ to compute it as well, that is
\begin{equation}
\label{eq:vem_proj2b} 
	\proj\vsh = \sum_{\beta=1}^{n_N} s^\beta m_\beta\,,
\end{equation}
where $s^\beta$ are the unknowns which define the projector operator. The above equation inserted into \eqref{eq:vem_proj2} and \eqref{eq:vem_proj_const} finally leads to the following linear system:
\begin{align}
	\label{eq:vem_proj3} 
	&\sum_{\beta=1}^{n_N} s^\beta \int_{\P} \nabla m_\alpha \cdot \nabla m_\beta \, d\xx = \int_{\P} \nabla m_\alpha \cdot \nabla \vsh \, d\xx \,, \qquad \alpha = 1,...,n_N\,,\\
	\label{eq:vem_proj_aux} 
	&\sum_{\beta=1}^{n_N} s^\beta P_0 m_\beta = P_0 \vsh \,.
\end{align}
The terms on the left side of \eqref{eq:vem_proj3} can be straightforwardly computed since they involve the integration of known polynomials over $\P$, while integration by parts is performed on the right hand side yielding
\begin{equation}
\label{eq:vem_proj4} 
	\int_{\P} \nabla m_\alpha \cdot \nabla \vsh \, d\xx = -\int_{\P}\Delta m_\alpha \vsh \, d\xx + \int_{\partial \P} \frac{\partial m_\alpha}{\partial n} \vsh \, dS.
\end{equation}
We notice that the first term can be directly computed at the sole aid of the internal degrees of freedom of $\vsh$, while the second one (which contains the derivative of the basis $m_\alpha$ in normal direction $n$ with respect to the cell boundary $\partial \P$) can be exactly computed through Gauss-Lobatto quadrature rules along each edge of the element. The values of $\vsh$ at these points are indeed known, and they coincide with some degrees of freedom of the local Virtual Element space. 
Once the solution vector $s^\beta$ of system \eqref{eq:vem_proj3}-\eqref{eq:vem_proj_aux} is obtained, it is employed to evaluate the projection $\proj\varphi_k$ of the virtual basis function $\vsh = \varphi_\ell$ recalling equation \eqref{eq:vem_proj2b}:
\begin{equation*} 
	\proj\varphi_\ell = \sum_{\beta=1}^{n_N} s^\beta_\ell m_\beta, \qquad \ell = 1,...,\NDOF\,.
\end{equation*}

System \eqref{eq:vem_proj3}-\eqref{eq:vem_proj_aux} can be expressed in matrix form as
\begin{equation}
\label{eq:proj_system} 
	\mathbf{G} \accentset{\star}{\mathbf{\Pi}}^{\nabla}_N = \mathbf{B}\,,
\end{equation}
where $\mathbf{G}$ and $\mathbf{B}$ are the following $n_N \times n_N$ and $n_N \times \NDOF$ matrices, respectively:
\begin{align}
\label{gb_matrices} 
	(\mathbf{G})_{\alpha \beta} &= P_0 m_\beta\,, \quad &&\textrm{for} \enskip \alpha = 1, \enskip \beta = 1,...,n_N\,, \\
	(\mathbf{G})_{\alpha \beta} &= \int_{\P}\nabla m_\alpha \cdot \nabla m_\beta \, d\xx\,, \quad &&\textrm{for} \enskip  \alpha \geq 2, \enskip \beta = 1,...,n_N\,,\\
	(\mathbf{B})_{\alpha  \ell} &= P_0 \varphi_\ell\,, \quad &&\textrm{for} \enskip \alpha = 1, \enskip \ell = 1,...,\NDOF\,, \\
	(\mathbf{B})_{\alpha  \ell} &= \int_{\P} \nabla m_\alpha \cdot \nabla \varphi_\ell \, d\xx\,, \quad &&\textrm{for} \enskip \alpha \geq 2, \enskip \ell = 1,...,\NDOF\,,
\end{align}
while $\accentset{\star}{\mathbf{\Pi}}^{\nabla}_N$ is the $n_N \times \NDOF$ matrix representation of the projection operator $\proj$ in the basis set $\mathcal{M}_N(\P)$, i.e., $(\accentset{\star}{\mathbf{\Pi}}^{\nabla}_N)_{\alpha \ell} = s_\ell^\alpha$. 

\subsubsection{$L_2$ projection operator} \label{ssec.L2proj}
Similarly to what done in the previous section, we introduce here the $L_2$ projection operator $\projL: \Vh_N(\P) \rightarrow \PS{N}(\P)$, again mapping functions from the virtual to the polynomial space of degree $N$. This projector arises from the following orthogonality condition:
\begin{equation}
\label{eq:vem_projL1} 
	\int_{\P} p (\projL\vsh - \vsh) \, d\xx = 0\,, \qquad \forall p \in \PS{N}(\P)\,,\quad \forall \vsh \in \Vh_N(\P)\, .
\end{equation}
Once more, it is possible to compute the projector operator only using the known degrees of freedom of $\vsh$. Indeed, since $\projL\vsh$ is an element of $\PS{N}(\P)$, we have that
\begin{equation}
\label{eq:vem_projL2} 
	\projL\vsh = \sum_{\beta=1}^{n_N} r^\beta m_\beta\,,
\end{equation}
which allows us to reformulate the condition \eqref{eq:vem_projL1} as
\begin{equation}
\label{eq:vem_projL3} 
	\sum_{\beta=1}^{n_N} r^\beta \int_{\P} m_\alpha m_\beta \, d\xx = \int_{\P} m_\alpha \vsh \, d\xx\,, \qquad \alpha = 1,...,n_N\,.
\end{equation}
Also in this case, the term on the left side of the above equation can be directly evaluated from the monomials \eqref{eq:m_alpha}, while the integral on the right side is still not computable. To solve the problem, we approximate the function $\vsh$ with its elliptic projection $\proj \vsh$ only for monomials of degree $N$ and $(N-1)$, since the moments for $m_\alpha \in \PS{N-2}(\P)$ are readily known as degrees of freedom of the virtual basis. For $\vsh = \varphi_\ell$, the resulting system can be written in the following matrix notation:
\begin{equation}
\label{eq:proj_system2} 
	\mathbf{H} \accentset{\star}{\mathbf{\Pi}}^0_N = \mathbf{C}\,,
\end{equation}
where $\mathbf{H}$ and $\mathbf{C}$ are an $n_N \times n_N$ and $n_N \times \NDOF$ matrix, respectively, defined as
\begin{align} 
\label{h_matrix}
	(\mathbf{H})_{\alpha \beta} &=\int_{\P}m_\alpha m_\beta \, d\xx\,, \quad\qquad\qquad \alpha,\beta = 1,...,n_N\,,  \\
\label{c_matrix}
	(\mathbf{C})_{\alpha  \ell} &= 
	\begin{cases}
		\int_{\P}m_\alpha \varphi_\ell \, d\xx\,, \qquad &\alpha= 1,...,n_{N-2}, \enskip \ell = 1,...,\NDOF\,, \\
		\int_{\P}m_\alpha \proj \varphi_\ell \, d\xx\,, \qquad &n_{N-2}+1 \leq \alpha \leq n_N, \enskip \ell = 1,...,\NDOF\,. 
	\end{cases} 
\end{align}
The matrix representation of the $L_2$ projection operator $\proj$ is given by $\accentset{\star}{\mathbf{\Pi}}^0_N$ with dimension $n_N \times \NDOF$.

\subsection{Nonconforming Virtual Element basis functions}
\label{ssec.basis}
The local Virtual Element space can be described by the elliptic and $L_2$ projection operators previously introduced (see definitions \eqref{eq:proj_system} and \eqref{eq:proj_system2}, respectively). The Virtual Element Method \cite{vem2} is then built upon the \emph{global} Virtual Element space that is constructed by gluing together all the local spaces $\Vh_N(P_i)$ for $i=1,...,N_P$. In this sense, the VEM methodology belongs to the class of \emph{continuous} finite element solvers (FEM). Differently, the nonconforming Virtual Element Method, originally presented in \cite{VEM_nc_org}, allows the solution to be discontinuous at element vertexes by defining nonconforming local virtual basis along the edges of each cell, meaning that continuity along the boundary $\partial \P$ is no longer guaranteed. A global Virtual Element space is finally still defined, formally collecting all the local virtual spaces, hence solving again the problem at the \emph{global} level.

Here, we aim at designing a new ansatz for the numerical representation of the solution within each control volume $P$, thus we are interested in the definition of local spaces, and corresponding basis functions, in the framework of discontinuous Galerkin methods. In practice, we need to provide an explicit formulation for the basis functions $\phi_\ell(\xx)$ appearing in the ansatz \eqref{eqn.uh}. To this end, the idea is to use the \emph{local} Virtual Element space in \eqref{eqn.uh} to approximate the numerical solution in the polynomial space $\PS{N}$, bearing in mind that we can actually employ the projectors previously described. Specifically, the basis functions $\phi_\ell(\xx)$ in \eqref{eqn.uh} are then \emph{approximated} using the $L_2$ projector \eqref{eq:proj_system2} as
\begin{equation}
	\label{eqn.vem_basis}
	\phi_\ell(\xx) := \varphi_\ell(\xx) \simeq \left( \accentset{\star}{\mathbf{\Pi}}^0_N\right)_{\alpha,\ell} \cdot m_{\alpha}, \qquad \ell=1,...,\NDOF,
\end{equation}    
where $\NDOF$ is the total number of degrees of freedom, hence implying $\mathcal{N}=\NDOF$ in \eqref{eqn.uh}. The \emph{spatial} discrete approximation of the generic quantity $q(\xx,t)$ writes then
\begin{equation}
	q_h^n:=q\left(\xx |_{P},t^n \right) = \sum \limits_{\ell=1}^{\NDOF} \left( \accentset{\star}{\mathbf{\Pi}}^0_N\right)_{\alpha,\ell} \cdot m_{\alpha} \, \hat{q}^{n}_{\ell}, \qquad \xx \in {P}.
	\label{eqn.uhvem}
\end{equation}
We refer to this ansatz as \emph{nonconforming Virtual Element basis}, since we admit jumps in the definition of the global solution space along the entire edges of the computational elements, like in classical Godunov-type DG or finite volume solvers. The numerical solution remains locally defined within each cell, and it spans the local Virtual Element space \eqref{eq:Vh:def}. 

In order to deal with a monolithic space-time discretization, we also need to construct a consistent space-time basis. This is simply achieved by a tensor product of the nonconforming Virtual Element basis in space \eqref{eqn.vem_basis} with one-dimensional nodal basis functions in time. In particular, the time nodal basis consists of $(N+1)$ linearly independent Lagrange interpolating polynomials belonging to the space $\PS{N}$, i.e. $\left\{\psi_k\right\}_{k=1}^{N+1}$, passing through a set of $(N+1)$ nodal points $\left\{\tau_k\right\}_{k=1}^{N+1}$, which are assumed to be the Gauss-Legendre nodes \cite{stroud} referred to the unit time interval given by the mapping \eqref{eqn.dt}. The resulting space-time basis functions $\theta(\xx,\tau)$ are expressed as
\begin{equation}
	\label{eqn.vem_basis_st}
	\theta_r(\xx,\tau) := \varphi_\ell(\xx) \cdot \psi_k(\tau) \simeq \left( \left( \accentset{\star}{\mathbf{\Pi}}^0_N\right)_{\alpha,\ell} \cdot m_{\alpha} \right) \cdot \psi_k, \qquad \ell=1,...,\NDOF, \quad k=1,\,...,N+1, \quad r=1,...,\NDOFST,
\end{equation}
where the total number of space-time degrees of freedom is $\NDOFST=\NDOF \cdot (N+1)$. Consequently, the numerical representation of the generic \emph{space-time} quantity $q(\xx,t)$ explicitly reads
\begin{equation}
	q_h:=q\left(\xx |_{P},t \in \dt \right) = \sum \limits_{\ell=1}^{\NDOFST} \theta_\ell(\xx,\tau) \, \hat{q}_{\ell} , \qquad \xx \in {P}, \quad \tau \in [0;1].
	\label{eqn.qhvem}
\end{equation} 

We underline that both spatial \eqref{eqn.vem_basis} and space-time \eqref{eqn.vem_basis_st} basis functions are actually computed at the aid of the $L_2$ projector detailed in Section \ref{ssec.L2proj}, hence implying an \emph{approximation} of the virtual basis functions $\{\varphi_\ell\}_{\ell=1}^{\NDOF}$ that remain unknown.

\begin{remark}[Orthogonalization of the Virtual Element basis]
If the shape of the polygonal cell is very irregular, meaning that the mesh regularity assumptions defined in Section \ref{ssec.domain} are barely satisfied, the construction of the projection matrices for the the elliptic and $L_2$ operators might lead to very ill-conditioned linear systems \cite{Berrone2017}. This poses serious limitations to the practical usage of the local Virtual Element space and the associated basis, especially for degree $N>2$, as analyzed in \cite{Mascotto2018}. To overcome this problem, in \cite{Berrone2017} an orthogonalization of the projectors is proposed, showing a remarkable improvement in the condition number of the algebraic systems related to the construction of the VEM projectors, namely \eqref{eq:proj_system} and \eqref{eq:proj_system2}. 

Therefore, for $N>2$ we propose a Gram-Schmidt orthogonalization along the lines of \cite{Berrone2017}. However, the entire basis is orthogonalized up to degree $N$, and not only up to degree $N-1$ as forwarded in \cite{Berrone2017}. To that aim, let us introduce matrix $\mathbf{R}$ such that
\begin{equation}
	\label{eqn.L_matrix}
	\mathbf{R} \, \mathbf{H} \, \mathbf{R}^{\top} = \boldsymbol{\Lambda},
\end{equation}
where matrix $\mathbf{H}$ is given by \eqref{h_matrix} and $\boldsymbol{\Lambda}$ is the diagonal matrix with the eigenvalues of $\mathbf{H}$. It follows that the columns of matrix $\mathbf{R}$ contain the right-eigenvectors of $\mathbf{H}$. The orthogonal matrix $\mathbf{Z}$ is computed as
\begin{equation}
	\label{eqn.T_matrix}
	\mathbf{Z} = \sqrt{\mathbf{\boldsymbol{\Lambda}}^{-1}} \, \mathbf{R}.
\end{equation}
In this way, the set of orthonormal polynomials $\{\mathbf{z_\alpha}\}_{\alpha=1}^{n_N}$ spanning the polynomial space $\PS{N}$ are defined by
\begin{equation}
	z_{\alpha} = (\mathbf{Z})_{\alpha \beta} \, m_{\beta}, \qquad \textrm{for} \enskip \alpha = 1,...,n_N, \enskip \beta = 1,...,n_N,
\end{equation}
where $m_{\beta}$ are the monomials introduced in \eqref{eq:m_alpha}. It is easy to verify that the orthonormal polynomials yield an identity mass matrix $\tilde{\mathbf{H}}$ of dimension $n_N \times n_N$, indeed
\begin{eqnarray}
	\label{eqn.mass_matrix_orth}
	(\tilde{\mathbf{H}})_{\alpha \beta} &=& \int_{\P}z_\alpha z_\beta \, d\xx \nonumber \\
	&=& \int_{\P}(\mathbf{Z})_{\alpha k} m_k \,  (\mathbf{Z})_{\beta \ell} m_\ell \, d\xx \nonumber \\
	&=& \int_{\P}(\mathbf{Z})_{\alpha k} m_k \,  m_\ell (\mathbf{Z})_{\beta \ell} \, d\xx \nonumber \\
	&=& \int_{\P}(\mathbf{Z})_{\alpha k} \, (\mathbf{H})_{k \ell} \, (\mathbf{Z})_{\beta \ell} \, d\xx \nonumber \\
	&=& \int_{\P}\sqrt{\mathbf{(\Lambda)}_{\alpha \alpha}^{-1}} \, (\mathbf{R})_{\alpha k} \, (\mathbf{H})_{k \ell} \, \sqrt{\mathbf{(\Lambda)}_{\beta \beta}^{-1}} \, (\mathbf{R})_{\beta \ell} \, d\xx \nonumber \\
	&=& \int_{\P}\sqrt{\mathbf{(\Lambda)}_{\alpha \alpha}^{-1}} \, (\mathbf{R})_{\alpha k} \, (\mathbf{H})_{k \ell} \, (\mathbf{R})_{\ell \beta} \, \sqrt{\mathbf{(\Lambda)}_{\beta \beta}^{-1}} \, d\xx \nonumber \\
	&=& \int_{\P}\sqrt{\mathbf{(\Lambda)}_{\alpha \alpha}^{-1}} \, \mathbf{(\Lambda)}_{\alpha \beta} \, \sqrt{\mathbf{(\Lambda)}_{\beta \beta}^{-1}} \, d\xx = \delta_{\alpha \beta}.
\end{eqnarray}
The nonconforming Virtual Element basis can thus be orthogonalized by modifying the definition \eqref{eqn.vem_basis} as follows:
\begin{equation}
	\label{eqn.vem_basis_orth}
	\phi_\ell(\xx) = \left( \accentset{\star}{\mathbf{\Pi}}^0_N\right)_{\alpha,\ell} \cdot \left( (\mathbf{Z})_{\alpha \beta} \, m_{\beta} \right), \qquad \ell=1,...,\NDOF, \quad \alpha,\beta=1,...,n_N.
\end{equation}
The space-time basis functions are always constructed as a tensor product of the orthogonalized spatial basis \eqref{eqn.vem_basis_orth} and the one-dimensional Lagrange nodal basis $\left\{\psi_k\right\}_{k=1}^{N+1}$ according to \eqref{eqn.vem_basis_st}.
\end{remark}

\begin{remark}[Derivatives of the Virtual Element basis] The discrete spatial derivatives of the nonconforming Virtual Element basis \eqref{eqn.vem_basis} are evaluated as follows:
\begin{equation}
	\label{eqn.vem_basis_der}
	\frac{\partial \phi_\ell}{\partial x} = \left( \accentset{\star}{\mathbf{\Pi}}^{0,x}_{N-1}\right)_{\alpha,\ell} \cdot m_{\alpha}, \quad \frac{\partial \phi_\ell}{\partial y} = \left( \accentset{\star}{\mathbf{\Pi}}^{0,y}_{N-1}\right)_{\alpha,\ell} \cdot m_{\alpha}, \qquad \ell=1,...,\NDOF, \quad \alpha=1, ..., n_{N-1}.
\end{equation} 
The $L_2$ operators $\accentset{\star}{\mathbf{\Pi}}^{0,x}_{N-1}$ and $\accentset{\star}{\mathbf{\Pi}}^{0,y}_{N-1}$ project the derivative of the solution in $x$ and $y$ direction from the local Virtual Element space $\Vh_{N-1}(\P)$ to the polynomial space $\PS{N-1}$, thus they are the matrix representation of the projectors $(\Pi^{0,x}_{\P,N-1},\Pi^{0,y}_{\P,N-1}): \Vh_{N-1}(\P) \rightarrow \PS{N-1}$. They are constructed following the same methodology used for the $L_2$ projector described in Section \ref{ssec.L2proj}. However, they refer to a polynomial degree $N-1$ and they project the derivatives of the virtual basis functions. In particular, the derivative operators introduced in \eqref{eqn.vem_basis_der} are computed as
\begin{align}
	\accentset{\star}{\mathbf{\Pi}}^{0,x}_{N-1} &= \hat{\mathbf{H}}^{-1}\mathbf{E}^x, \\
	\accentset{\star}{\mathbf{\Pi}}^{0,y}_{N-1} &= \hat{\mathbf{H}}^{-1}\mathbf{E}^y, 
\end{align}
where matrix $\hat{\mathbf{H}}$ is obtained by taking the first $n_{N-1}$ rows and columns of matrix $\mathbf{H}$ defined in \eqref{h_matrix}. Matrices $\mathbf{E}^x$ and $\mathbf{E}^y$ account for the derivatives of the virtual basis and they write
\begin{align}
	(\mathbf{E}^x)_{k\alpha} = \int_{\P} \varphi_{k,x}m_\alpha \, d\xx, \quad\quad	(\mathbf{E}^y)_{k\alpha} = \int_{\P} \varphi_{k,y}m_\alpha \, d\xx, \quad \quad \alpha = 1,...,n_{N-1}.
\end{align}
In practice, they can be computed performing integration by parts:
\begin{equation}
	\label{eqn.Ex}
	\int_{\P} \varphi_{k,x}m_\alpha \, d\xx = -\int_{\P} \varphi_{k} m_{\alpha,x} \, d\xx + \int_{\partial \P} \varphi_{k} m_\alpha n_x \, dS, \quad \quad \alpha = 1,...,n_{N-1}. 
\end{equation}	
The first integral on the right hand side of \eqref{eqn.Ex} can be computed relying on the internal degrees of freedom of the local Virtual Element space (scaled internal momenta which are available up to degree $(N-2)$), while the second integral can also be readily evaluated using the degrees of freedom lying on the boundary $\partial \P$, which coincide with the Gauss-Lobatto quadrature nodes along each edge of the cell. The same holds true for the matrix $\mathbf{E}^y$.

The time derivative of the space-time Virtual Element basis \eqref{eqn.vem_basis_st} does not pose any problem, and it can be simply evaluated by taking the derivative of the Lagrange basis in time. Therefore the space-time derivatives are computed at the aid of the definitions \eqref{eqn.vem_basis_der} as follows:
\begin{equation}
	\frac{\partial \theta_r}{\partial x} = \frac{\partial \phi_\ell}{\partial x} \cdot \psi_k, \quad \frac{\partial \theta_r}{\partial y} = \frac{\partial \phi_\ell}{\partial y} \cdot \psi_k, \quad \frac{\partial \theta_r}{\partial \tau} = \phi_\ell \cdot \frac{\partial \psi_k}{\partial \tau}.
\end{equation}
\end{remark}

%

\subsection{Space-time VEM-DG scheme}
We discretize the Navier--Stokes equations defined by system \eqref{eqn.PDE}-\eqref{eqn.NS} designing an ADER VEM-DG scheme, which permits to reach arbitrary accuracy in space and time. This numerical scheme consists in two main steps. 
\begin{enumerate}
\item The computation of a local space-time ADER \emph{predictor} starting from the known solution at the current time step $t^n$ for each mesh element $\P$. This step gives a high order space-time approximation of the discrete solution that is valid locally inside each cell within the time interval $[t^n;t^{n+1}]$. There is no need for information exchanges between neighbor cells.
\item A \emph{corrector} step that leads to the effective solution at the new time step $t^{n+1}$ by directly integrating the system of equations over the space-time control volume $\P \times [t^n;t^{n+1}]$. The spatial discretization is built upon the novel nonconforming Virtual Element basis used to devise a discontinuous Galerkin scheme. Thanks to the usage of Riemann solvers at the cell interfaces, the corrector step will also take into account the interaction between the predictors of neighboring cells.
\end{enumerate}

\subsubsection{ADER local space-time predictor} \label{ssec.ader}
The ADER methodology was firstly introduced in \cite{toro3,titarevtoro} aiming at computing a predictor solution $\q_h(\xx,t)$ by solving the generalized Riemann problem ``in the small", hence neglecting interactions between neighbor cells. To evolve the solution in time, time derivatives must be known. They can be evaluated relying on an element-local \textit{weak formulation} of the governing equations according to \cite{Dumbser2008}. 

We start by multiplying the PDE \eqref{eqn.PDE} by a set of test functions $\theta_k(\xx,\tau)$ of the same form of the space-time basis \eqref{eqn.vem_basis_st}, and then we integrate the governing system over a space-time control volume given by $\P \times [0;1]$, obtaining
\begin{equation}
	\label{eqn.ader1}
	\int_{0}^{1} \int_{\P} \theta_k \, \frac{\partial \Q}{\partial \tau}  \, d\xx \, d\tau = - \dt \, \int_{0}^{1} \int_{\P} \theta_k \, \left( \frac{\partial \f}{\partial x} + \frac{\partial \g}{\partial y} \right) \, d\xx \, d\tau, \qquad k = 1, ..., \NDOFST. 
\end{equation}
Let us notice that the PDE has been reformulated in the \emph{reference time coordinate} $\tau$, hence implying the presence of the Jacobian $1/\dt$ of the transformation according to the mapping \eqref{eqn.dt}. The unknown discrete solution $\q_h$ as well as the fluxes $(\f_h,\g_h)$ are approximated using the nonconforming Virtual Element space-time basis \eqref{eqn.vem_basis_st}, thus
\begin{equation}
	\label{eqn.qh}
	\q_h(\xx,t) = \sum \limits_{\ell=1}^{\NDOFST} \theta_\ell(\xx,\tau) \, \hat{\q}_{\ell}, \quad \f_h(\xx,t) = \sum \limits_{\ell=1}^{\NDOFST} \theta_\ell(\xx,\tau) \, \hat{\f}_{\ell}, \quad \g_h(\xx,t) = \sum \limits_{\ell=1}^{\NDOFST} \theta_\ell(\xx,\tau) \, \hat{\g}_{\ell}.
\end{equation}
Thanks to the interpolation property \eqref{eqn.interp_prop}, it holds that 
\begin{equation}
	\label{eqn.f_interp}
	\hat{\f}_\ell=\f(\hat{\q}_\ell) \quad \textrm{and} \quad \hat{\g}_\ell=\g(\hat{\q}_\ell),
\end{equation}
so the degrees of freedom of the fluxes can be simply evaluated \emph{pointwise} from $\q_h$. The above definitions are inserted into the weak formulation \eqref{eqn.ader1}, hence obtaining
\begin{equation}
	\label{eqn.ader2}
	\int_{0}^{1} \int_{\P} \theta_k \, \frac{\partial \theta_\ell}{\partial \tau} \hat{\q}_\ell  \, d\xx \, d\tau = - \dt \, \int_{0}^{1} \int_{\P} \theta_k \, \left( \frac{\partial \theta_\ell}{\partial x} \hat{\f}_\ell + \frac{\partial \theta_\ell}{\partial y} \hat{\g}_\ell \right) \, d\xx \, d\tau, \qquad k = 1, ..., \NDOFST.
\end{equation}
To satisfy the causality principle accounting only for information coming from the past in each cell, the term on the left hand side of \eqref{eqn.ader2} is integrated by parts in time, yielding
\begin{eqnarray}
	\label{eqn.ader3}
	&& \left( \int_{\P} \theta_k(\xx,1) \theta_\ell(\xx,1) \, d\xx - \int_{0}^{1} \int_{\P} \frac{\partial \theta_k}{\partial \tau} \, \theta_\ell \, d\xx \, d\tau \right) \hat{\q}_\ell  =\nonumber \\
	&& \int_{\P} \theta_k(\xx,0) \phi_\ell(\xx,0) \, \hat{\uu}^n_\ell \, d\xx -
	 \dt \, \int_{0}^{1} \int_{\P} \theta_k \, \left( \frac{\partial \theta_\ell}{\partial x} \hat{\f}_\ell + \frac{\partial \theta_\ell}{\partial y} \hat{\g}_\ell \right) \, d\xx \, d\tau, \qquad k = 1, ..., \NDOFST,
\end{eqnarray}
where the numerical solution at the current time $t^n$ (i.e. $\tau=0$) is expressed with its definition given by \eqref{eqn.uh} in terms of the spatial basis $\phi_\ell$ defined in \eqref{eqn.vem_basis}. We can rewrite the weak formulation \eqref{eqn.ader3} in compact matrix-vector notation as
\begin{equation}
	\label{eqn.ader4}
	\mathbf{K}_1 \, \hat{\q}_\ell = \mathbf{F}_0 \, \uu^n_\ell - \dt \left( \mathbf{K}_x\, \f(\hat{\q}_\ell) + \mathbf{K}_y  \,\g(\hat{\q}_\ell)\right) ,
\end{equation} 
with the definitions
\begin{eqnarray}
	\label{eqn.st_matrices}
	&& \mathbf{K}_1 = \int_{\P} \theta_k(\xx,1) \theta_\ell(\xx,1) \, d\xx - \int_{0}^{1} \int_{\P} \frac{\partial \theta_k}{\partial \tau} \, \theta_\ell \, d\xx \, d\tau, \qquad \mathbf{F}_0 = \int_{\P} \theta_k(\xx,0) \phi_\ell(\xx,0) \, d\xx, \nonumber \\
	&& \mathbf{K}_x = \int_{0}^{1} \int_{\P} \theta_k \, \frac{\partial \theta_\ell}{\partial x} \, d\xx \, d\tau, \qquad 	\mathbf{K}_y = \int_{0}^{1} \int_{\P} \theta_k \, \frac{\partial \theta_\ell}{\partial y} \, d\xx \, d\tau.
\end{eqnarray}
The nonlinear algebraic equation system \eqref{eqn.ader4} is solved locally for the unknown space-time expansion coefficients $\hat{\q}_\ell$ with a simple fixed-point iterative scheme:
\begin{equation}
	\label{eqn.ader5}
	\hat{\q}_\ell^{r+1} = (\mathbf{K}_1)^{-1} \left( \mathbf{F}_0 \, \uu^n_\ell - \dt \left( \mathbf{K}_x \, \f(\hat{\q}_\ell^r) + \mathbf{K}_y  \, \g(\hat{\q}_\ell^r) \right) \right), 
\end{equation}
where the superscript $r$ indicates the iteration number. The iteration stops when the residual of \eqref{eqn.ader5} is less than a prescribed tolerance that guarantees precision, typically set to $10^{-12}$. However, it could also be sufficient a smaller number of iterations to achieve at least the formal order of accuracy, as recently investigated in \cite{AdaADER23}.

\begin{remark}[Computation of the space-time stiffness matrices]
The space-time matrices given by \eqref{eqn.st_matrices} involve integration of the nonconforming Virtual Element space-time basis functions \eqref{eqn.vem_basis_st} over the space-time control volume $\P \times [0;1]$. They are directly computed by means of very efficient quadrature rules on general polygonal cells forwarded in \cite{SOMMARIVA2009886,SOMMARIVA2020}. This gives the so-called consistency term of the matrices, which is enough for those matrices which do not need to be inverted, namely $\mathbf{F_0}$, $\mathbf{K_x}$ and $\mathbf{K_y}$. However, the space-time stiffness matrix $\mathbf{K_1}$ has to be inverted in \eqref{eqn.ader5}, and a stabilization term must be added as usually done in the Virtual Element framework \cite{vem2,vem1}. Indeed, mass and stiffness VEM matrices can become highly ill-conditioned, especially when $N>2$ (see \cite{Mascotto2018}).

There are no analysis in the literature for space-time Virtual Element matrices, thus we mimic and modify what has been done in the purely spatial setting \cite{vem2} in order to account for novel space-time VEM operators. Let us write matrix $\mathbf{K_1}$ using the space-time virtual basis functions $\tilde{\varphi}_{ik}$, which are unknown:
\begin{equation}
	\label{eqn.K1_virt}
	(\mathbf{K_1})_{k \ell} = \int_{0}^{1} \int_{\P} \tilde{\varphi}_{ik} \frac{\partial \tilde{\varphi}_{j\ell}}{\partial \tau} \, d\xx d\tau, \qquad \tilde{\varphi}_{ik} = \varphi_i \, \psi_k, \quad \tilde{\varphi}_{j\ell} = \varphi_j \, \psi_\ell,
\end{equation}
where $\psi_k,\psi_\ell$ are the time basis functions described in Section \ref{ssec.basis}. We then introduce the expansions
\begin{eqnarray}
	\label{eqn.expansion_st}
	\tilde{\varphi}_{ik} &=& \projL \varphi_i \, \psi_k + (I-\projL) \varphi_i \, \psi_k \\
	\frac{\partial \tilde{\varphi}_{j\ell}}{\partial \tau} &=& \frac{\partial}{\partial \tau} \left( \projL \varphi_j \, \psi_\ell + (I-\projL) \varphi_j \, \psi_\ell \right),
\end{eqnarray}
with $I$ being the identity matrix of dimension $\NDOF \times \NDOF$. The above expressions are plugged into the definition \eqref{eqn.K1_virt} leading to
\begin{eqnarray}
	\label{eqn.K1_virt2}
	\int_{0}^{1} \int_{\P} \tilde{\varphi}_{ik} \frac{\partial \tilde{\varphi}_{j\ell}}{\partial \tau} \, d\xx d\tau &=& \int_{0}^{1} \int_{\P} \projL \varphi_i \, \psi_k \cdot \frac{\partial}{\partial \tau} \left( \projL \varphi_j \, \psi_\ell \right) \, d\xx d\tau \nonumber \\
	&+& \int_{0}^{1} \int_{\P} \projL \varphi_i \, \psi_k \cdot \frac{\partial}{\partial \tau} \left( (I-\projL) \varphi_j \, \psi_\ell \right) \, d\xx d\tau \nonumber \\
	&+& \int_{0}^{1} \int_{\P} (I-\projL) \varphi_i \, \psi_k \cdot \frac{\partial}{\partial \tau} \left( \projL \varphi_j \, \psi_\ell \right) \, d\xx d\tau \nonumber \\
	&+& \int_{0}^{1} \int_{\P} (I-\projL) \varphi_i \, \psi_k \cdot \frac{\partial}{\partial \tau} \left( (I-\projL) \varphi_j \, \psi_\ell \right) \, d\xx d\tau.
\end{eqnarray}
Since the basis functions $\psi_k,\psi_\ell$ only depend on time, and $\varphi_i,\varphi_j$ only depend on space, we can rearrange the above equation as
\begin{eqnarray}
	\label{eqn.K1_virt3}
	\int_{0}^{1} \int_{\P} \tilde{\varphi}_{ik} \frac{\partial \tilde{\varphi}_{j\ell}}{\partial \tau} \, d\xx d\tau &=& \int_{0}^{1} \int_{\P} \projL \varphi_i \, \psi_k \cdot \frac{\partial}{\partial \tau} \left( \projL \varphi_j \, \psi_\ell \right) \, d\xx d\tau \nonumber \\
	&+& \int_{0}^{1} \psi_k \cdot \frac{\partial \psi_\ell}{\partial \tau} \, d\tau \, \int_{\P} \projL \varphi_i \cdot (I-\projL) \varphi_j \, d\xx \nonumber \\
	&+& \int_{0}^{1} \psi_k \cdot \frac{\partial \psi_\ell}{\partial \tau} \, d\tau \, \int_{\P} (I-\projL) \varphi_i \cdot \projL \varphi_j \, d\xx \nonumber \\
	&+& \int_{0}^{1} \psi_k \cdot \frac{\partial \psi_\ell}{\partial \tau} \, d\tau \, \int_{\P} (I-\projL) \varphi_i \cdot (I-\projL) \varphi_j \, d\xx.
\end{eqnarray}
The first term on the right hand side ensures consistency, while the second and the third term vanish because of the orthogonality condition \eqref{eq:vem_projL1} related to the construction of the projector $\projL$. The last term ensures stability and the spatial integral is approximated using the dof--dof stabilization \cite{vem2}. Thus, recalling the definition of the space-time Virtual Element basis \eqref{eqn.vem_basis_st}, the space-time stiffness matrix $\mathbf{K}_1$ can eventually be computed by
\begin{eqnarray}
	\label{eqn.K1_final}
	\int_{0}^{1} \int_{\P} \tilde{\varphi}_{ik} \frac{\partial \tilde{\varphi}_{j\ell}}{\partial \tau} \, d\xx d\tau &\simeq& \int_{0}^{1} \int_{\P} \left( \left( \accentset{\star}{\mathbf{\Pi}}^0_N\right)_{\alpha,i} m_{\alpha} \right) \cdot \psi_k \cdot \frac{\partial}{\partial \tau} \left( \left( \left( \accentset{\star}{\mathbf{\Pi}}^0_N\right)_{\alpha,j} m_{\alpha} \right) \cdot \psi_\ell \right) \, d\xx d\tau \nonumber \\
	&+& \int_{0}^{1} \psi_k \cdot \frac{\partial \psi_\ell}{\partial \tau} \, d\tau \, \, \left((\mathbf{I}-\accentset{\star}{\mathbf{\Pi}}^0_N)_{\alpha,i}\right)^{\top} \cdot (\mathbf{I}-\accentset{\star}{\mathbf{\Pi}}^0_N)_{\alpha,j}.
\end{eqnarray}
This formulation provides both consistency and stability in a fully space-time setting, so that matrix $\mathbf{K}_1$ can be inverted and used in the iterative scheme \eqref{eqn.ader5}. 
\end{remark}

\subsubsection{Fully discrete one-step VEM-DG corrector}
The predictor solution is used to carry out the corrector step, which accounts for the numerical fluxes between neighbor cells. The corrector is based on a discontinuous Galerkin scheme that is directly applied to the integrated form of the governing equations \eqref{eqn.PDE} in space and time. Since the predictor solution has already been computed within the local space-time control volumes, the corrector involves a one-step time integration, yielding a fully discrete scheme. This improves the computational efficiency in parallel computation over multi-step schemes, like Runge-Kutta DG methods, because only one communication is needed among the threads, which exchange the predictor solution. Furthermore, being the predictor locally evaluated, no information exchange is needed therein.

The variational formulation is obtained upon multiplication of the PDE \eqref{eqn.PDE} by a spatial test function $\phi_k(\xx)$ of the same form of the basis functions $\phi_\ell(\xx)$ used to approximate the numerical solution \eqref{eqn.uh}, followed by integration on the space-time control volume $\P \times [t^n;t^{n+1}]$ for each $\P_i$, $i=1,...,\NP$:
\begin{equation}
	\label{eqn.pdeweak}
\int_{t^n}^{t^{n+1}} \int_{\P} \phi_k \,\partial_t \Q \,d\xx dt + \int_{t^n}^{t^{n+1}} \int_{\P} \phi_k \,\nabla \cdot \F \,d\xx dt = \mathbf{0}\,.
\end{equation}
Integration by parts in space of the second term leads to
\begin{equation}
	\label{eqn.pdeweak_int_parts}
\int_{t^n}^{t^{n+1}} \int_{\P} \vbase_k \,\partial_t \Q \,d\xx dt + \int_{t^n}^{t^{n+1}} \int_{\partial \P} \vbase_k \,\F\cdot \nor \,dS dt - \int_{t^n}^{t^{n+1}} \int_{\P} \nabla \vbase_k \cdot \F\,d\xx dt = \mathbf{0}\,,
\end{equation}
where $\nor$ is the outward pointing unit normal vector of the element boundary $\partial \P$. Using the solution representation \eqref{eqn.uh} as well as the predictor solution $\q_h(\xx,t)$ and its gradient, the weak form \eqref{eqn.pdeweak_int_parts} becomes
\begin{eqnarray}
	\int_{\P_i} \vbase_k \vbase_\ell \, d\xx \, \uu_{\ell}^{n+1} &=& \int_{\P} \vbase_k \vbase_\ell \, d\xx \, \uu_{\ell}^{n} - \int_{t^n}^{t^{n+1}} \int_{\partial \P} \vbase_k \,\mathcal{G}\left( (\q_h^-,\nabla \q_h^-), (\q_h^+,\nabla \q_h^+) \right) \cdot \nor \,dS dt \nonumber \\
	&+& \int_{t^n}^{t^{n+1}} \int_{\P} \nabla \vbase_k \cdot \F(\q_h,\nabla \q_h) \,d\xx dt, \label{eqn.DGscheme}
\end{eqnarray}
where $\mathcal{G}\left( (\q_h^-,\nabla \q_h^-), (\q_h^+,\nabla \q_h^+) \right) \cdot \nor$ is the numerical flux function, which involves the left $(\q_h^-,\nabla \q_h^-)$ and right $(\q_h^+,\nabla \q_h^+)$ high order boundary-extrapolated data and gradients with respect to the cell boundary $\partial \P$. To compute it, the Rusanov flux~\cite{Rusanov:1961a} is employed, modified to simultaneously include both the convective and the viscous terms:
\begin{equation}
	\mathcal{G}\left( (\q_h^-,\nabla \q_h^-), (\q_h^+,\nabla \q_h^+) \right) \cdot \nor = \frac{1}{2} \left( \F(\q_h^-,\nabla \q_h^-) + \F(\q_h^+,\nabla \q_h^+) \right) \cdot \nor - \frac{1}{2} \left( |\lambda^c_{\max}| + 2 \eta |\lambda^v_{\max}| \right) \left( {\q}_h^+ - {\q}_h^- \right).
	\label{eqn.rusanov}
\end{equation} 
The factor $\eta$ in the numerical dissipation is estimated from the solution of the generalized Riemann problem (GRP) for the diffusion equation \cite{MunzDiffusionFlux} and it is given by 
\begin{equation}
	\eta = \frac{2N+1}{(h_{P^-}+h_{P^+})\sqrt{\frac{\pi}{2}}},
\end{equation}
where $h_{P^-}$ and $h_{P^+}$ are the characteristic length sizes of the two mesh elements (left and right, respectively) that share the interface for which the flux is being evaluated.

The nonconforming Virtual Element basis functions \eqref{eqn.vem_basis} are used to define both $\phi_k$ and $\phi_\ell$, thus the one-step fully discrete DG scheme \eqref{eqn.DGscheme} is referred to as VEM-DG scheme. 

\begin{remark}[Computation of the mass matrix]
The VEM-DG scheme \eqref{eqn.DGscheme} requires the computation of the mass matrix for the generic element $\P$, which is assumed to be defined relying on the unknown virtual basis functions $\varphi(\xx)$, that is
\begin{equation}
	\label{eqn.mass_matrix}
	(\mathbf{M})_{k \ell} := \int_{\P} \vbase_k \vbase_\ell \, d\xx=\int_{\P} \varphi_k \varphi_\ell \, d\xx.
\end{equation}
Following \cite{vem2}, we can introduce the expansion
\begin{equation}
	\varphi_k = \projL \varphi_k  + (I-\projL) \varphi_k,
\end{equation}
that is plugged into the definition \eqref{eqn.mass_matrix}, hence obtaining
\begin{eqnarray}
	\label{eqn.proj_M} 
	(\mathbf{M})_{k \ell} &=&\int_{\P} \projL \varphi_k  \projL \varphi_\ell  \, d\xx + \int_{\P} (I-\projL)\varphi_k  (I-\projL) \varphi_\ell \, d\xx  \nonumber \\ 
	&+& \int_{\P} \projL\varphi_k  (I-\projL) \varphi_\ell \, d\xx + \int_{\P} (I-\projL)\varphi_k  \projL\varphi_\ell \, d\xx.
\end{eqnarray}
The last two terms on the right hand side are exactly zero because of the orthogonality condition \eqref{eq:vem_projL1}, while the second term accounts for the stabilization, and it is approximated again with a dof--dof stabilization \cite{vem2} as for the space-time stiffness matrix in \eqref{eqn.K1_final}, thus yielding
\begin{equation}
	(\mathbf{M})_{k \ell} \simeq \int_{\P} \left( \left( \accentset{\star}{\mathbf{\Pi}}^0_N\right)_{\alpha,k} m_{\alpha} \right)  \, \left( \left( \accentset{\star}{\mathbf{\Pi}}^0_N\right)_{\alpha,\ell} m_{\alpha} \right)  \, d\xx + \abs{\P} \, \left((\mathbf{I}-\accentset{\star}{\mathbf{\Pi}}^0_N)_{\alpha,k}\right)^{\top} \cdot (\mathbf{I}-\accentset{\star}{\mathbf{\Pi}}^0_N)_{\alpha,\ell}.
	\label{eqn.proj_M3}
\end{equation}
The first term in \eqref{eqn.proj_M3} is the consistency term, while the second term is responsible for stabilization, which is crucial for the inversion of the mass matrix in the VEM-DG scheme \eqref{eqn.DGscheme}.
\end{remark}

\subsection{Runge-Kutta VEM-DG scheme}
Alternatively, the nonconforming Virtual Element basis \eqref{eqn.vem_basis} can be used to devise a standard Runge-Kutta scheme, that is a multi-step algorithm based on the method of lines. A Runge-Kutta method with a total number of $s$ sub-stages is described by a Butcher tableau of the form shown in Table \ref{tab.BTRK}.
\begin{table}[h!]
	\caption{Butcher tableau for Runge-Kutta explicit methods.}
	\begin{center} 
		\begin{tabular}{c|ccccc} 
			0 &  & & & & \\
			$\alpha_2$ & $\beta_{21}$ & & & & \\
			$\alpha_3$ & $\beta_{31}$ & $\beta_{32}$ & & & \\
			$\vdots$  & $\vdots$      & $\vdots$  & $\ddots$ & & \\
			$\alpha_s$ & $\beta_{s1}$ & $\beta_{s2}$ & $...$ & $\beta_{s (s-1)}$ & \\
			\hline
			& $c_1$ & $c_2$ & $...$ & $c_{s-1}$ & $c_s$
			\label{tab.BTRK}
		\end{tabular}		
	\end{center}	
\end{table}	

The semi-discrete scheme for the governing PDE is then given by
\begin{equation}
	\mathbf{M} \, \frac{\textrm{d} \Q}{\textrm{d} t} = - \mathcal{L}_h(\uu_h),
\end{equation}
where $\mathbf{M}$ is the Virtual Element spatial mass matrix \eqref{eqn.proj_M3} and the term $\mathcal{L}_h(\uu_h)$ contains the spatial discretization:
\begin{equation}
	\mathcal{L}_h(\uu_h) = \int_{\partial \P} \vbase_k \,\mathcal{G}\left( (\uu_h^-,\nabla \uu_h^-), (\uu_h^+,\nabla \uu_h^+) \right) \cdot \nor \,dS - \int_{\P} \nabla \vbase_k \cdot \F(\uu_h,\nabla \uu_h) \,d\xx.
\end{equation}
Here, the test functions $\phi_k$ are assumed to be of the form of \eqref{eqn.vem_basis}, and the numerical flux term is computed according to the definition \eqref{eqn.rusanov}. The numerical solution is determined at the next time step as
\begin{equation}
	\U^{n+1} = \U^{n} + \mathbf{M}^{-1} \, \dt \, \sum \limits_{i=1}^s c_i \cdot \kappa_i,
\end{equation}
with the generic Runge-Kutta stage $\kappa_i$ evaluated at time level $t^{(i)}=t^n+\alpha_i \dt$ by
\begin{equation}
	\kappa_i = -\mathcal{L}_h\left( \uu_h^n + \dt \sum \limits_{j=1}^{i} \beta_{ij} \cdot \kappa_j \right).
\end{equation}
In this case, the space-time basis functions \eqref{eqn.vem_basis_st} are no longer needed, and the only matrix that must be inverted is the spatial mass matrix $\mathbf{M}$.

\subsection{VEM-DG limiter with artificial viscosity}
The DG scheme \eqref{eqn.DGscheme} is linear in the sense of Godunov \cite{Godunov1959}, thus spurious oscillations might arise when dealing with shock waves and other discontinuities. To limit these instabilities, we rely on a simple artificial viscosity method, see for instance \cite{PerssonAV,TavelliCNS,Hesthaven_LimiterAV2011,Bassi_LimiterAV2018,ADERAFEDG}. The limiter does not act on the entire mesh, but only on those cells which are crossed by discontinuities. Therefore, one first needs to detect the so-called \emph{troubled} cells, and then to limit the numerical solution only there.

At each time step, the troubled elements are detected with the flattener indicator $\beta_{\P}$ proposed in \cite{BalsaraFlattener}, that is computed for each element as
\begin{equation}
	\beta_{\P}= \min {\left[ 1, \max {\left(0, -\frac{\nabla \cdot \vv + \bar{g} c_{\min}}{\bar{g} c_{\min}}\right)_{\P}}\right]},
	\label{eqn.flattener}
\end{equation}
with the coefficient $\bar{g}=0.1$  according to \cite{BalsaraFlattener}. The minimum of the sound speed $c_{\min}$ is evaluated considering the element $\P$ itself and its neighborhood, that is
\begin{equation}
	c_{\min} = \min_{\partial \P} (c^+,c^-), \qquad c^{\pm}=\sqrt{\gamma R T^{\pm}},
\end{equation}
while the divergence of the velocity field $\nabla \cdot \vv$ is estimated as 
\begin{equation}
	\nabla \cdot \vv = \frac{1}{|\P|}\sum \limits_{\partial \P}{ |\partial \P|^{\pm} \left(\vv^+ - \vv^- \right) \cdot \nor },
	\label{eqn.divV}
\end{equation}
where the velocity vector is computed as a cell average quantity from the numerical solution $\uu_h$, and the quantity $|\partial \P|^{\pm}$ is the length of the edge shared by the left $\P^-$ and the right $\P^+$ cell. A cell is marked as troubled if $\beta_{\P}>10^{-10}$, and some artificial viscosity $\mu_{a}$ is added to the physical viscosity $\mu$, thus obtaining an effective viscosity $\tilde{\mu}_{\P}=\mu_{a}+\mu$ which is then used in the Navier-Stokes fluxes \eqref{eqn.NS}. The additional viscosity
$\mu_{a}$ is determined so that a resulting unity mesh Reynolds number is assigned to the troubled cells, which is
\begin{equation}
	\Rey = \frac{\rho |\lambda_{\max}^c| \hP}{\tilde{\mu}_{\P}}=1.
\end{equation}
To account also for heat conduction, a unit Prandtl number must be set, hence imposing an artificial heat conduction coefficient $\tilde{\kappa}_{\P}$ such that $\Pr = \tilde{\mu}_{\P} \gamma c_v/\tilde{\kappa}_{\P} = 1$.

\section{Numerical results} \label{sec.validation}
The new numerical schemes are tested against benchmarks for compressible gas dynamics to properly assess their accuracy and robustness. The label VEM-DG is used to refer to the methods presented in this work, with ADER time discretization unless otherwise indicated. Moreover, if not stated otherwise, we set the ratio of specific heats to $\gamma=1.4$ and the gas constant to $R=1$, hence retrieving a specific heat capacity at constant volume of $c_v = 2.5$. The initial condition of the flow field is typically given in terms of the vector of primitive variables $\Pv(\xx,t) = (\rho,u,v,p)$. The simulations are run on 64 CPUs with MPI parallelization.

\subsection{Numerical convergence study}
We consider the smooth isentropic vortex test case proposed in \cite{HuShuTri} to study the numerical convergence of the VEM-DG schemes. The computational domain is given by $\Omega=[0;10] \times [0;10]$ with periodic boundaries. The fluid is inviscid and heat conduction is neglected, thus we set $\mu=\kappa=0$, and some perturbations are initially assigned on the top of a homogeneous background field:
\begin{equation}
	\Pv(\xx,0) = (1+\delta \rho, 1+\delta u, 1+\delta v, 1+\delta p).
	\label{eq.ConvEul-IC}
\end{equation}
The perturbations for temperature $\delta T$, density $\delta \rho$ and pressure $\delta p$ read
\begin{equation}
	\label{ShuVortDelta1}
	\delta T = -\frac{(\gamma-1)\varepsilon^2}{8\gamma\pi^2}e^{1-r^2}, \quad
	\delta \rho = (1+\delta T)^{\frac{1}{\gamma-1}}-1,  \quad 
	\delta p = (1+\delta T)^{\frac{\gamma}{\gamma-1}}-1,
\end{equation}
where the radius of the vortex has been defined as $r=\sqrt{(x-5)^2+(y-5)^2}$ and the vortex strength has been set to $\varepsilon=5$. The velocity perturbations are given by
\begin{equation}
	\label{ShuVortDelta2}
	\left(\begin{array}{c} \delta u \\ \delta v \\ \end{array}\right) = \frac{\varepsilon}{2\pi}e^{\frac{1-r^2}{2}} \left(\begin{array}{c} -(y-5) \\ \phantom{-}(x-5) \end{array}\right).
\end{equation}
The final time of the simulation is chosen to be $t_f=0.1$, and the exact solution $\Pv_e(\xx,t)$ can be simply computed as the time-shifted initial condition by the convective velocity $\vv_c=(1,1)$, that is $\Pv_e(\xx,t_f)=\Pv_e(\xx-\vv_c \, t_f,0)$. The error is measured at the final time in $L_2$ norm as
\begin{equation}
	\epsilon_{L_2}= \sqrt{\int \limits_{\Omega} \left( \Pv_e(\xx,t_f)-\uu_h(\xx,t_f)\right)^2 \, \text{d}\xx}.
\end{equation}

A sequence of successively refined unstructured Voronoi meshes of characteristic mesh size $h_{\Omega}$ is used to perform the convergence analysis. The VEM-DG schemes are compared against two different space-time DG methods, namely: i) the classical modal DG (M-DG) schemes which make use of modal basis functions, which are the monomials \eqref{eq:m_alpha}, and ii) the Agglomerated Finite Element DG (AFE-DG) methods recently forwarded in \cite{ADERAFEDG}, which are based on a local finite element basis within each control volume. We use a $\CFL$ number of $\CFL=0.25$.

\begin{figure}[!htbp]
	\begin{center}
		\begin{tabular}{cc} 
			\includegraphics[width=0.47\textwidth]{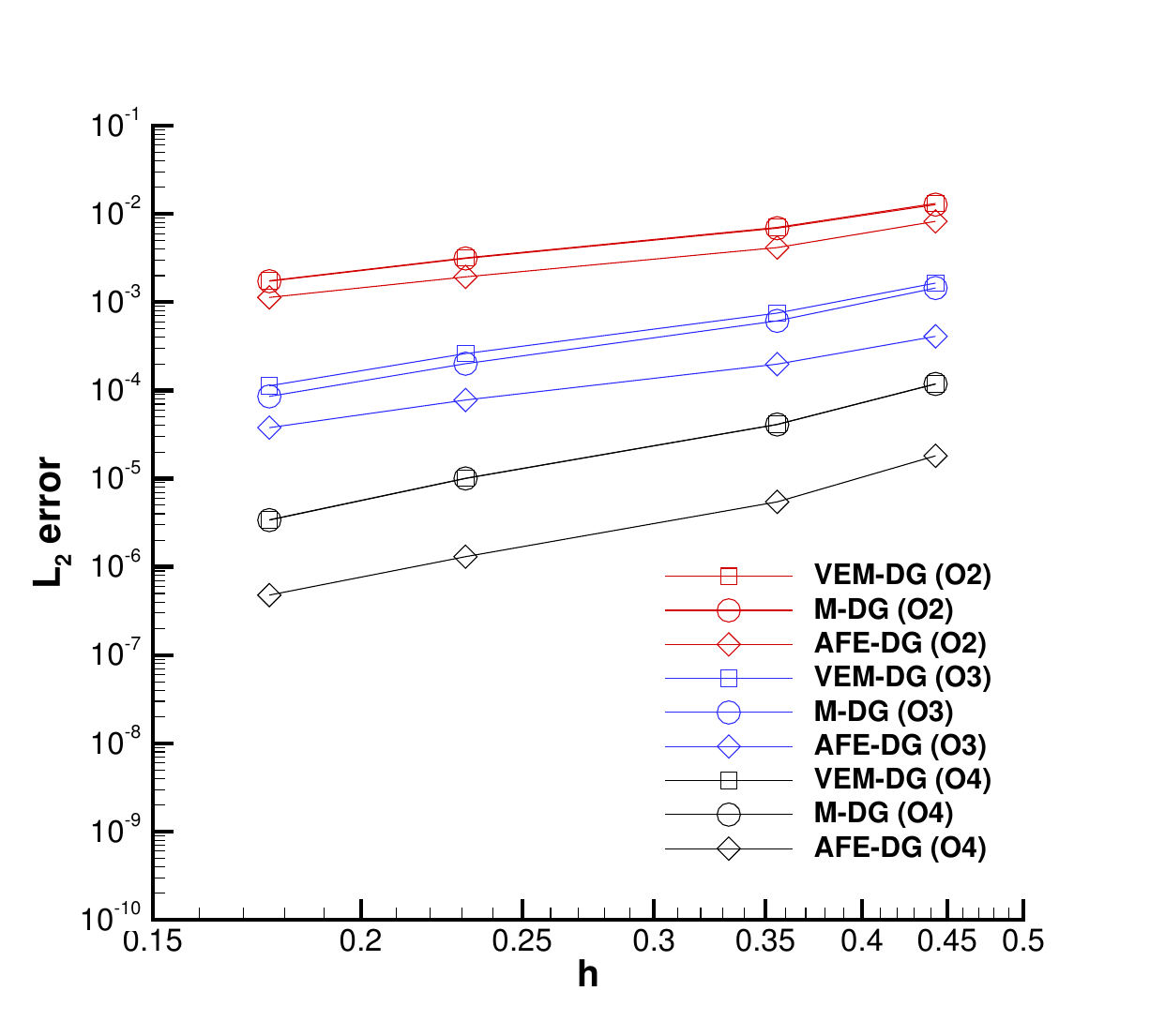} & 
			\includegraphics[width=0.47\textwidth]{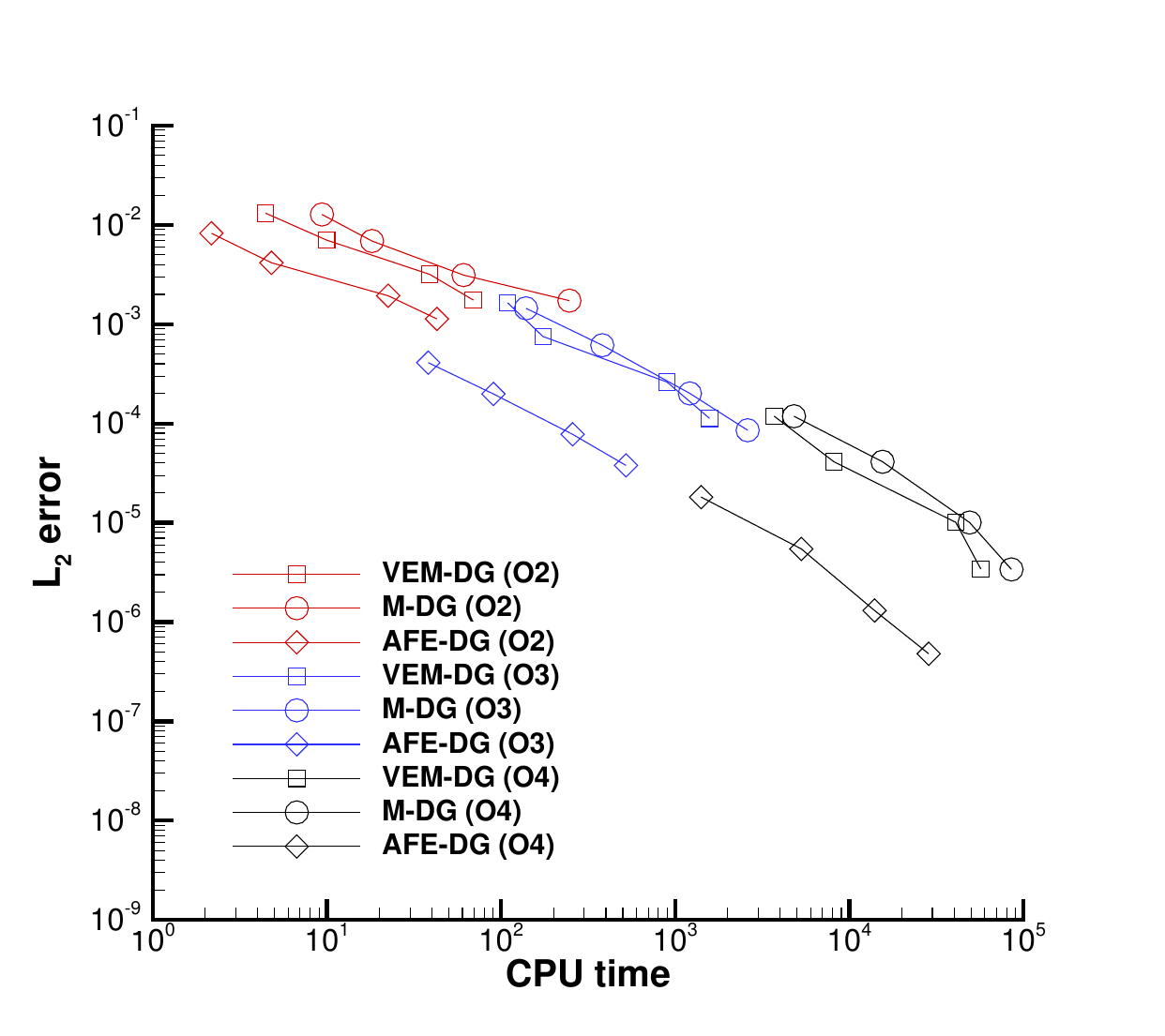}\\   
		\end{tabular} 
		\caption{Comparison between VEM-DG (squares), M-DG (circles) and AFE-DG (diamonds) schemes from second up to fourth order of accuracy. Left: dependency of the error norm on the mesh size. Right: dependency of the error norm on the CPU time.}
		\label{fig.comparison_efficiency}
	\end{center}
\end{figure}

The convergence results are collected in Table \ref{tab.conv}, which demonstrate that the formal order of convergence is attained by the VEM-DG schemes. Figure \ref{fig.comparison_efficiency} compares the three different methods, showing that the VEM-DG schemes are more efficient than the standard M-DG methods, while being less efficient and less accurate than the AFE-DG schemes. This is expected since the AFE-DG schemes are built upon a sub-triangulation of the Voronoi elements, thus allowing each sub-triangle to be mapped to a reference element and ultimately to design a fully quadrature-free DG scheme, as detailed in \cite{ADERAFEDG}. Nevertheless, compared to the commonly used M-DG methods, the novel VEM-DG schemes are more efficient and, potentially, they also easily allows mixed elements and nonconforming grids to be handled.

\begin{table}[!htbp]  
	\caption{Numerical convergence results for the compressible Euler equations using VEM-DG, M-DG and AFE-DG schemes from second up to fourth order of accuracy in space and time. The errors are measured in the $L_2$ norm and refer to density ($\rho$) at time $t_f=0.1$. The absolute CPU time of each simulation is also reported in seconds $[s]$. The characteristic mesh size is given by $h_{\Omega}=\max \limits_{\mathcal{T}_\Omega} \hP$.}  
	\begin{center} 
		\begin{small}
			\renewcommand{\arraystretch}{1.1}
			\begin{tabular}{c|ccc|ccc|ccc} 
				\multicolumn{1}{c}{} &  \multicolumn{3}{c}{VEM-DG} & \multicolumn{3}{c}{M-DG} & \multicolumn{3}{c}{AFE-DG} \\
				$h_{\Omega}$ & $\rho_{L_2}$ & $\mathcal{O}(\rho_{L_2})$ & CPU time & $\rho_{L_2}$ & $\mathcal{O}(\rho_{L_2})$ & CPU time & $\rho_{L_2}$ & $\mathcal{O}(\rho_{L_2})$ & CPU time \\
				\hline
				\multicolumn{1}{c|}{} & \multicolumn{9}{c}{Order of accuracy: $\mathcal{O}(2)$} \\
				4.428E-01 & 1.315E-02 &    - & 4.438E+00 & 1.277E-02 &    - & 9.359E+00 & 8.239E-03 & - & 2.172E+00 \\
				3.557E-01 & 7.039E-03 & 2.85 & 9.953E+00 & 6.892E-03 & 2.81 & 1.814E+01 & 4.161E-03 & 3.12 & 4.797E+00 \\
				2.311E-01 & 3.178E-03 & 1.85 & 3.894E+01 & 3.130E-03 & 1.83 & 6.086E+01 & 1.936E-03 & 1.78 & 2.247E+01 \\
				1.762E-01 & 1.751E-03 & 2.20 & 6.930E+01 & 1.728E-03 & 2.19 & 2.467E+02 & 1.132E-03 & 1.98 & 4.278E+01 \\
				\multicolumn{1}{c|}{} & \multicolumn{9}{c}{Order of accuracy: $\mathcal{O}(3)$} \\
				4.428E-01 & 1.646E-03 &    - & 1.090E+02 & 1.447E-03 &    - & 1.390E+02 & 4.102E-04 & - & 2.172E+00 \\
				3.557E-01 & 7.524E-04 & 3.57 & 1.740E+02 & 6.141E-04 & 3.91 & 3.813E+02 & 1.985E-04 & 3.31 & 4.797E+00 \\
				2.311E-01 & 2.615E-04 & 2.45 & 8.972E+02 & 2.009E-04 & 2.59 & 1.214E+03 & 7.785E-05 & 2.17 & 2.247E+01 \\
				1.762E-01 & 1.128E-04 & 3.10 & 1.575E+03 & 8.547E-05 & 3.15 & 2.606E+03 & 3.783E-05 & 2.66 & 4.278E+01 \\
				\multicolumn{1}{c|}{} & \multicolumn{9}{c}{Order of accuracy: $\mathcal{O}(4)$} \\
				4.428E-01 & 1.184E-04 &    - & 3.699E+03 & 1.184E-04 &    - & 4.809E+03 & 1.814E-05 & - & 1.409E+03 \\
				3.557E-01 & 4.113E-05 & 4.83 & 8.167E+03 & 4.111E-05 & 4.83 & 1.552E+04 & 5.451E-06 & 5.49 & 5.301E+03 \\
				2.311E-01 & 1.009E-05 & 3.26 & 4.088E+04 & 1.005E-05 & 3.27 & 4.906E+04 & 1.305E-06 & 3.32 & 1.397E+04 \\
				1.762E-01 & 3.417E-06 & 3.99 & 5.692E+04 & 3.390E-06 & 4.00 & 8.539E+04 & 4.790E-07 & 3.69 & 2.857E+04 \\
			\end{tabular}
		\end{small}
	\end{center}
	\label{tab.conv}
\end{table}

\begin{figure}[!htbp]
	\begin{center}
		\begin{tabular}{ccc} 
			\includegraphics[trim=2 0 0 5,clip,width=0.32\textwidth]{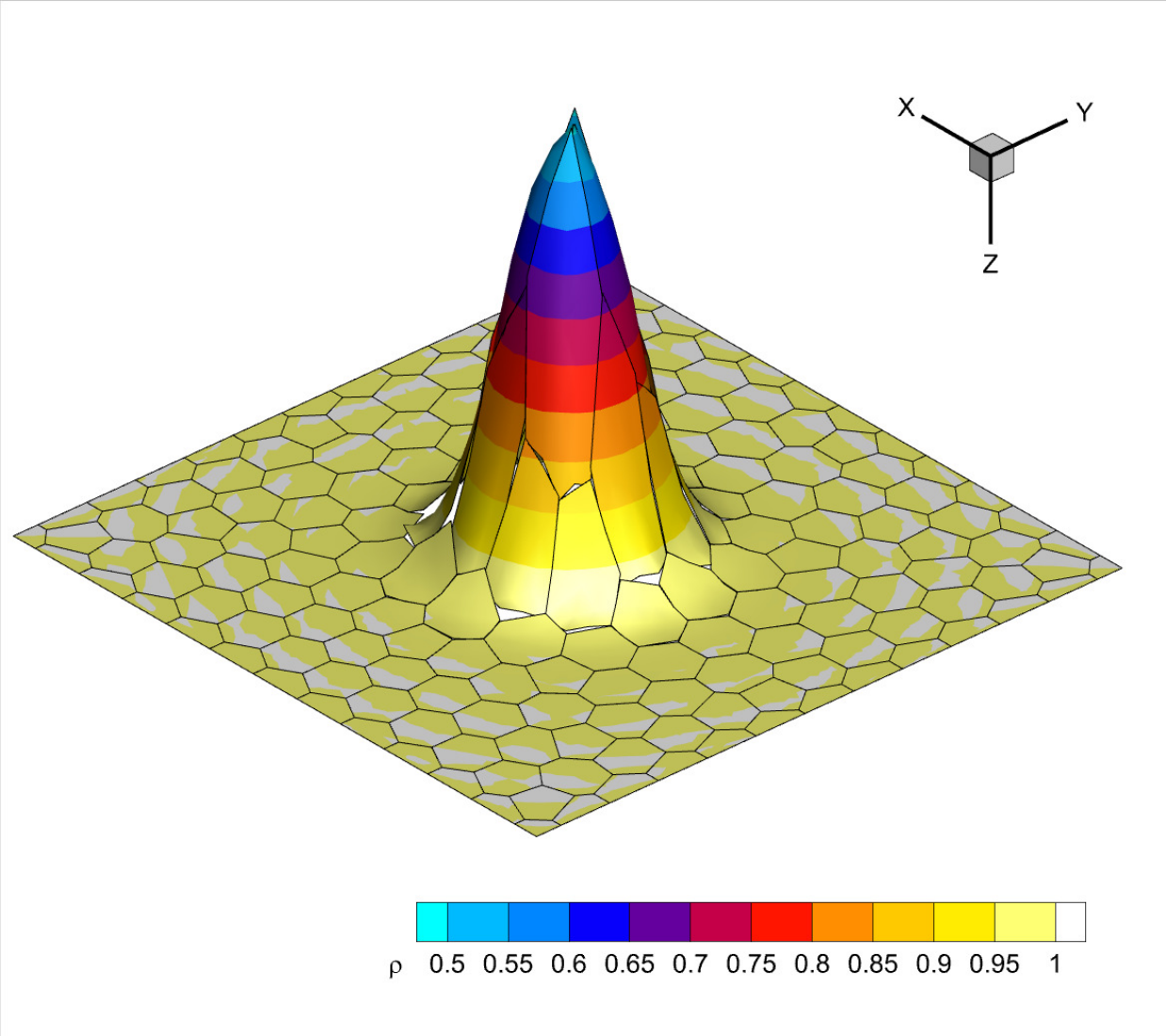} & 
			\includegraphics[trim=2 0 0 5,clip,width=0.32\textwidth]{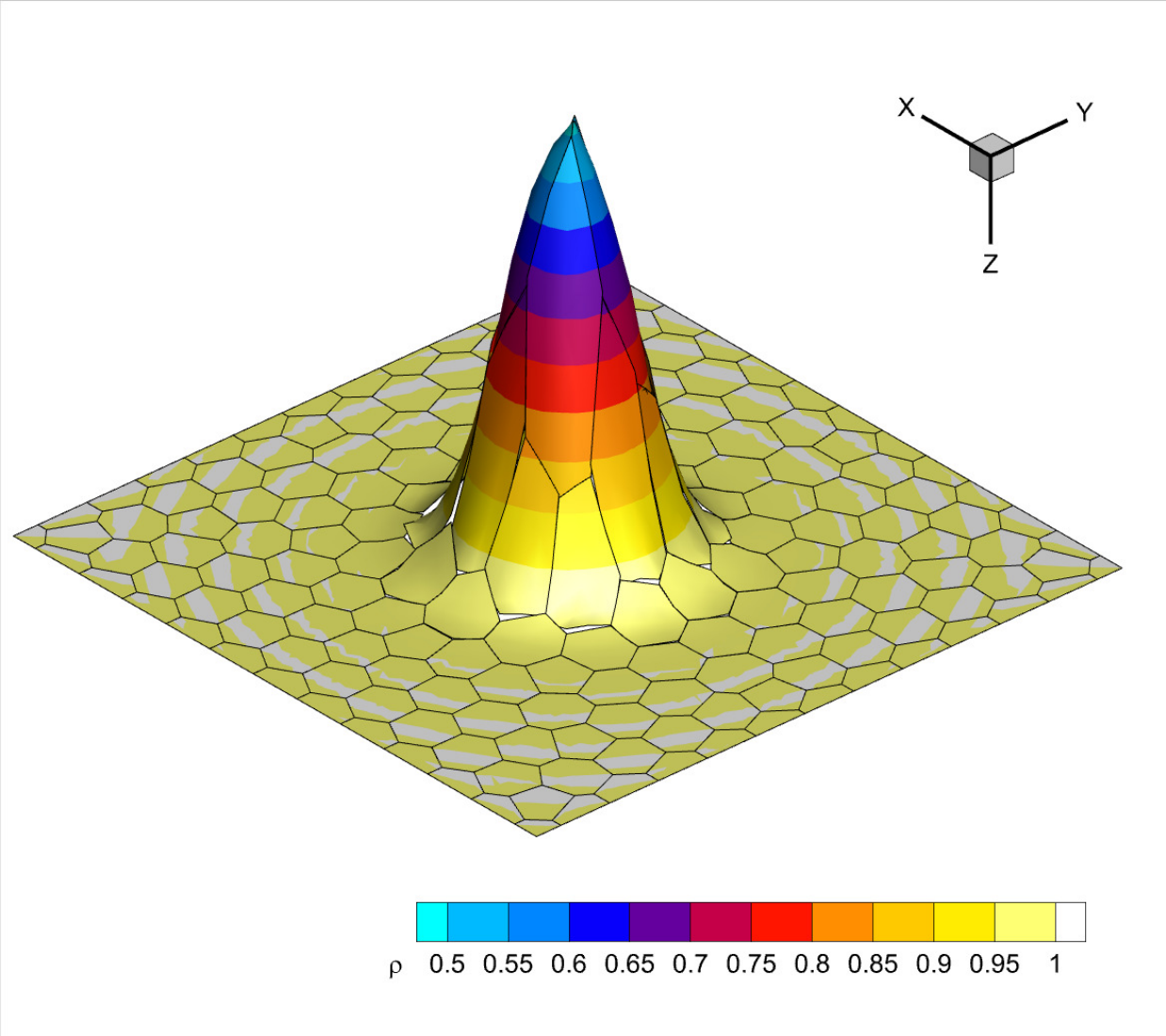} & 
			\includegraphics[trim=2 0 0 5,clip,width=0.32\textwidth]{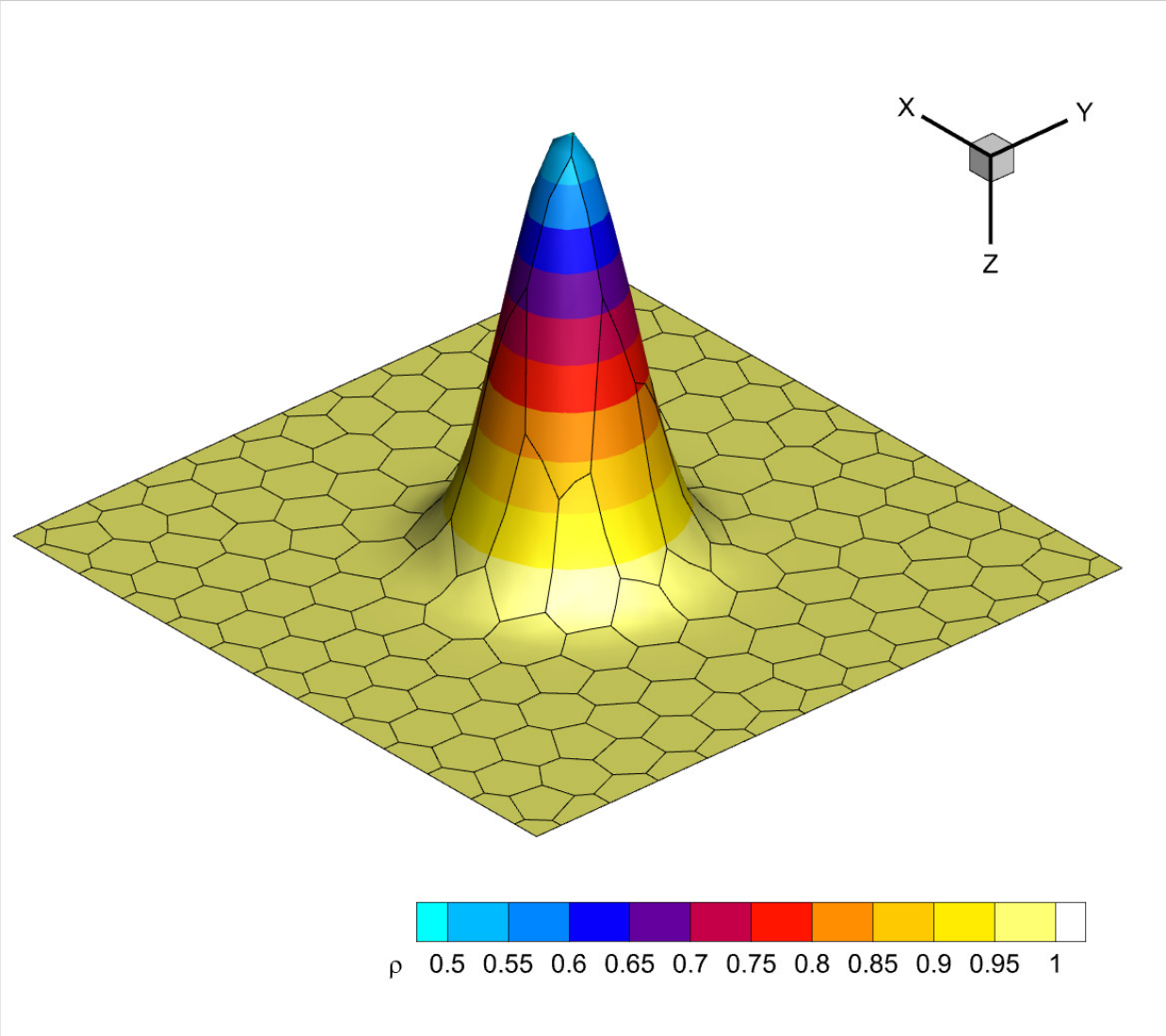}\\  
			\includegraphics[trim=2 0 0 5,clip,width=0.32\textwidth]{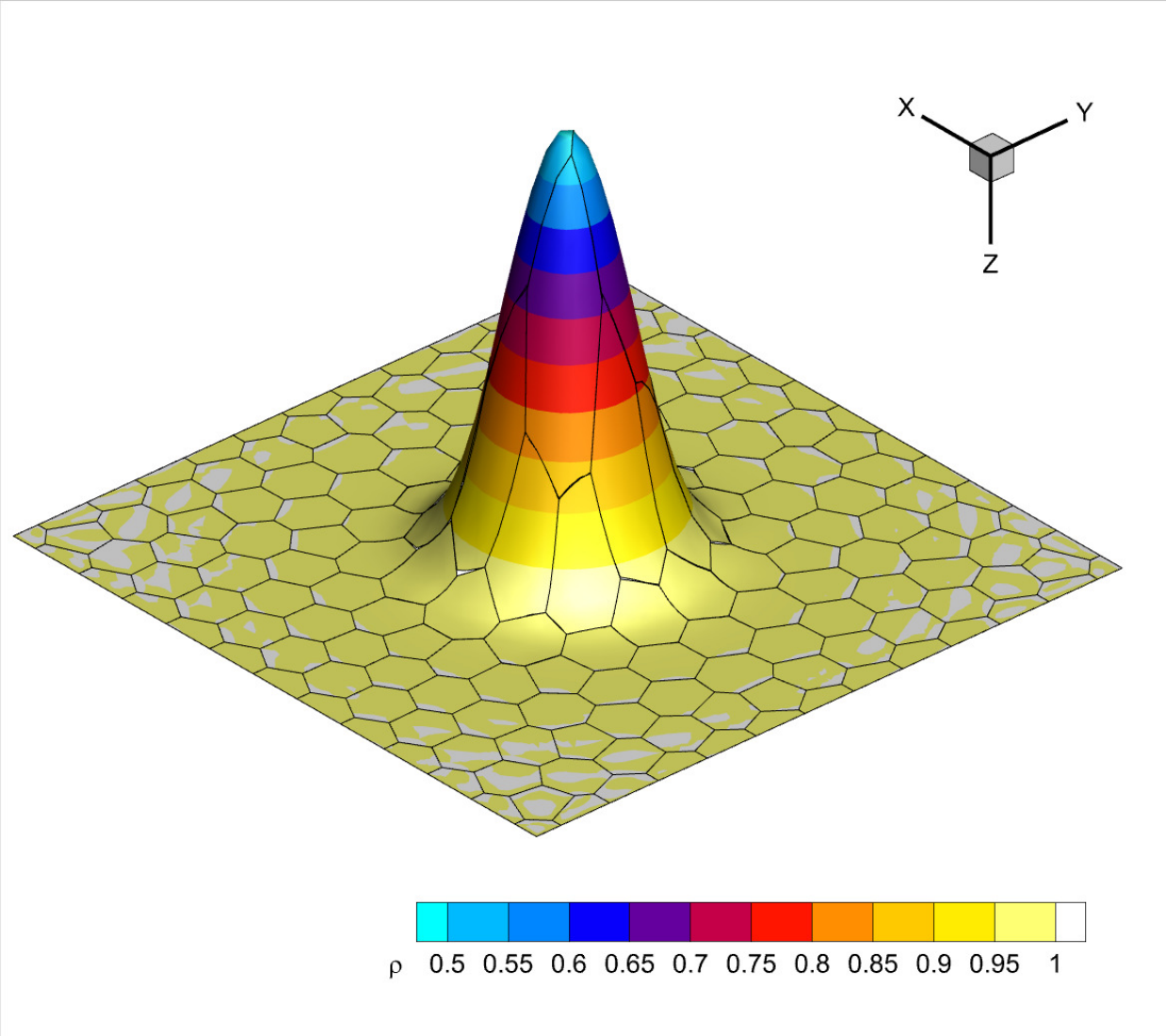} & 
			\includegraphics[trim=2 0 0 5,clip,width=0.32\textwidth]{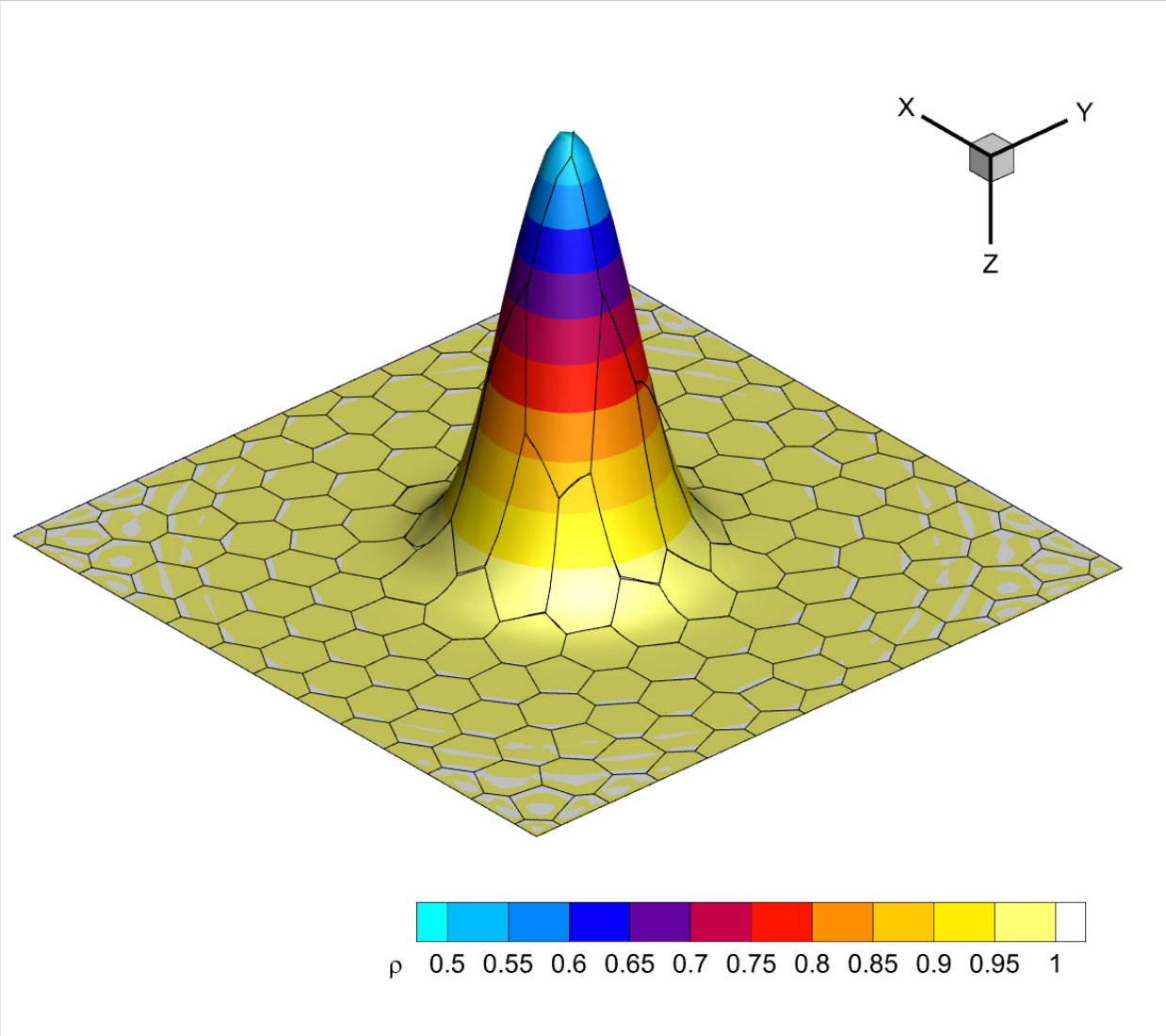} & 
			\includegraphics[trim=2 0 0 5,clip,width=0.32\textwidth]{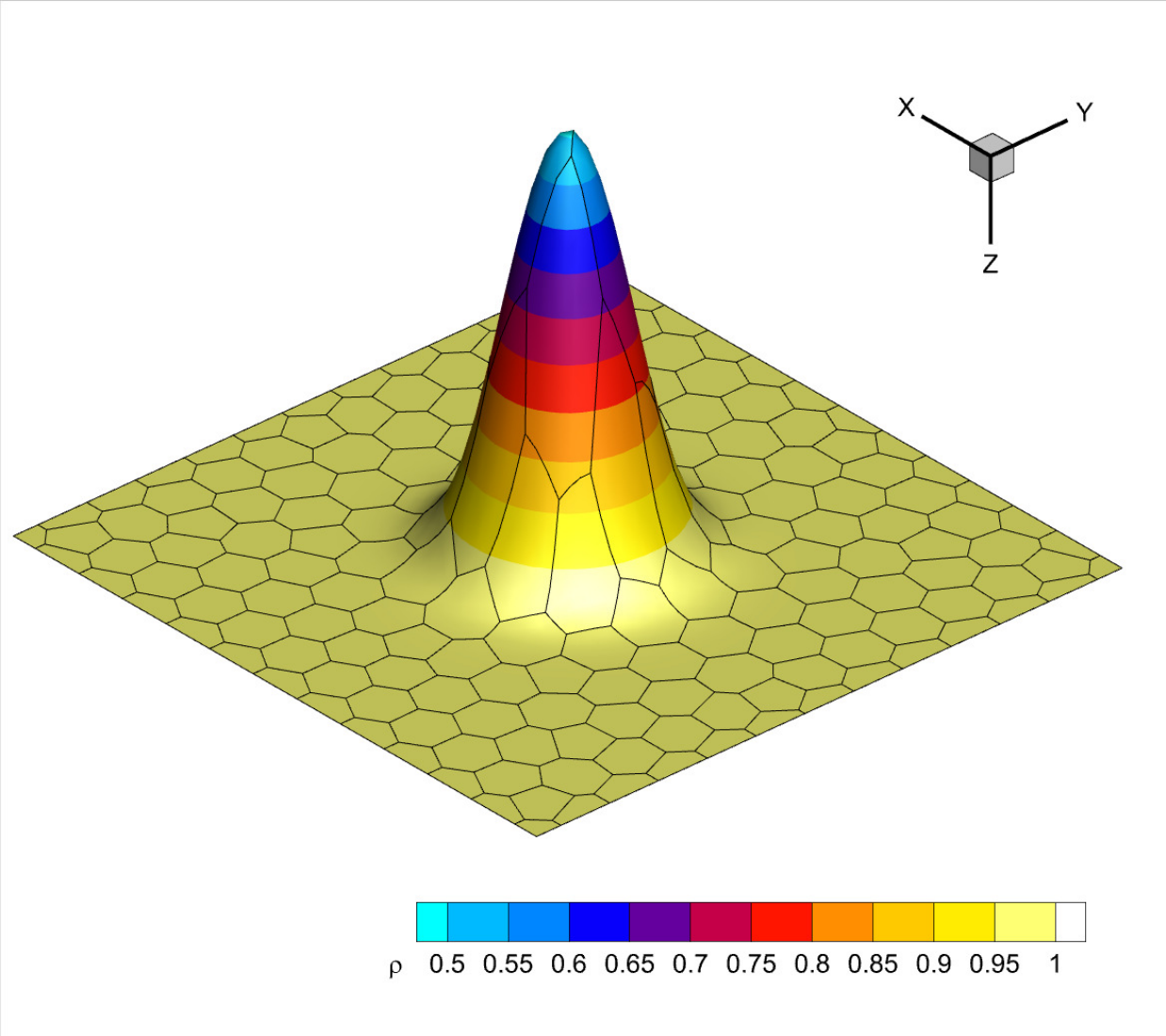}\\ 
		\end{tabular} 
		\caption{$L_2$ projection of the density distribution of the initial condition for the isentropic vortex test case with $N=2$ (top row) and $N=3$ (bottom row) on a mesh with characteristic size $h_{\Omega}=10/12$. Left: VEM-DG. Middle: M-DG. Right: AFE-DG.}
		\label{fig.ShuVortex-IC}
	\end{center}
\end{figure}

Table \ref{tab.IC} reports the errors related to the $L_2$ projection of the initial condition \eqref{eq.ConvEul-IC} for the three methods, showing that the VEM-DG schemes behave quite similar to the M-DG methods. Figure \ref{fig.ShuVortex-IC} depicts the initial condition on a very coarse mesh of characteristic mesh size $h_{\Omega}=5/6$ for all the DG methods considered here, which qualitatively confirms the results obtained from the analysis of Table \ref{tab.IC}.

\begin{table}[!htbp]  
	\caption{Errors related to the $L_2$ projection of the initial condition for the isentropic vortex test case measured in $L_2$ norm for density ($\rho$), horizontal velocity ($u$) and pressure ($p$). The characteristic mesh size is $h_{\Omega}=10/12$ and the errors are reported for third and fourth order VEM-DG, M-DG and AFE-DG schemes.}  
	\begin{center} 
		\begin{small}
			\renewcommand{\arraystretch}{1.0}
			\begin{tabular}{c|ccc|ccc|ccc} 
				\multicolumn{1}{c}{} &  \multicolumn{3}{c}{VEM-DG} & \multicolumn{3}{c}{M-DG} & \multicolumn{3}{c}{AFE-DG} \\
				$N$ & $\rho_{L_2}$ & $u_{L_2}$ & $p_{L_2}$ & $\rho_{L_2}$ & $u_{L_2}$ & $p_{L_2}$ & $\rho_{L_2}$ & $u_{L_2}$ & $p_{L_2}$ \\
				\hline
				2 & 1.145E-02 & 2.348E-02 & 1.726E-02 & 9.575E-03 & 1.787E-02 & 1.448E-02 & 2.041E-03 & 3.872E-03 & 2.990E-03 \\
				3 & 2.119E-03 & 5.276E-03 & 3.460E-03 & 1.506E-03 & 3.669E-03 & 2.493E-03 & 1.573E-04 & 3.751E-04 & 2.433E-04 \\
			\end{tabular}
		\end{small}
	\end{center}
	\label{tab.IC}
\end{table}

\subsection{Analysis of the condition number of the VEM matrices} \label{ssec.cond_num}
\begin{figure}[!htbp]
	\begin{center}
		\begin{tabular}{cc} 
			\includegraphics[trim= 2 2 2 2, clip,width=0.47\textwidth]{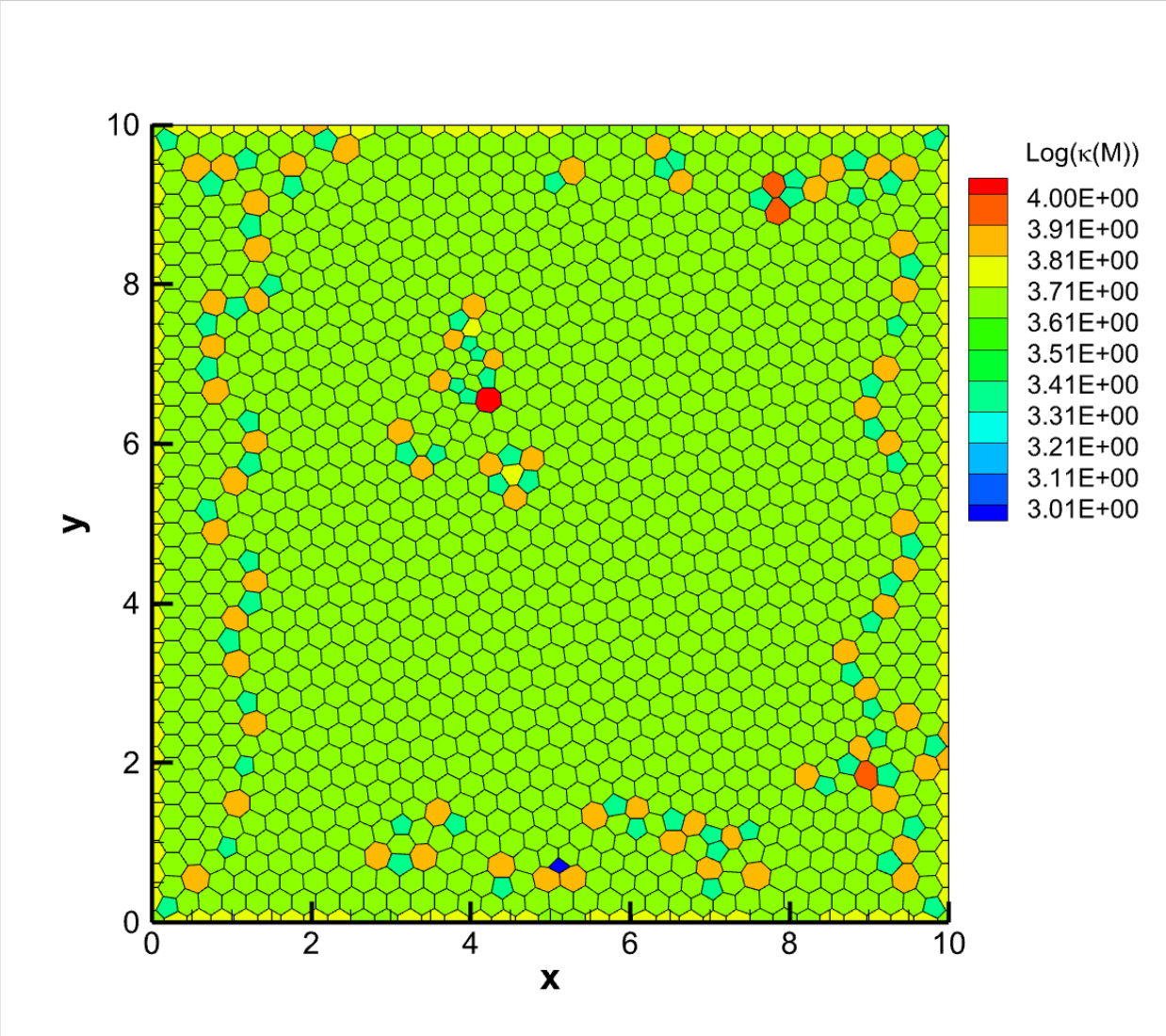} & 
			\includegraphics[trim= 2 2 2 2, clip,width=0.47\textwidth]{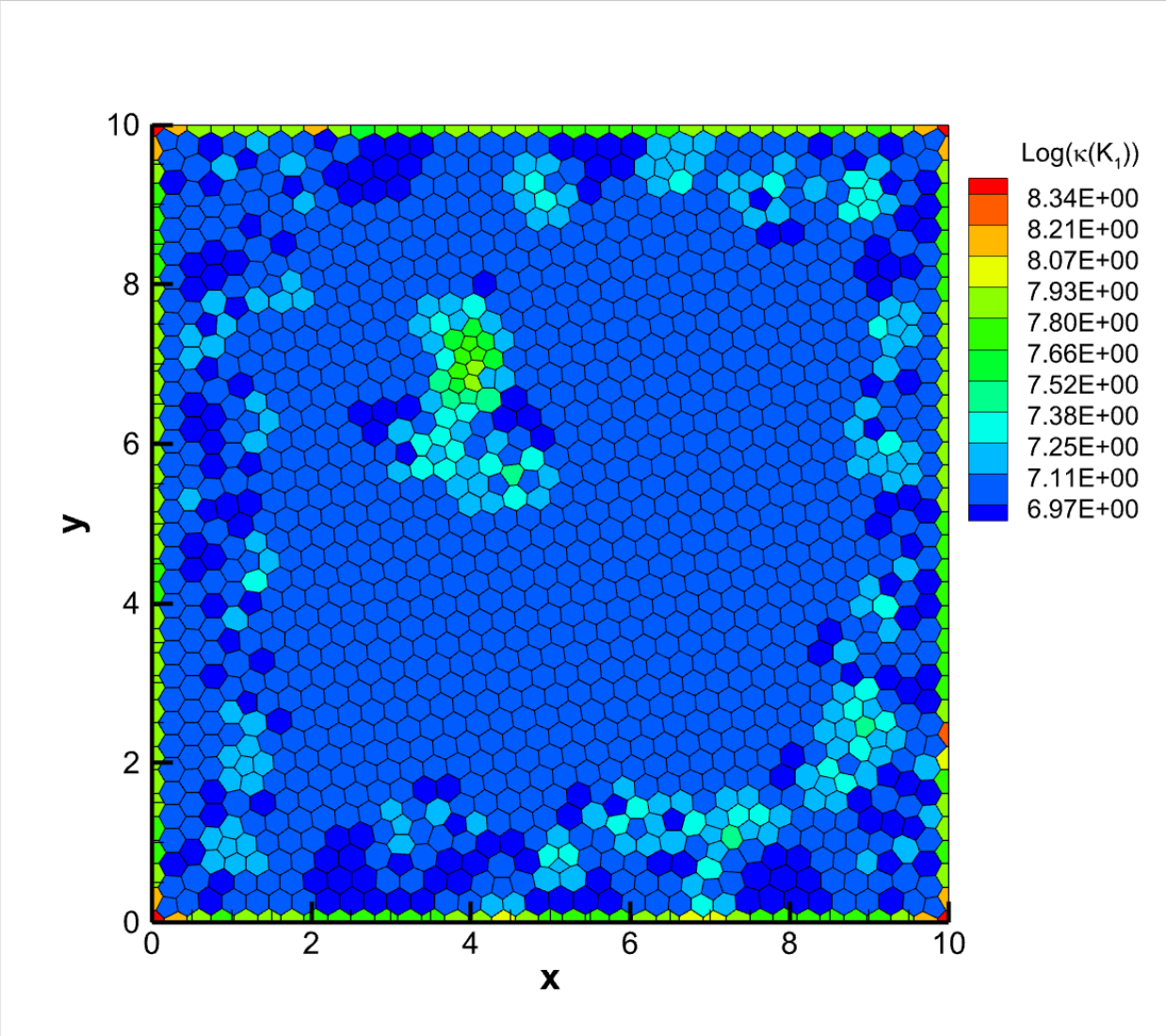} \\
			\includegraphics[trim= 2 2 2 2, clip,width=0.47\textwidth]{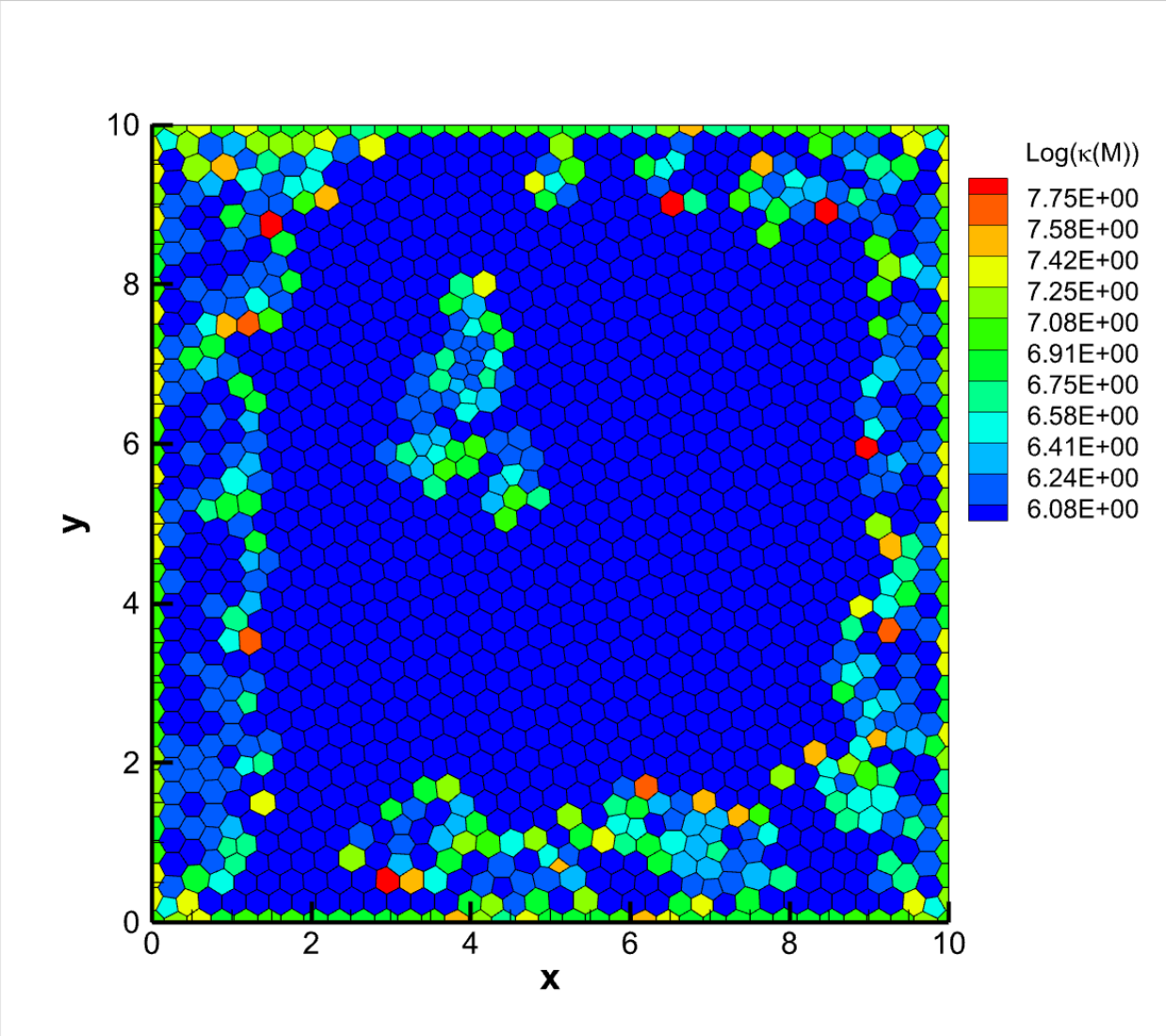} & 
			\includegraphics[trim= 2 2 2 2, clip,width=0.47\textwidth]{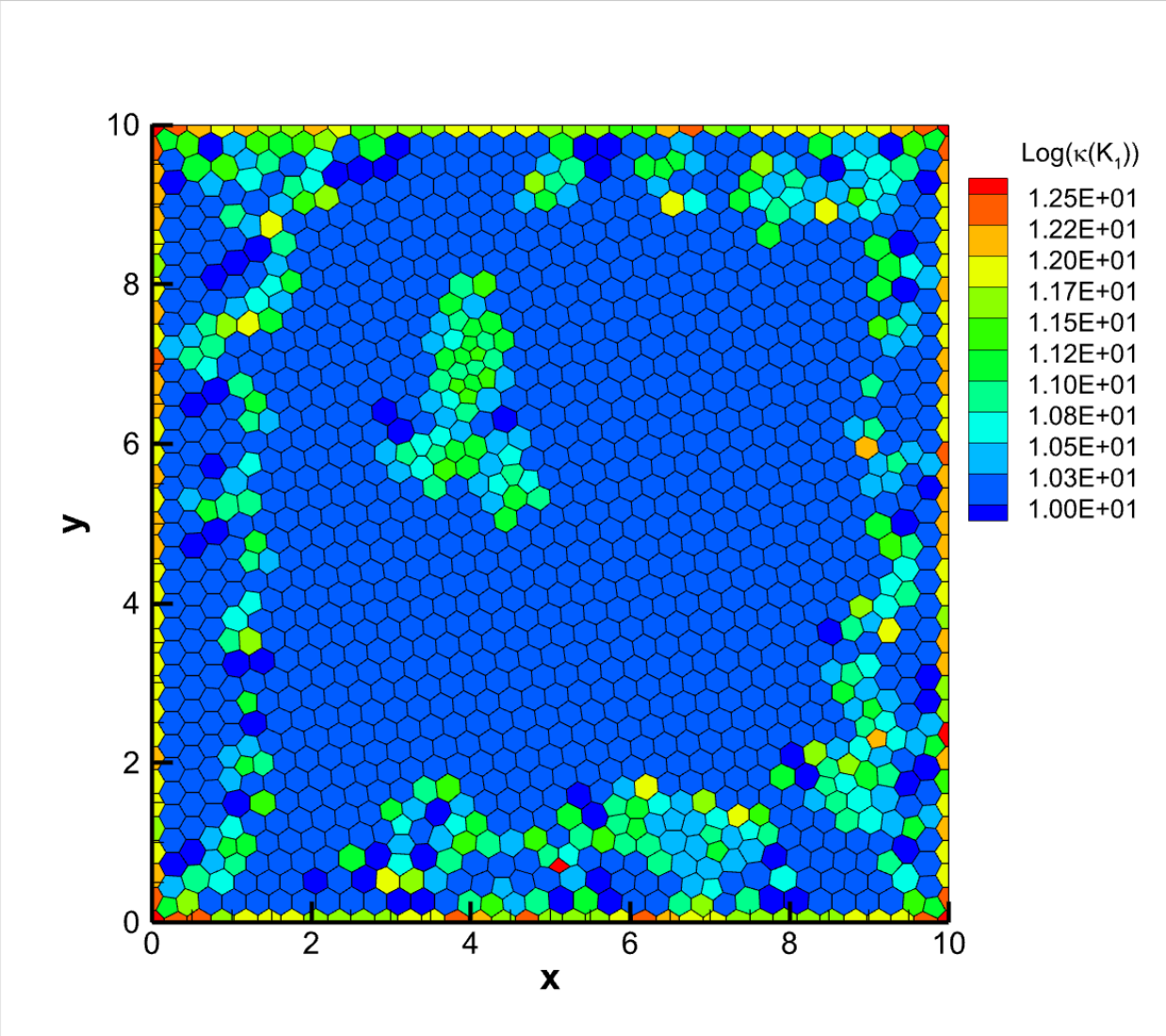} \\
			\includegraphics[trim= 2 2 2 2, clip,width=0.47\textwidth]{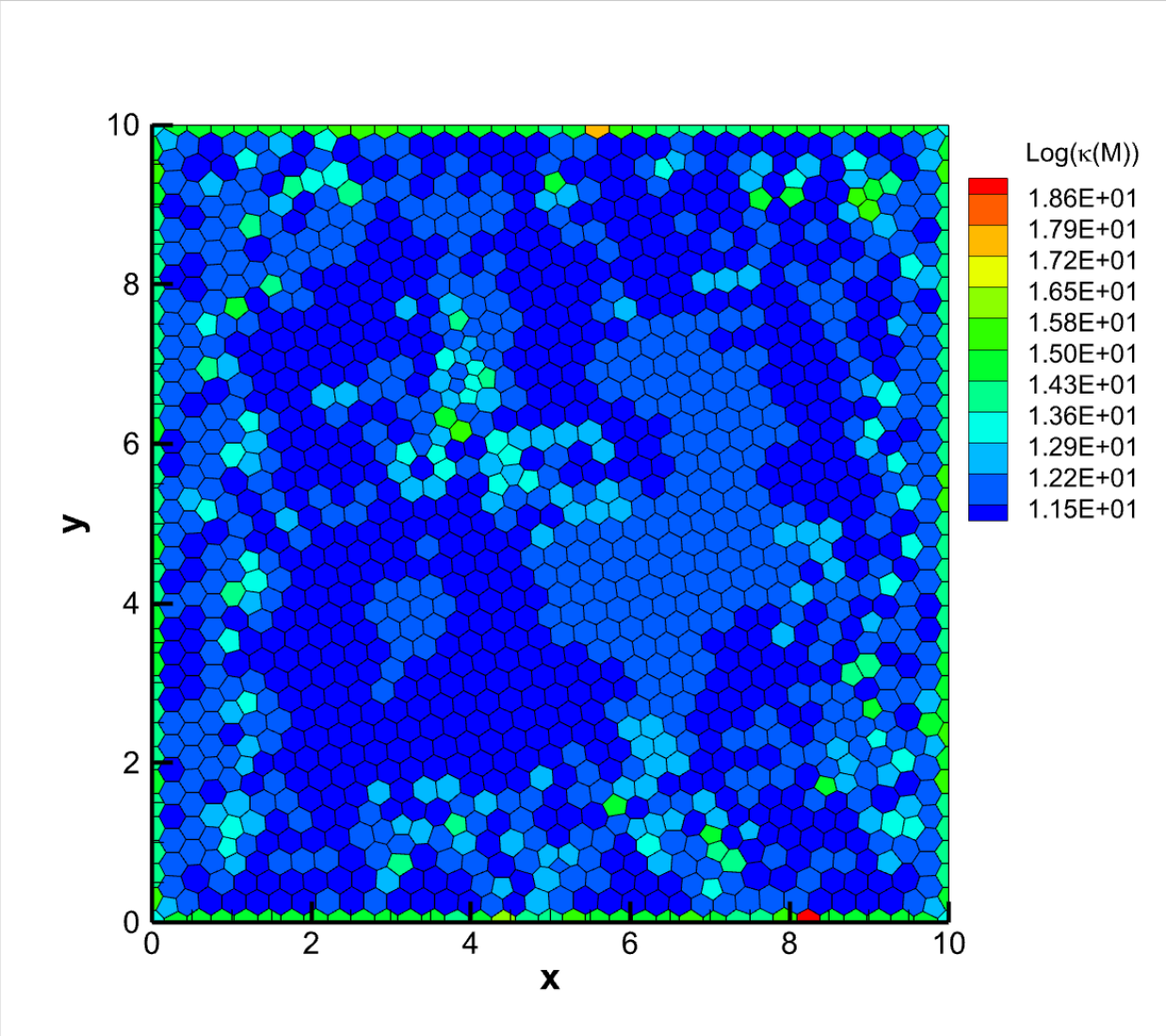} & 
			\includegraphics[trim= 2 2 2 2, clip,width=0.47\textwidth]{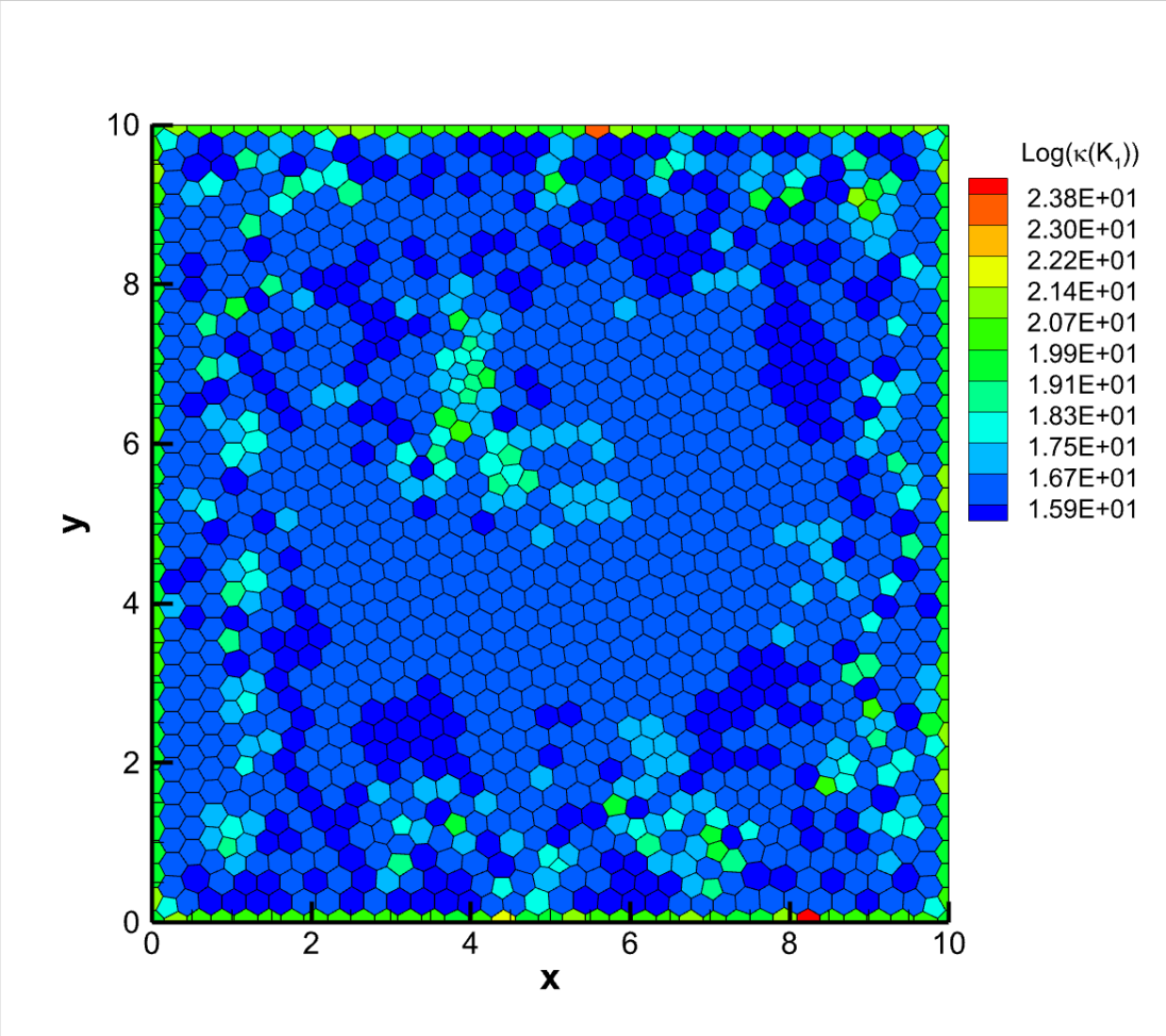} \\
		\end{tabular} 
		\caption{Analysis of the condition number of the VEM matrices employed in the new VEM-DG schemes for $N=1$ (top row), $N=2$ (middle row) and $N=3$ (bottom row) for a mesh with characteristic size $h_{\Omega}=1/3$. Left column: logarithm of the condition number of the space mass matrix $\mathbf{M}$. Right column: logarithm of the condition number of the space-time stiffness matrix $\mathbf{K}_1$.  }
		\label{fig.cond_num}
	\end{center}
\end{figure}

Next, we analyze the condition number of the matrices that need to be inverted in the novel VEM-DG schemes. Specifically, we consider the space mass matrix $\mathbf{M}$ in the corrector step given by \eqref{eqn.proj_M3} and the space-time stiffness matrix $\mathbf{K}_1$ in the predictor step defined in \eqref{eqn.K1_final}. The condition number $\kappa(\mathbf{A})$ of a generic matrix $\mathbf{A}$ of size $n \times n$ is defined as
\begin{equation} \label{eqn.cond_num}
	\kappa(\mathbf{A}) = ||\mathbf{A}||_F \, ||\mathbf{A}^{-1}||_F, 
\end{equation}
where we use the Frobenius norm
\begin{equation}
	||\mathbf{A}||_F =  \sqrt{\sum_{i=1}^{n} \sum_{j=1}^{n} a_{ij}^2},
\end{equation}
with $a_{ij}$ denoting the entries of matrix $\mathbf{A}$. Two meshes are used, taken from the previous convergence test case, with characteristic mesh size $h_{\Omega,1}=1/3$ and $h_{\Omega,2}=1/6$, and the condition numbers are computed for virtual basis of degree $N=\{1,2,3\}$. The results are gathered in Table \ref{tab.cond_num}, where it is evident that the higher the polynomial degree of the virtual basis, the worse the condition number of the matrices. In particular, the space-time stiffness matrix is always worse conditioned compared to the space mass matrix of approximately two orders of magnitude. Furthermore, for $N=3$ the condition numbers become very large, hence making the corresponding VEM-DG schemes less robust. These results are in agreement with the findings reported in \cite{Berrone2017,Mascotto2018}, where it is indeed shown that the condition number of the virtual basis deteriorates very rapidly starting from $N=3$. In Table \ref{tab.cond_num} we measure the minimum, maximum and average condition number of each matrix. The average quantity is simply given by the arithmetic average of the condition number over all the cells of the computational domain. Figure \ref{fig.cond_num} depicts the condition number for each cell of the coarser mesh employed for this analysis (we use the logarithm of the condition number to improve the readability of the results due to the involved very large scales). It is interesting to notice that the worst ill-conditioning is mainly located at boundary elements or at internal cells with a quite irregular shape. In both cases, the characteristic mesh size $h_{P}$, as defined in equation \eqref{eqn.h}, is no longer a suitable measure of the size of the cell, since the aspect ratio is quite far from unity, hence yielding a local bad mesh quality. On the other side, we notice that the internal cells which exhibit a bad condition number are characterized by high differences in the length of the boundary edges, hence involving polygons with more than six edges or highly stretched control volumes with only four sides.

This simple analysis suggests that: i) more investigations are needed in the stabilization terms of the VEM matrices, i.e. the space mass matrix $\mathbf{M}$ and the space-time stiffness matrix $\mathbf{K}_1$; ii) the aspect ratio of the control volume highly affects the condition number of the VEM matrices, hence requiring novel approaches to incorporate this uncertainty in the definition of the local VEM space.  

\begin{table}[!htbp]  
	\caption{Minimum $\kappa_{\min}$, maximum $\kappa_{\max}$ and average $\kappa_{\text{av}}$ condition number of the space mass matrix $\mathbf{M}$ and of the space-time stiffness matrix $\mathbf{K}_1$ for two different meshes of characteristic size $h_{\Omega,1}=1/3$ and $h_{\Omega,2}=1/6$. We consider Virtual Element basis of degree $N=\{1,2,3\}$.}  
	\begin{center} 
		\begin{small}
			\renewcommand{\arraystretch}{1.0}
			\begin{tabular}{c|ccc|ccc}
				\hline
				\multicolumn{7}{c}{$N=1$} \\
				\multicolumn{1}{c}{} &  \multicolumn{3}{c}{$\kappa(\mathbf{M})$} & \multicolumn{3}{c}{$\kappa(\mathbf{K}_1)$} \\
				\hline
				$h_{\Omega}$ & $\kappa_{\min}$ & $\kappa_{\max}$ & $\kappa_{\text{av}}$ & $\kappa_{\min}$ & $\kappa_{\max}$ & $\kappa_{\text{av}}$ \\
				\hline
				1/3 & 1.835E+01 & 6.059E+01 & 3.774E+01 & 8.788E+02 & 4.687E+03 & 1.301E+03 \\
				1/6 & 1.812E+01 & 6.424E+01 & 3.751E+01 & 2.938E+03 & 2.093E+04 & 5.012E+03 \\
				\multicolumn{7}{c}{} \\
				\hline
				\multicolumn{7}{c}{$N=2$} \\
				\multicolumn{1}{c}{} &  \multicolumn{3}{c}{$\kappa(\mathbf{M})$} & \multicolumn{3}{c}{$\kappa(\mathbf{K}_1)$} \\
				\hline
				$h_{\Omega}$ & $\kappa_{\min}$ & $\kappa_{\max}$ & $\kappa_{\text{av}}$ & $\kappa_{\min}$ & $\kappa_{\max}$ & $\kappa_{\text{av}}$ \\
				\hline
				1/3 & 3.616E+02 & 2.802E+03 & 5.666E+02 & 1.855E+04 & 3.151E+05 & 4.382E+04 \\
				1/6 & 3.584E+02 & 2.001E+04 & 5.107E+02 & 6.902E+04 & 3.563E+06 & 1.448E+05 \\
				\multicolumn{7}{c}{} \\
				\hline
				\multicolumn{7}{c}{$N=3$} \\
				\multicolumn{1}{c}{} &  \multicolumn{3}{c}{$\kappa(\mathbf{M})$} & \multicolumn{3}{c}{$\kappa(\mathbf{K}_1)$} \\
				\hline
				$h_{\Omega}$ & $\kappa_{\min}$ & $\kappa_{\max}$ & $\kappa_{\text{av}}$ & $\kappa_{\min}$ & $\kappa_{\max}$ & $\kappa_{\text{av}}$ \\
				\hline
				1/3 & 4.599E+04 & 2.566E+08 & 5.548E+05 & 2.899E+06 & 6.082E+10 & 1.183E+08 \\
				1/6 & 4.295E+04 & 1.777E+09 & 6.075E+05 & 1.010E+07 & 1.622E+12 & 5.079E+08 \\
			\end{tabular}
		\end{small}
	\end{center}
	\label{tab.cond_num}
\end{table}

\subsection{First problem of Stokes}
\begin{figure}[!htbp]
	\begin{center}
		\begin{tabular}{cc} 
			\includegraphics[trim= 2 2 2 2, clip,width=0.47\textwidth]{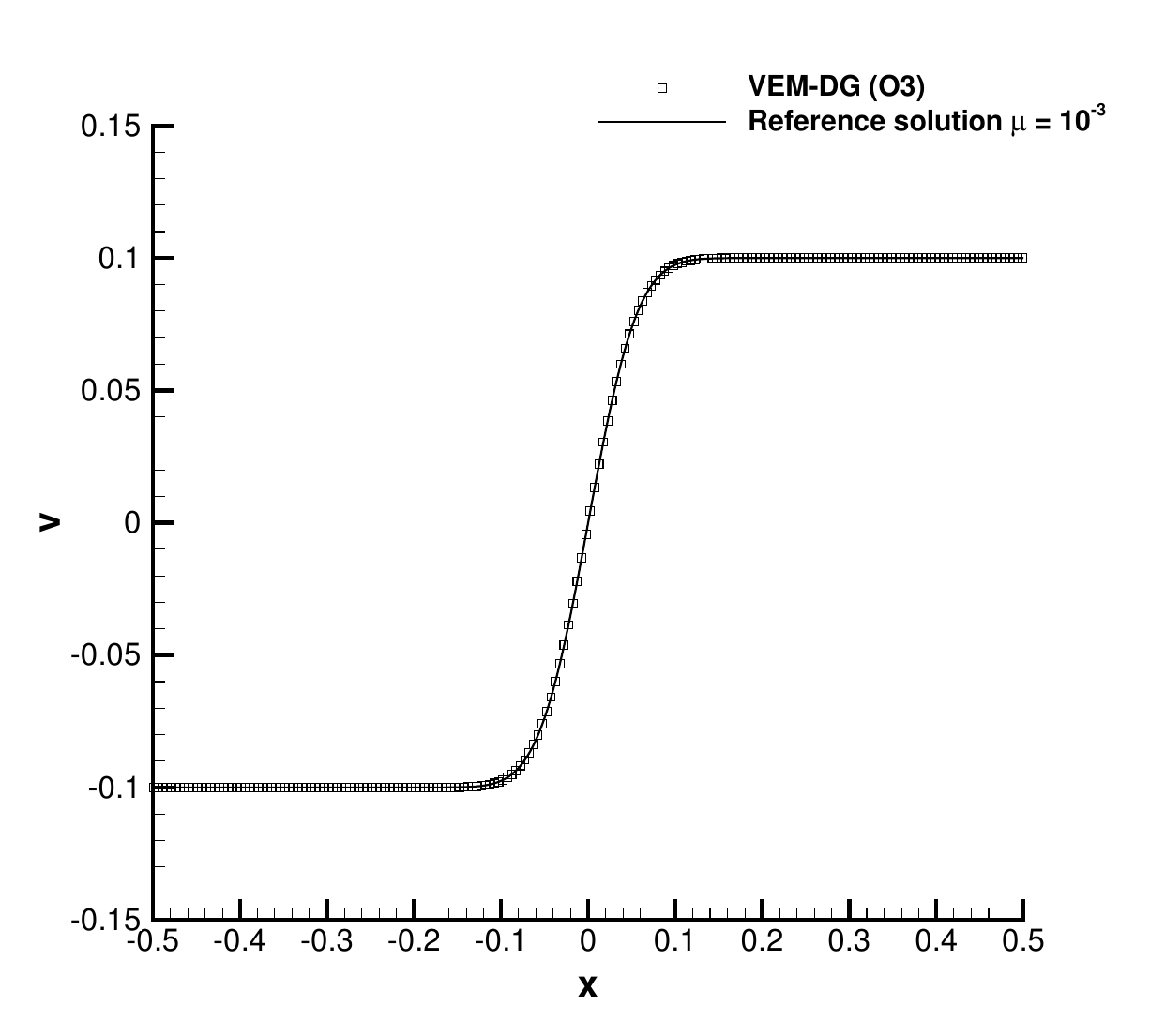} & 
			\includegraphics[trim= 2 2 2 2, clip,width=0.47\textwidth]{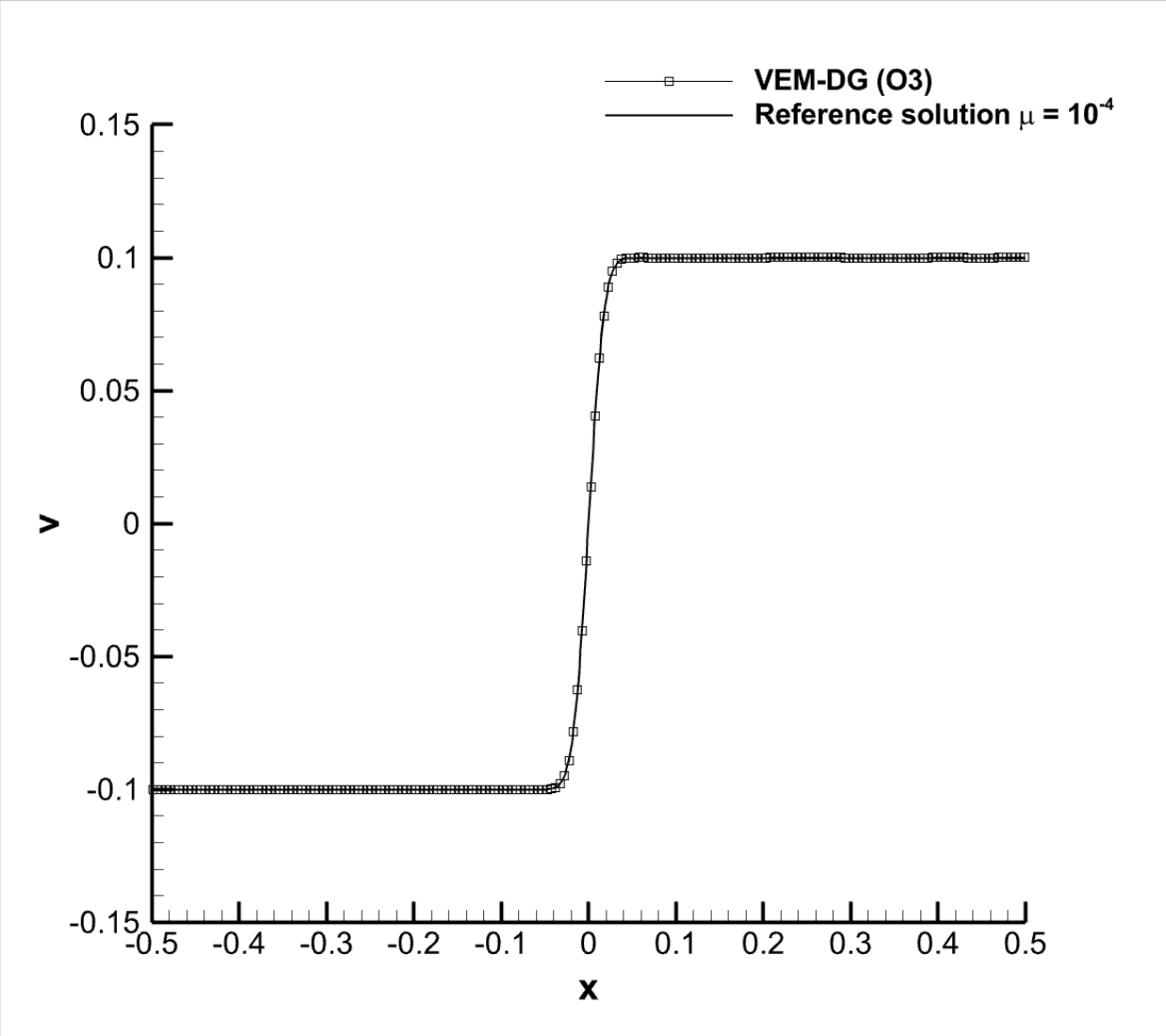} \\
		\end{tabular} 
		\caption{First problem of Stokes at time $t_f=1$. Third order numerical results for the vertical component of the velocity obtained with the VEM-DG scheme and compared against the reference solution by extracting a one-dimensional cut of 200 equidistant points along the $x-$direction at $y=0$. Viscosity $\mu=10^{-3}$ (left) and $\mu=10^{-4}$ (right).}
		\label{fig.Stokes}
	\end{center}
\end{figure}

The first problem of Stokes deals with the time-evolution of an infinite incompressible shear layer with a flow dominated by viscous effects, and for this problem an exact analytical solution of the unsteady Navier-Stokes equations is available \cite{Schlichting}. The computational domain is the channel $\Omega=[-0.5,0.5]\times[-0.05,0.05]$, which is discretized with a coarse grid made of $N_P=358$ cells. The initial condition reads
\begin{equation} 
	\label{eqn.StokesExact}
	\rho = 1, \quad u=0, \quad v=\left\{ \begin{array}{rl}
		v_0 & x \leq 0 \\ -v_0 & x > 0 
	\end{array} \right. , \quad p=\frac{1}{\gamma}, \qquad v_0 = 0.1.
\end{equation}
This initial condition is used to enforce boundary conditions in $x-$direction, while periodic boundaries are prescribed in $y-$direction. The flow is characterized by a low Mach number $\M=10^{-1}$ and heat conduction is neglected ($\kappa=0$). The final time is chosen to be $t_f=1$ and we set $\CFL=0.5$. The simulation is carried out with the third order VEM-DG schemes for two different values of viscosity, namely $\mu=10^{-3}$ and $\mu=10^{-4}$. The comparison between the reference solution and the numerical results is presented in Figure \ref{fig.Stokes}, showing an excellent matching.

\subsection{Explosion problem}
The explosion problem is a benchmark for the compressible Euler equations, hence we consider again an inviscid fluid with no heat conduction. The computational domain is the square $\Omega=[-1;1]^2$ with transmissive boundaries that is paved with a Voronoi mesh of characteristic size $h_{\Omega}=1/128$. The initial condition consists in two different states, separated by the circle of radius $R=0.5$: 
\begin{equation}
	\mathbf{P}(\xx,0) = \left\{ \begin{array}{clcc} \mathbf{P}_i = & (1.0, 0.0, 0.0, 1.0)    & \textnormal{ if } & r \leq R, \\ 
		\mathbf{P}_o = & (0.125, 0.0, 0.0, 0.1)  & \textnormal{ if } & r > R,        
	\end{array}  \right. 
\end{equation}
where $\mathbf{P}_i$ and $\mathbf{P}_o$ denote the inner and the outer state, respectively, while $r=\sqrt{x^2+y^2}$ is the generic radial position. As done in \cite{TavelliCNS}, the initial discontinuity in the $L_2$ projection of the DG solution is slightly smoothed, so that nonphysical oscillations at the initial time are eliminated:
\begin{equation}
	\mathbf{P}(\xx,0) = \frac{1}{2} \left(\mathbf{P}_o+\mathbf{P}_i\right) + \frac{1}{2} \left(\mathbf{P}_o-\mathbf{P}_i\right) \textnormal{erf} \left( \frac{r-R}{\alpha_0} \right), \qquad \alpha_0=10^{-2}.
\end{equation}

\begin{figure}[!htbp]
	\begin{center}
		\begin{tabular}{cc}  
			\includegraphics[trim= 5 5 5 5, clip, width=0.47\textwidth]{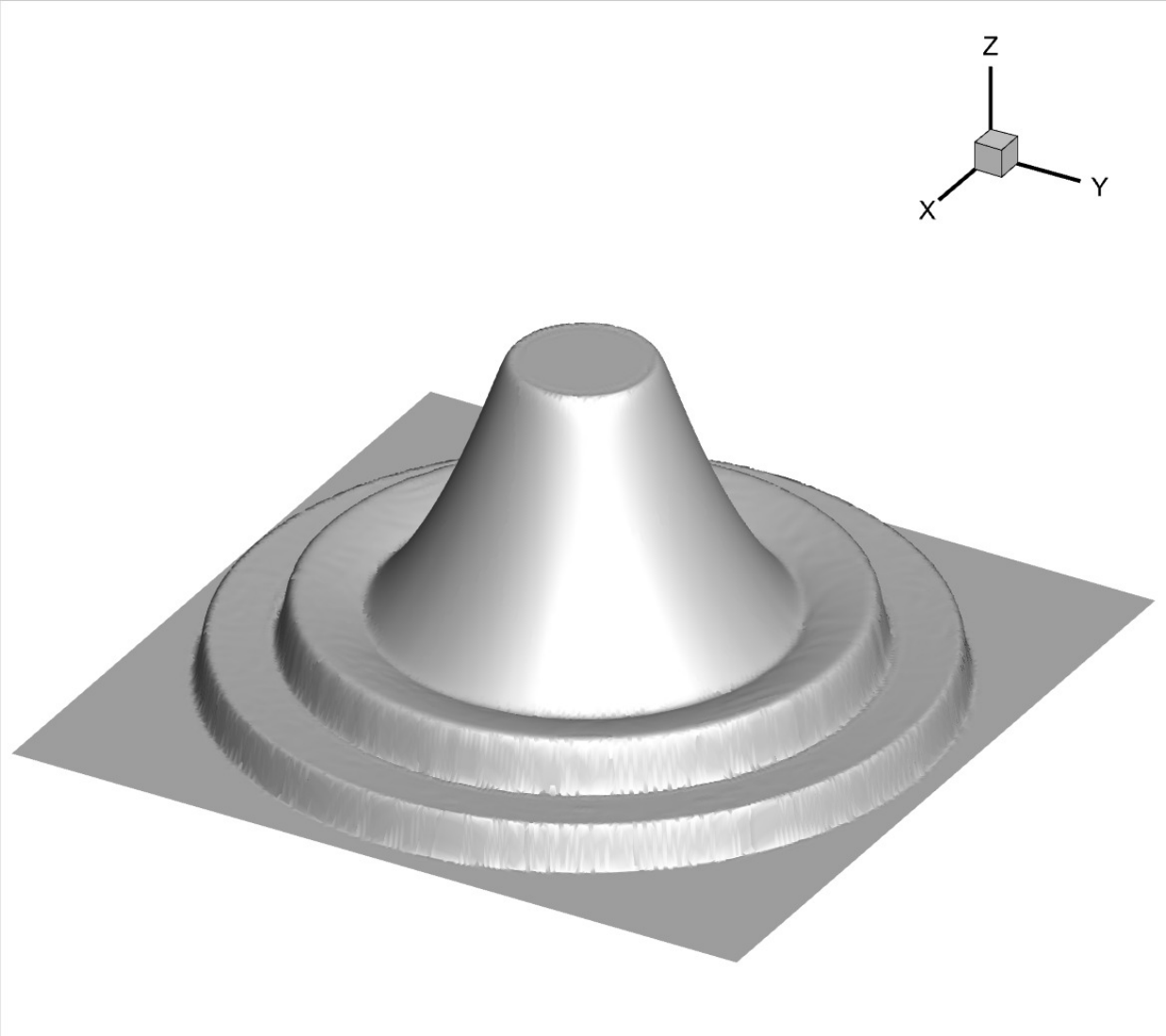}&
			\includegraphics[width=0.47\textwidth]{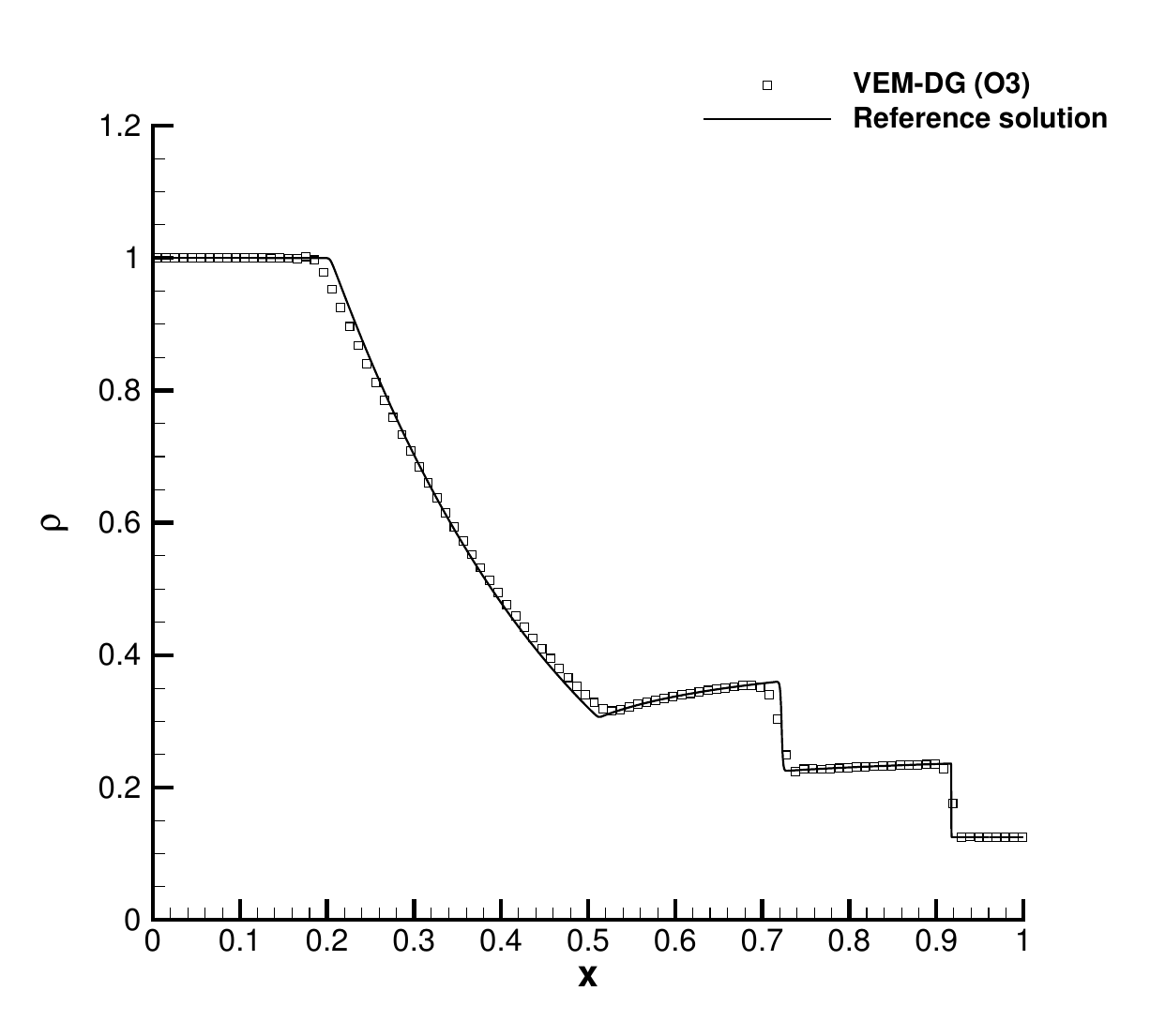} \\				
			\includegraphics[width=0.47\textwidth]{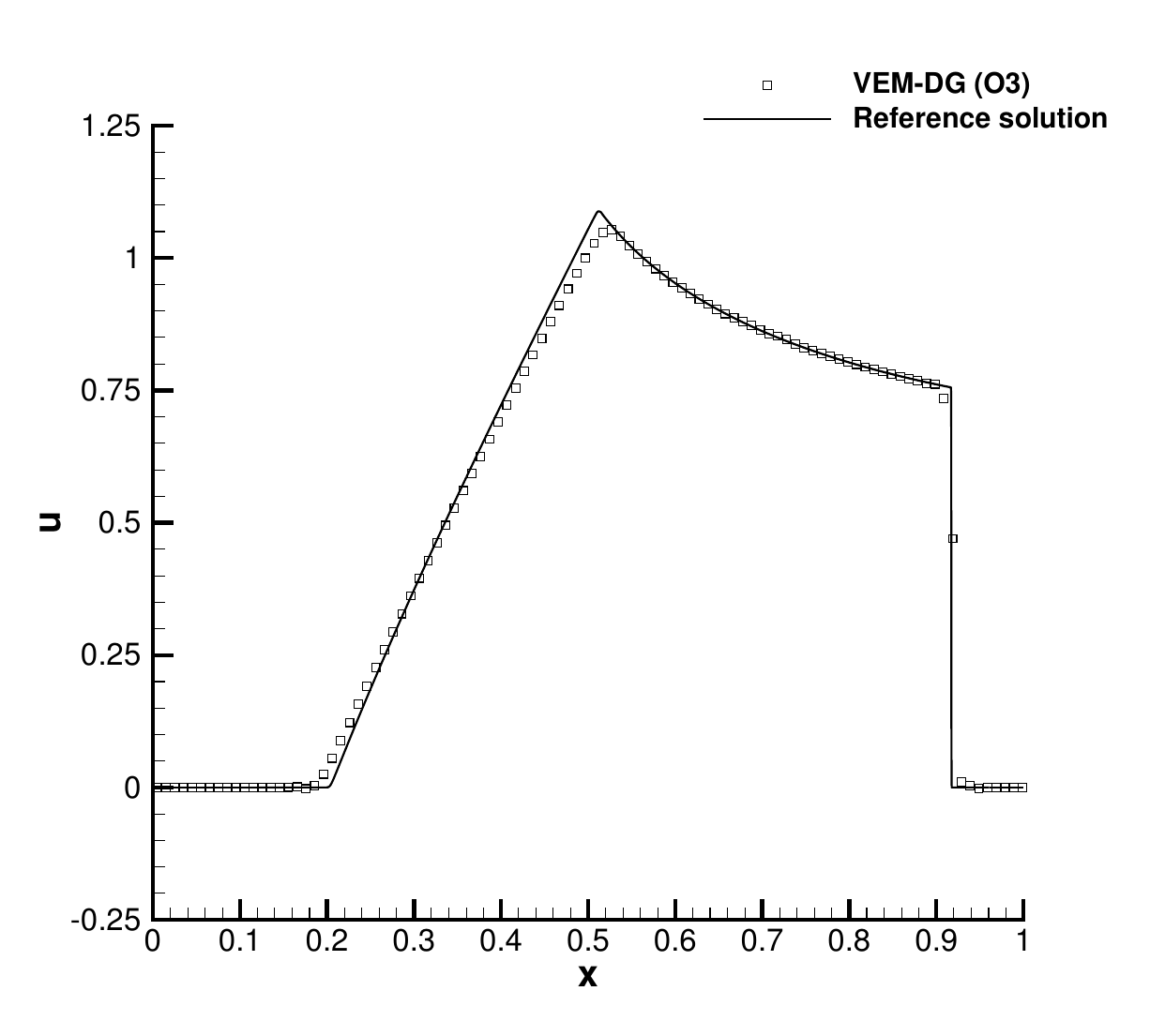} & 
			\includegraphics[width=0.47\textwidth]{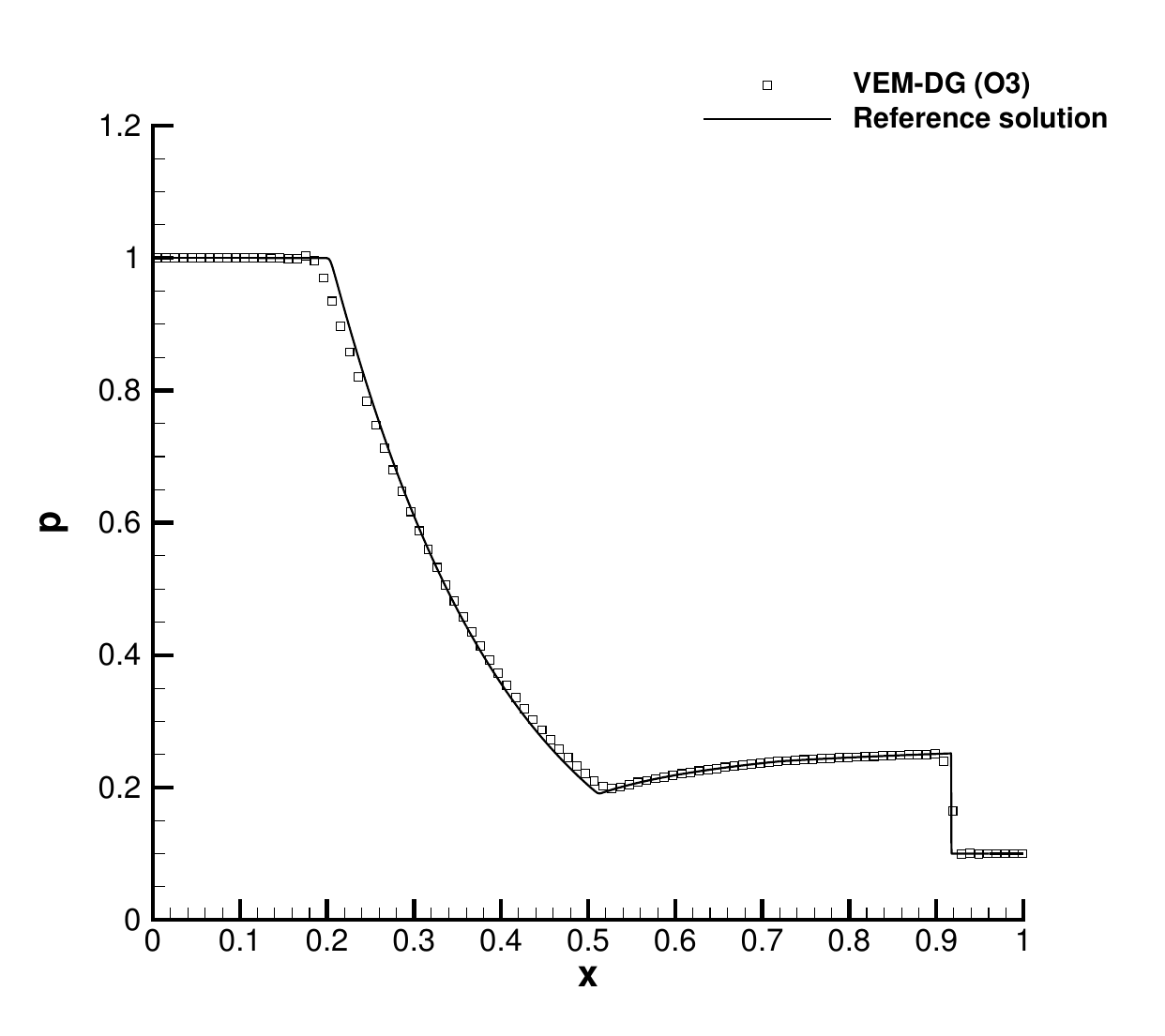} \\
		\end{tabular} 
		\caption{Explosion problem at time $t_f=0.25$. Third order numerical results with VEM-DG scheme for density, horizontal velocity and pressure compared against the reference solution
			extracted with a one-dimensional cut of 200 equidistant points along the $x-$direction at $y=0$. Three-dimensional view of the density distribution is shown in the top left panel.}
		\label{fig.EP2D}
	\end{center}
\end{figure}

The solution involves three different types of waves, namely an outward traveling shock front, an inward moving rarefaction fan and a contact wave in between. We set $\CFL=0.5$ and the final time $t_f=0.25$. Figure \ref{fig.EP2D} depicts the numerical results obtained with the third order version of the novel VEM-DG schemes, which are compared against the reference solution, computed as an equivalent one-dimensional problem in radial direction $r$ with a geometric source term (see \cite{ToroBook}). An excellent agreement can be appreciated, noticing that the numerical solution preserves its symmetry even on general unstructured meshes, as evident from the three-dimensional view of the density distribution plot in Figure \ref{fig.EP2D}.

\begin{figure}[!htbp]
	\begin{center}
		\begin{tabular}{cc}   
			\includegraphics[trim= 2 2 2 2, clip,width=0.47\textwidth]{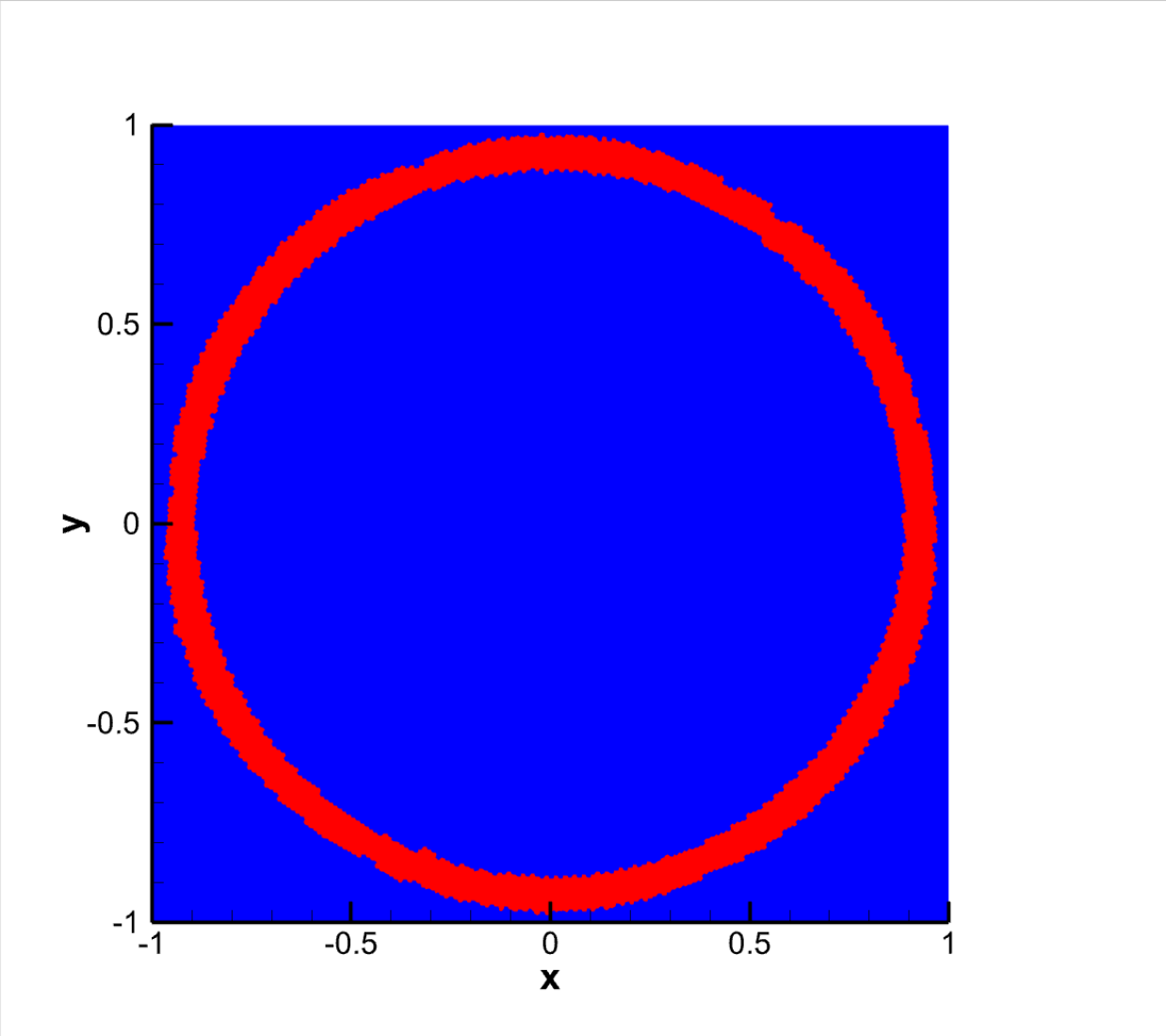} & 
			\includegraphics[width=0.47\textwidth]{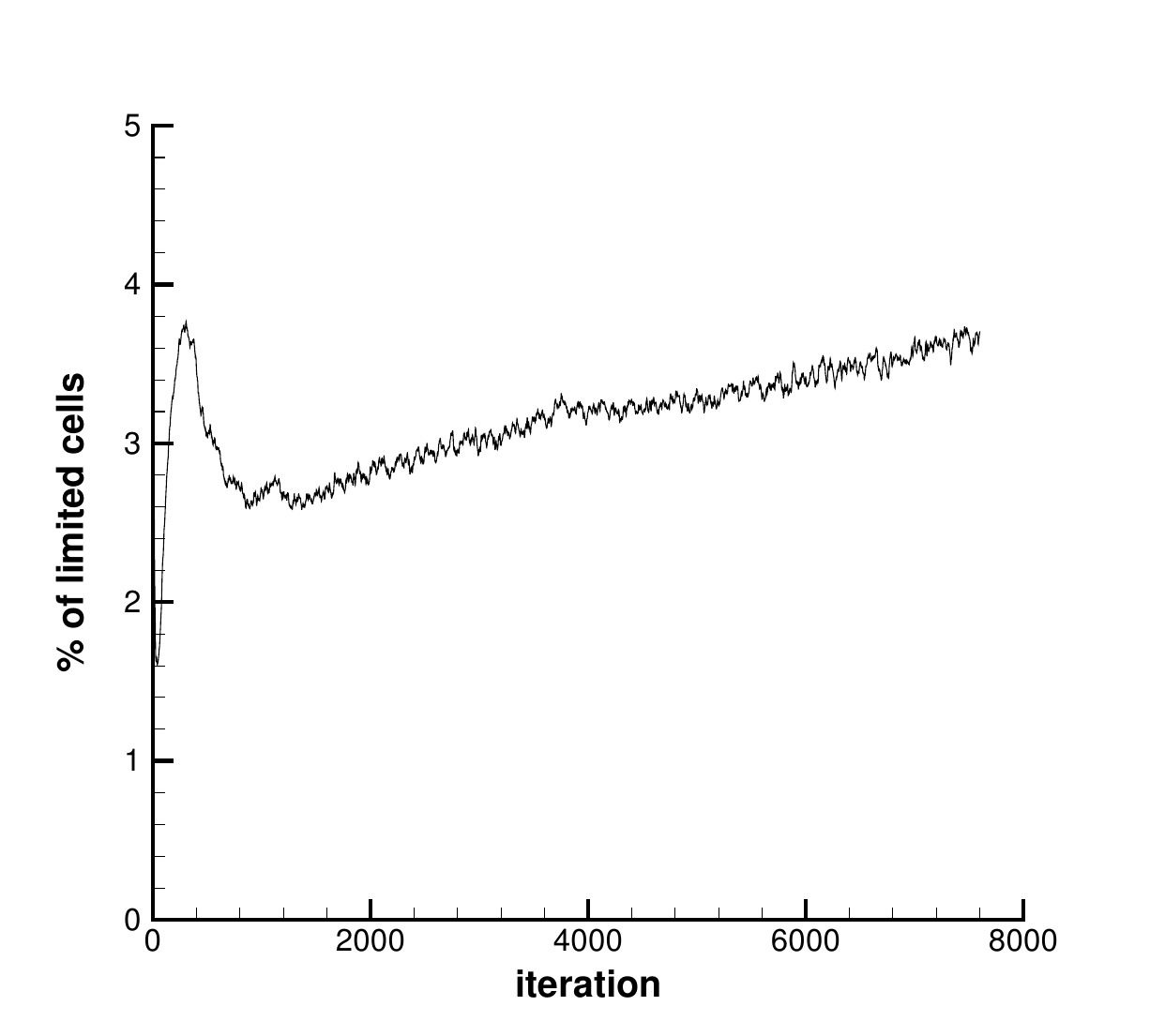} \\ 
		\end{tabular} 
		\caption{Explosion problem with third order VEM-DG scheme. Troubled cell map (left) with the limited cells highlighted in red and the unlimited cells colored in blue; percentage of limited cells at each time step (right).}
		\label{fig.EP2Dlim}
	\end{center}
\end{figure}

Since the solution exhibits a shock wave, the artificial viscosity limiter is activated. Figure \ref{fig.EP2Dlim} shows the limited cells, which are perfectly detected only across the shock wave. The time evolution of the total number of limited cells is also reported. We notice that very few cells need to be supplemented with artificial viscosity, mainly concentrated along the shock profile and involving less than $4\%$ of the total number of elements.

\subsection{Viscous shock profile}
\begin{figure}[!htbp]
	\begin{center}
		\begin{tabular}{cc} 
			\multicolumn{2}{c}{\includegraphics[trim=2 0 2 2,clip,width=0.9\textwidth]{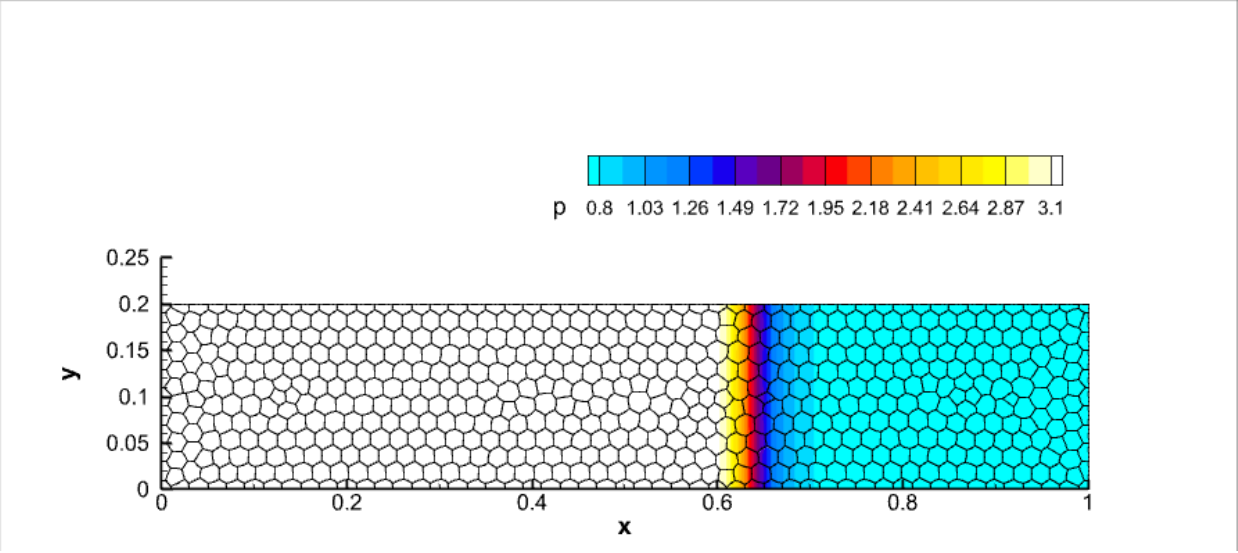}} \\
			\includegraphics[width=0.47\textwidth]{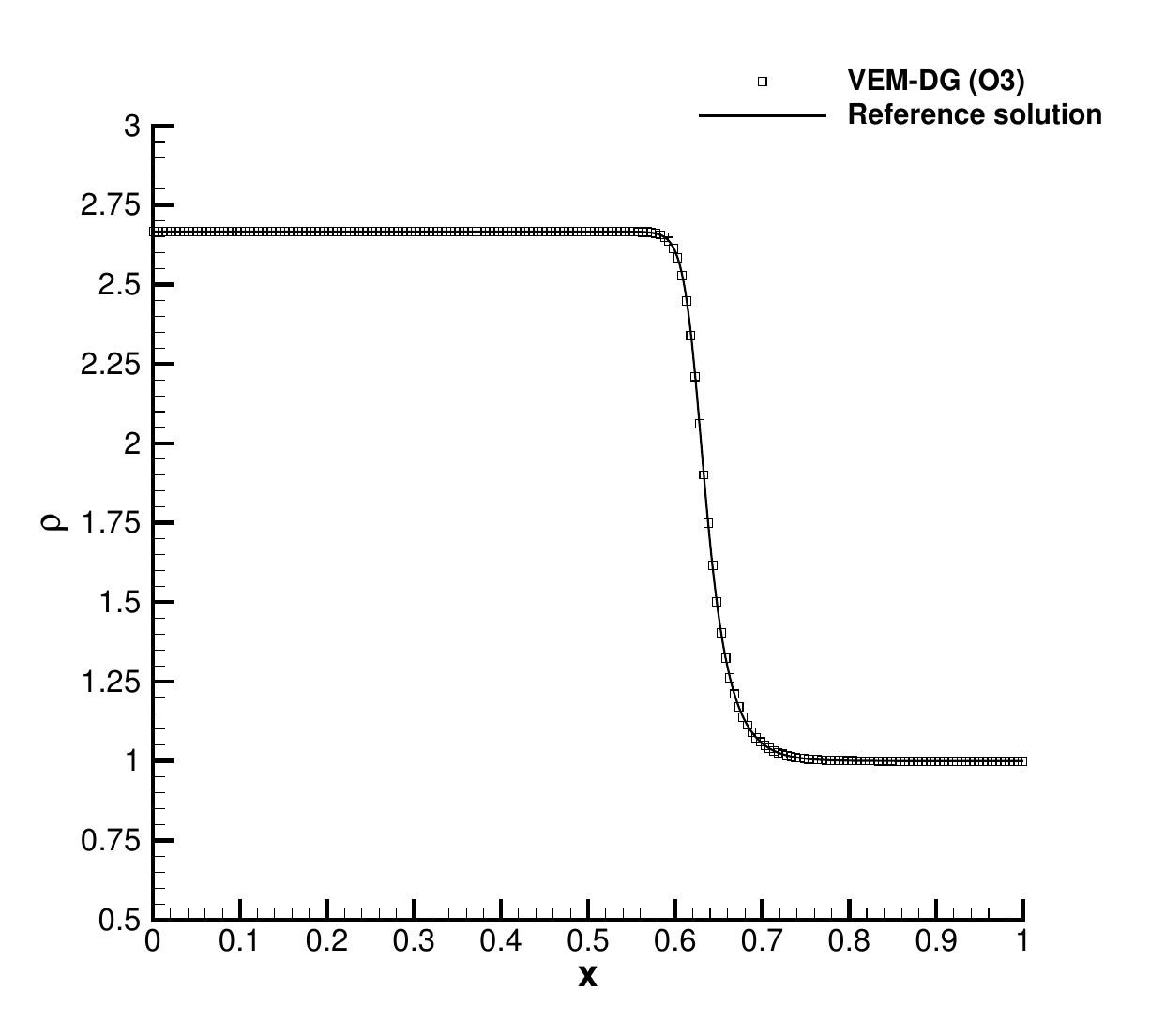} & 
			\includegraphics[width=0.47\textwidth]{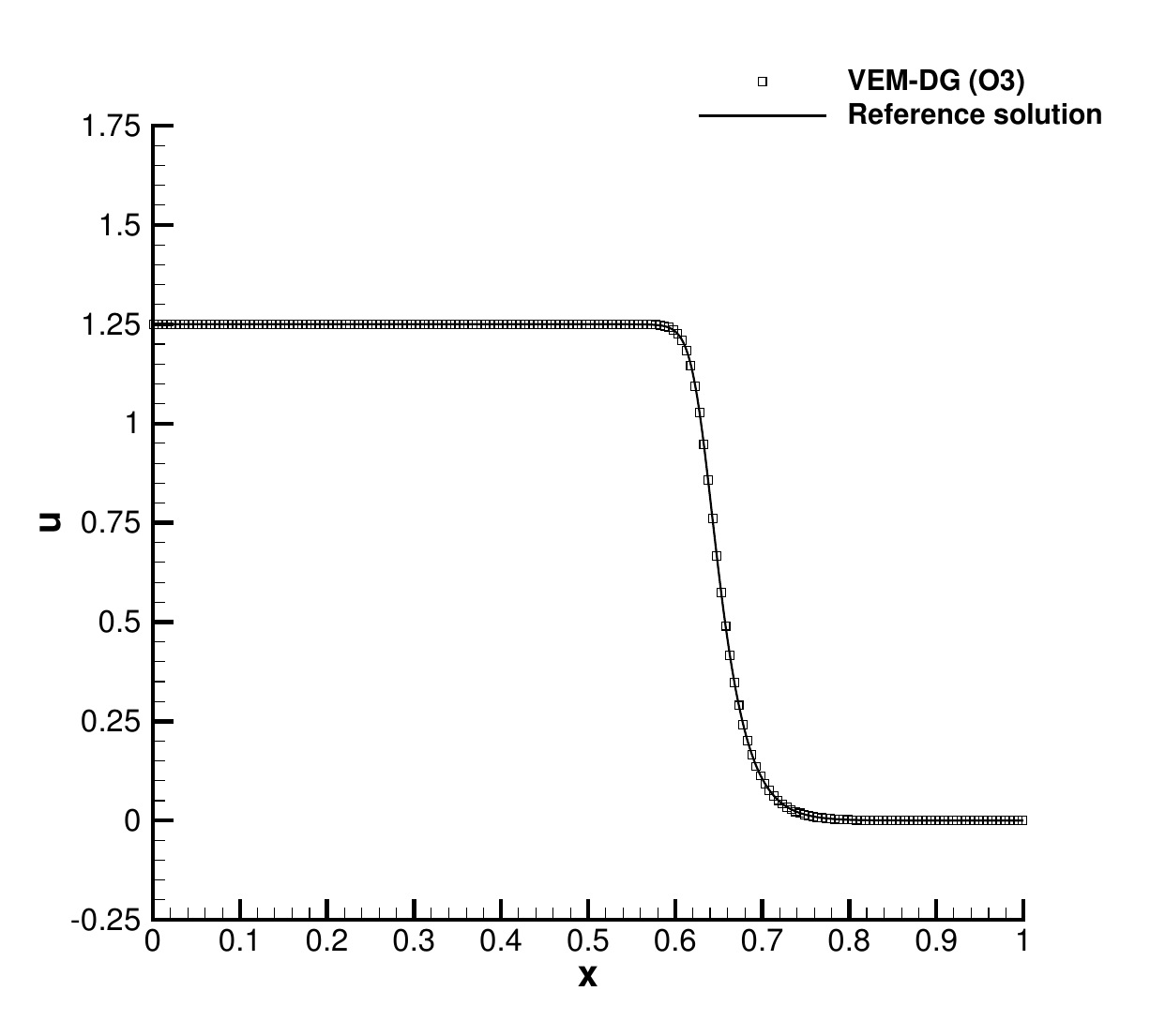} \\
			\includegraphics[width=0.47\textwidth]{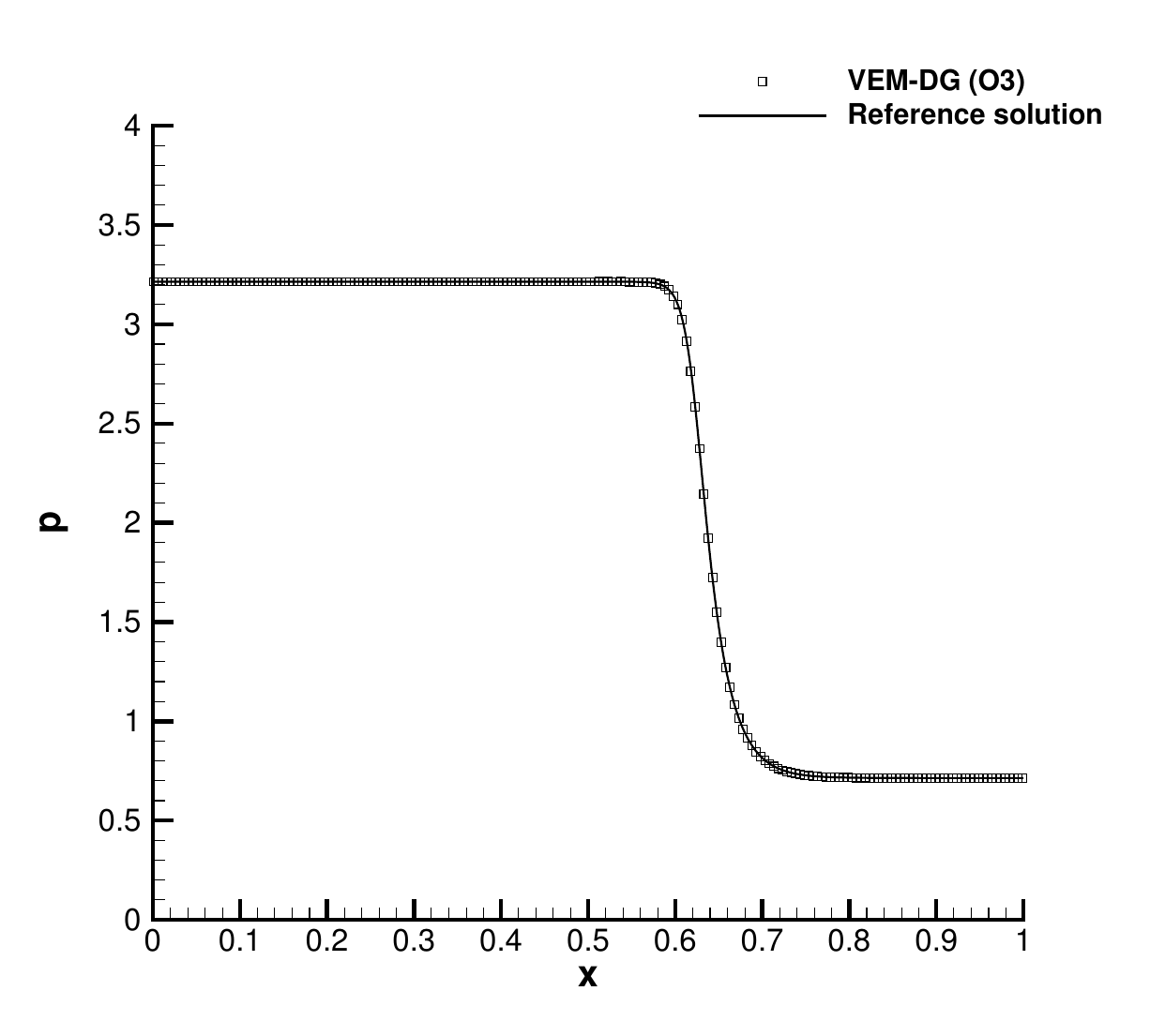} & 
			\includegraphics[width=0.47\textwidth]{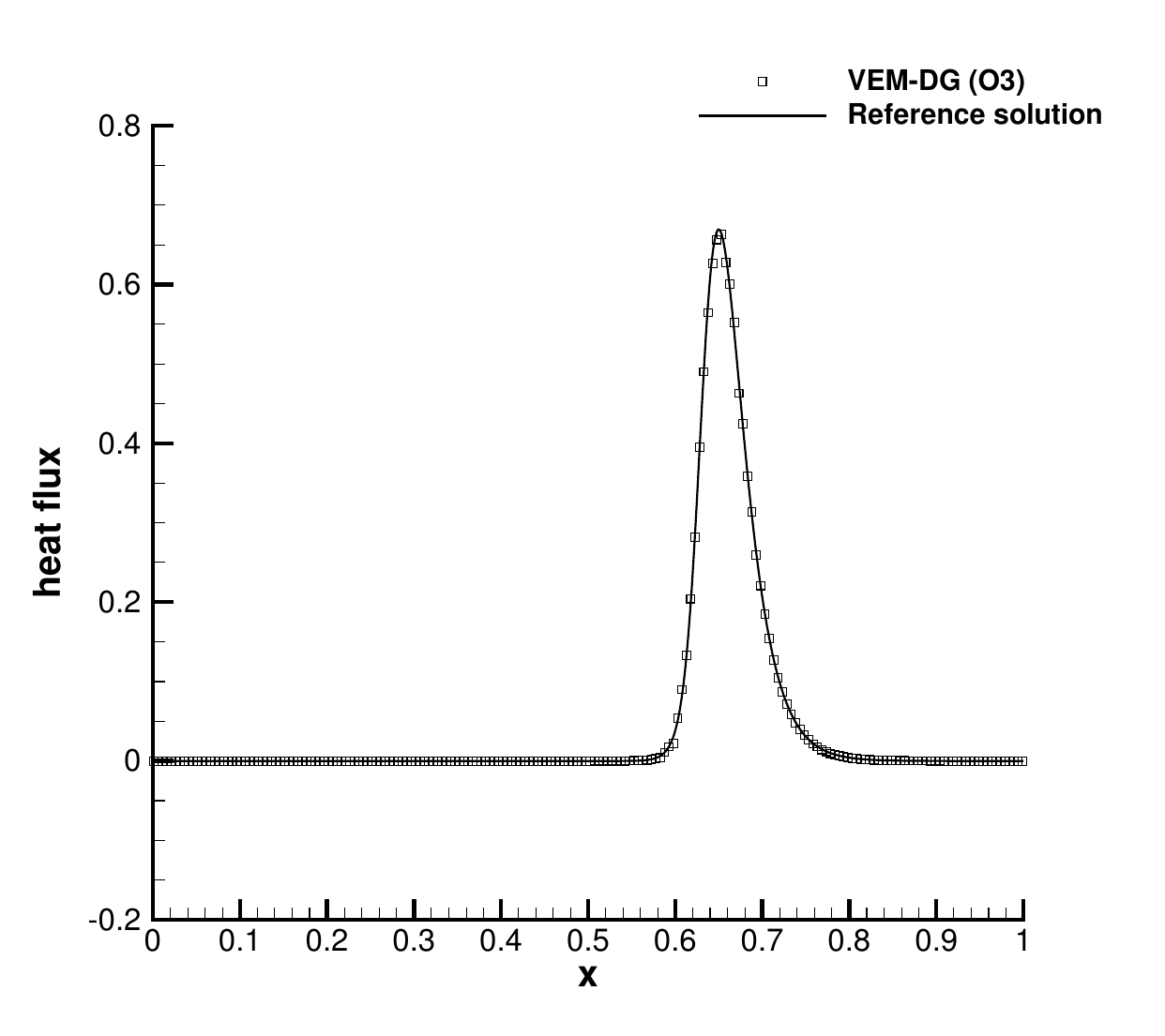} \\
		\end{tabular} 
		\caption{Viscous shock profile with shock Mach number $\mathrm{M_s}=2$ and Prandtl number $\Pr=0.75$ at time $t_f=0.2$. Top panel: Voronoi tessellation and pressure contours. Middle and bottom panel: third order numerical solution with VEM-DG scheme compared against the reference solution for density, horizontal velocity, pressure and heat flux (from middle left to bottom right panel). We show a one-dimensional cut of 200 equidistant points along the $x-$direction at $y=0.1$.}
		\label{fig.ViscousShock}
	\end{center}
\end{figure}

We consider an analytical solution of the compressible Navier-Stokes equation given by an isolated viscous shock wave which is traveling into a medium at rest with a shock Mach number of $\mathrm{M_s} > 1$. The exact solution is derived in \cite{Becker1923}, where the compressible Navier-Stokes equations are solved for the special case of a stationary shock wave at Prandtl number $\Pr= 0.75$ with constant viscosity. This test case is particularly interesting because it allows all terms contained in the Navier-Stokes system to be properly verified, since advection, thermal conduction and viscous stresses are all present. The computational domain is given by $\Omega=[0;1]\times [0;0.2]$ which is discretized by a coarse grid made of a total number of unstructured polygons $N_P=1120$, shown in the top panel of Figure \ref{fig.ViscousShock}. At $x=0$, inflow boundary conditions are imposed, while outflow conditions are considered at $x=1$. Periodic boundaries are set in $y-$direction. The details for the setup of this test case can be found in \cite{ADERNSE,ADERAFEDG}, thus we briefly recall that the initial condition is given by a shock wave centered at $x=0.25$, traveling at Mach $\mathrm{M_s}=2$ from left to right, with a Reynolds number $\Rey=100$. The viscosity coefficient is $\mu=2 \cdot 10^{-2}$ and the final time of the simulation is $t_f=0.2$, with the shock front located at $x=0.65$. We use a $\CFL$ number of $\CFL=0.5$ and the third order VEM-DG schemes. The results are depicted in Figure \ref{fig.ViscousShock}, where a comparison against the analytical solution is provided, retrieving an excellent matching. We compare the exact solution and the numerical solution for density, horizontal velocity component, pressure and heat flux $q_x=\kappa \, \frac{\partial T}{\partial x}$, with $T=p/(R\rho)$ being the temperature, as given in the thermal EOS \eqref{EOS}.  We also notice an excellent symmetry preservation of the solution along the $y-$direction, despite the use of unstructured Voronoi meshes, where in general the edges of the control volumes are not aligned with the main flow field in the $x-$direction. 

\subsection{2D Taylor-Green vortex}
\begin{figure}[!htbp]
	\begin{center}
		\begin{tabular}{cc} 
			\includegraphics[trim=2 2 2 2,clip,width=0.47\textwidth]{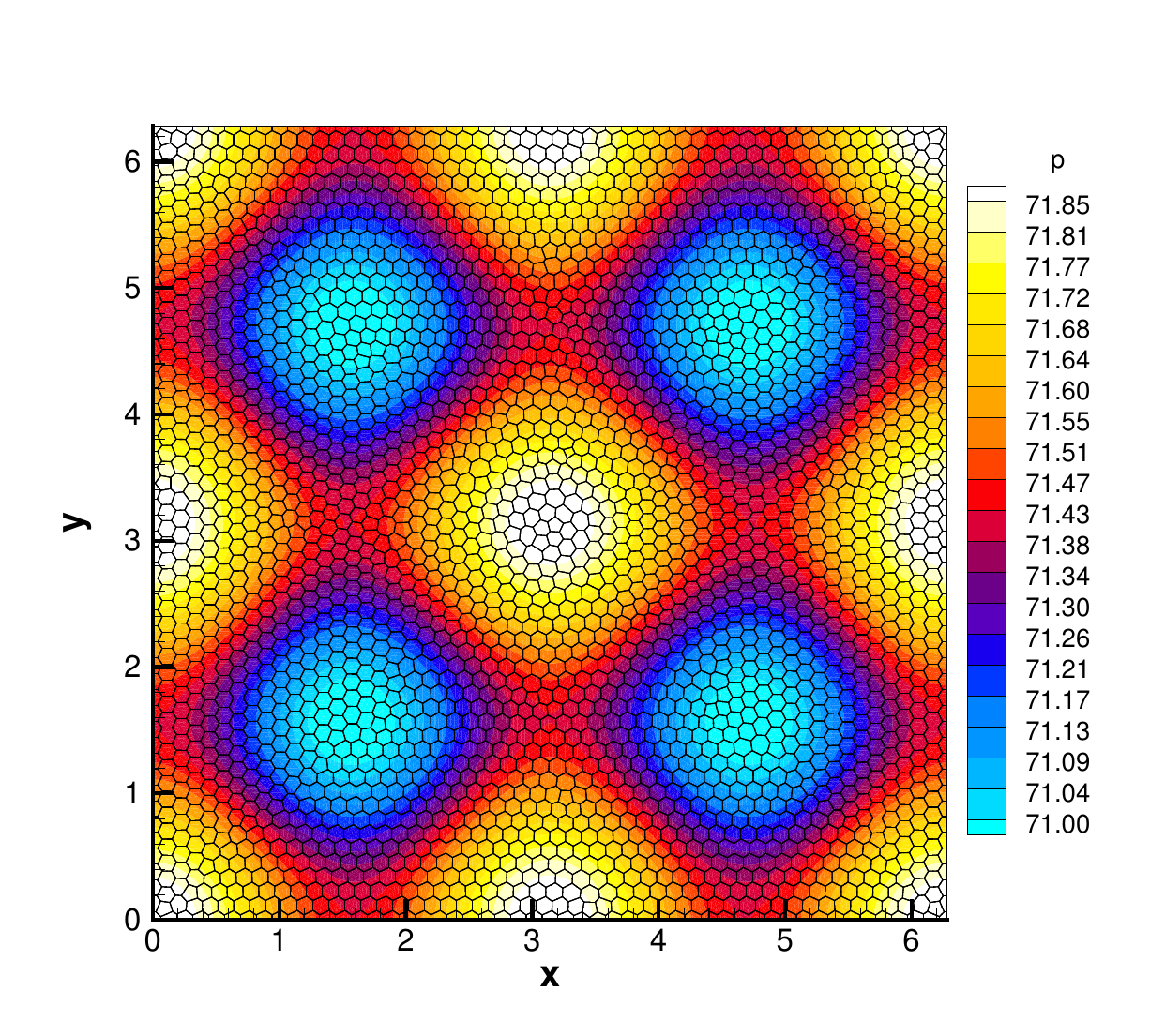} & 
			\includegraphics[trim=2 2 2 2,clip,width=0.47\textwidth]{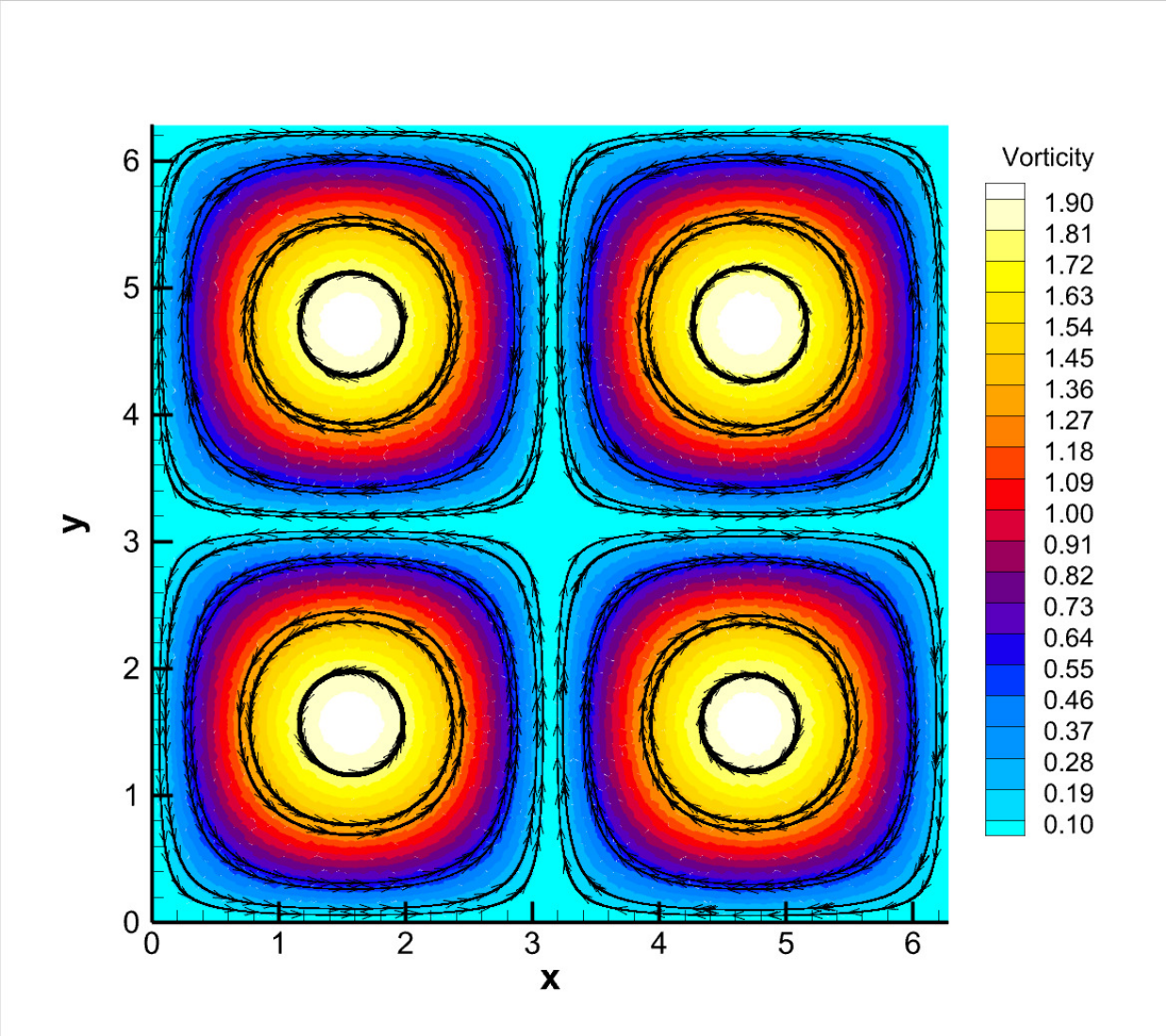} \\
			\includegraphics[width=0.47\textwidth]{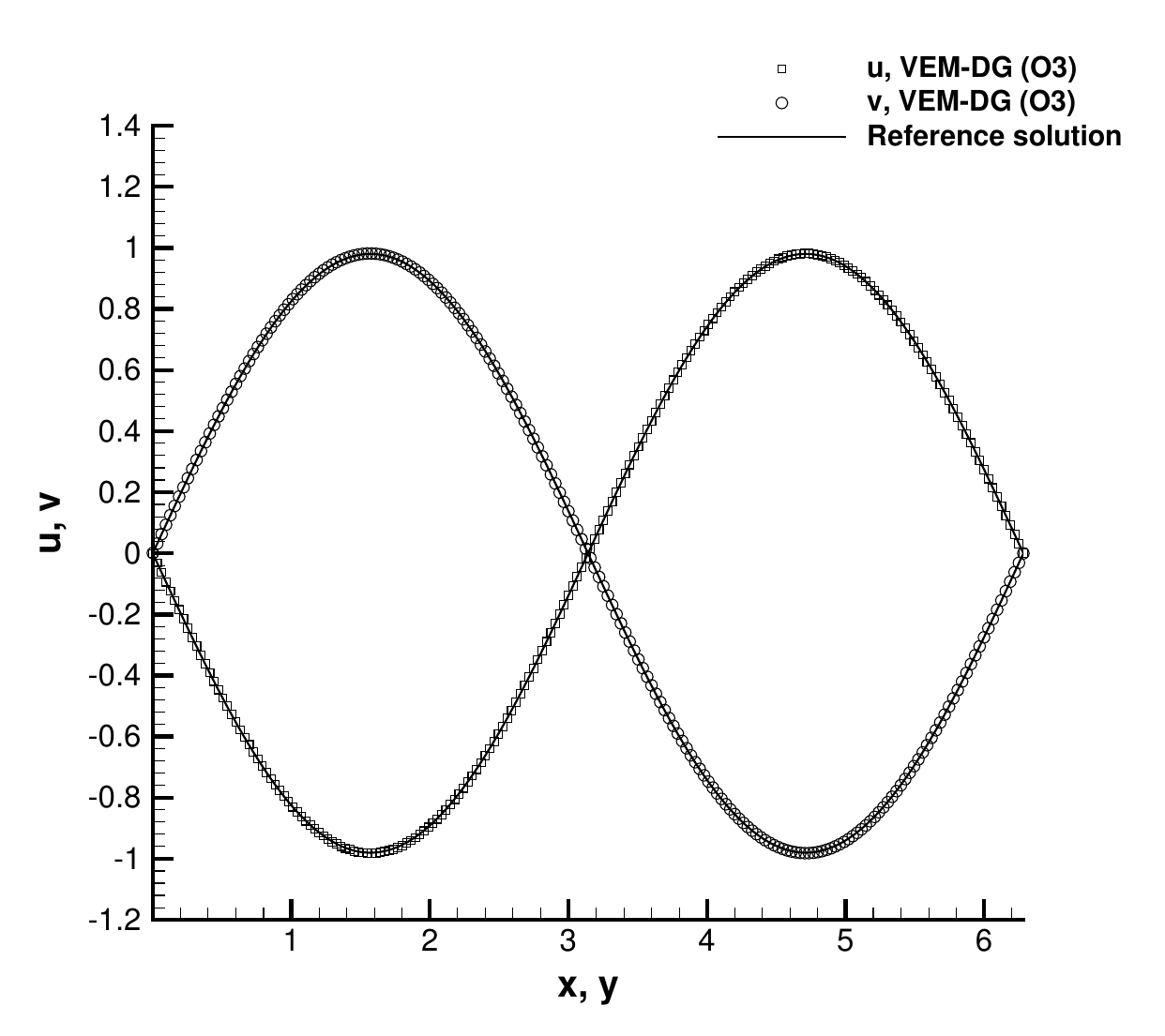} & 
			\includegraphics[width=0.47\textwidth]{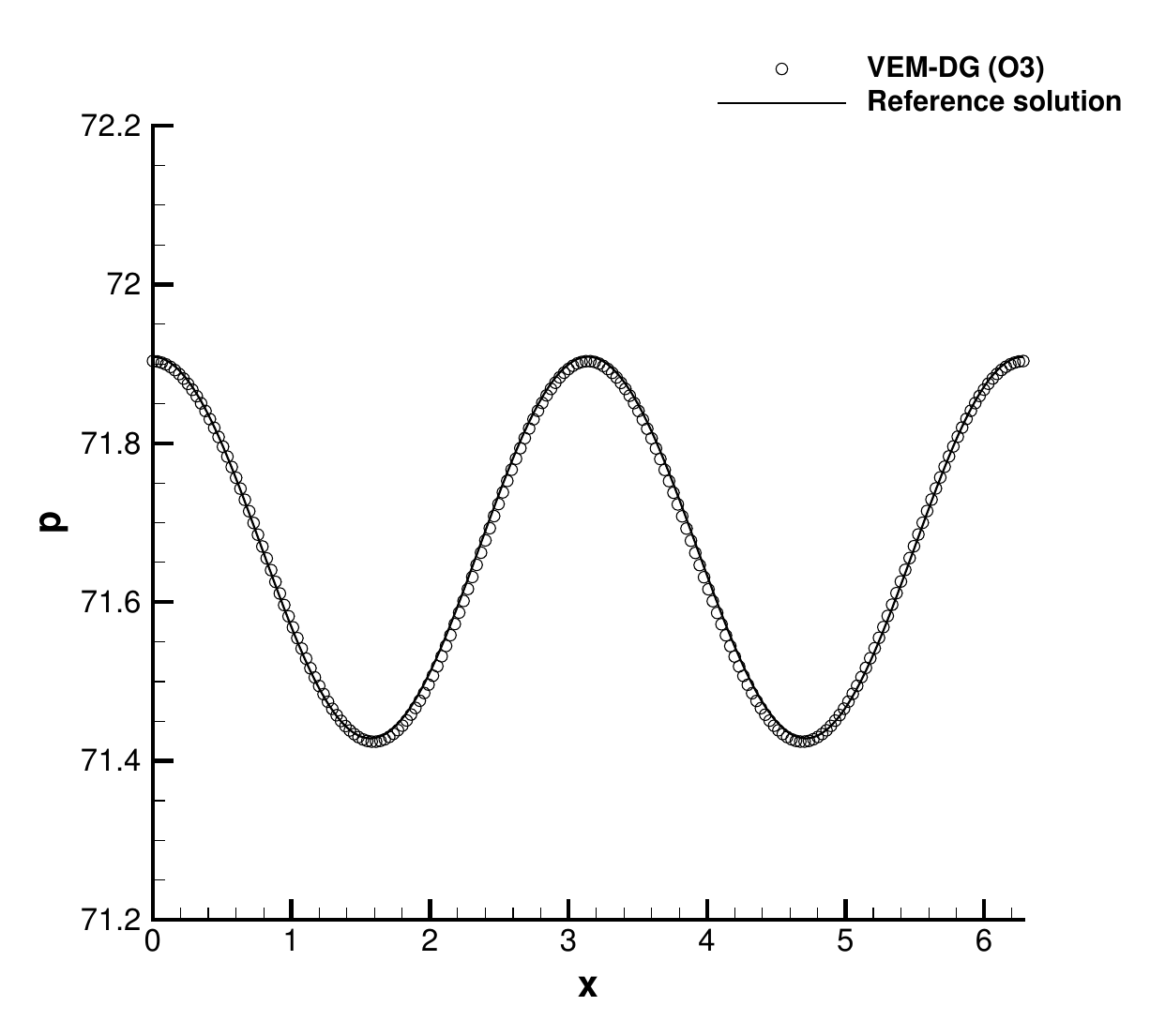} \\
		\end{tabular} 
		\caption{2D Taylor-Green vortex at time $t_f=1$ with viscosity $\mu=10^{-2}$. Exact solution of the Navier-Stokes equations and third order numerical solution with VEM-DG scheme. Top: mesh configuration with pressure distribution (left) and vorticity magnitude with stream-traces (right). Bottom: one-dimensional cut of 200 equidistant points along the $x$-axis and the $y-$axis for the velocity components $u$ and $v$ (left) and for the pressure $p$ (right).}
		\label{fig.TGV2D}
	\end{center}
\end{figure}

The Taylor-Green vortex problem is a benchmark for the incompressible Navier-Stokes equations, and an exact solution is known in two space dimensions:
\begin{eqnarray}
	u(\xx,t)&=&\phantom{-}\sin(x)\cos(y) \, e^{-2\nu t},  \nonumber \\
	v(\xx,t)&=&-\cos(x)\sin(y) \, e^{-2\nu t}, \nonumber \\
	p(\xx,t)&=& C + \frac{1}{4}(\cos(2x)+\cos(2y)) \, e^{-4\nu t},
	\label{eq:TG_ini}
\end{eqnarray}
with $\nu=\mu/\rho$ being the kinematic viscosity of the fluid. For modeling the low Mach regime of the compressible Navier-Stokes equations, we set $C=100/\gamma$ as additive constant for the pressure field. Heat conduction is neglected, the viscosity coefficient is chosen to be $\mu=10^{-2}$, and the initial density is $\rho(\xx,0)=1$. The computational domain is defined by $\Omega=[0;2\pi]^2$ with periodic boundaries everywhere. The mesh counts a total number of $N_P=2916$ Voronoi cells. The final time of the simulation is $t_f=1$ and we set $\text{CFL}=0.5$. Figure \ref{fig.TGV2D} shows the numerical results obtained running the third order version of the new VEM-DG schemes, where a comparison against the exact solution is depicted. Here, one can appreciate a very good agreement between the VEM-DG method in the low Mach number regime and the exact solution of the incompressible Navier-Stokes equations, both for velocity and pressure profiles. The stream-traces of the velocity field and the pressure distribution are also plotted.

\subsection{Compressible mixing layer} \label{ssec.MixLayer}
The study of an unsteady compressible mixing layer, as proposed in \cite{Colonius}, is widely used to test numerical methods for compressible viscous flows. The two-dimensional computational domain is the rectangular box $\Omega=[-200,200] \times [-50,50]$, discretized with a total number of $\NP=15723$ Voronoi cells. At the left side of the domain ($x=-200$), a time-dependent inflow boundary condition is prescribed in terms of a function $\delta(y,t)$:
\begin{eqnarray}
	\rho(y,t) &=& \rho_0 + 0.05 \, \delta(y,t), \nonumber \\
	\vv(y,t) &=& \vv_0 + \left( \begin{array}{c}
		1.0 \\ 0.6
	\end{array} \right)\, \delta(y,t), \\
	\quad p(y,t) &=& p_0 + 0.2 \, \delta(y,t), \nonumber 
\end{eqnarray}
with the background state 
\begin{equation}
\rho_0= 1,  \quad \vv_0 = \left( \begin{array}{c}
			\frac{1}{8} \tanh(2y) + \frac{3}{8} \\ 0
		\end{array} \right), \quad p_0 = \frac{1}{\gamma}, 
	\end{equation}  
and the perturbation 
{\small
	\begin{eqnarray}
		\delta(y,t) &=& -10^{-3} \exp(-0.25 y^2) \cdot \nonumber \\
		&\phantom{=}&\left[ \cos(\omega t) + \cos\left(\frac{1}{2}\omega t -0.028\right) + \cos\left(\frac{1}{4}\omega t +0.141\right) + \cos\left(\frac{1}{8}\omega t +0.391\right) \right], 
\end{eqnarray}}
with the fundamental frequency of the mixing layer $\omega = 0.3147876$. Initially, the background state is assigned to the flow, hence
\begin{equation*}
	\rho(\xx,0)=\rho_0,  \quad \vv(\xx,0)=\vv_0, \quad p(\xx,0)=p_0.
\end{equation*}

\begin{figure}[!htbp]
		\begin{center}
			\begin{tabular}{c} 
				\includegraphics[trim= 5 5 5 5, clip, width=0.9\textwidth]{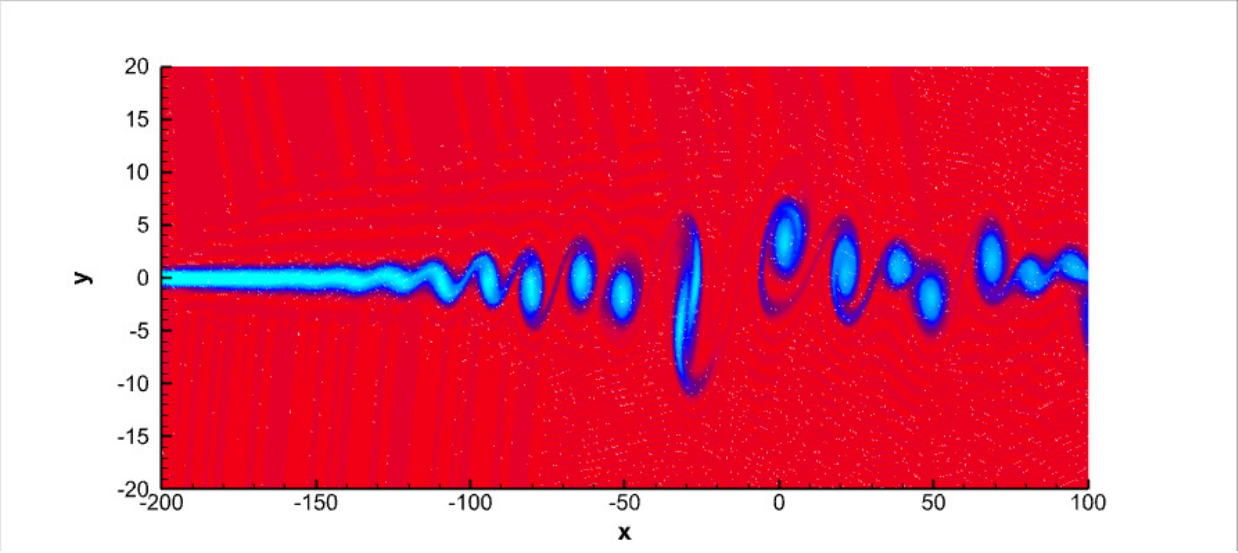}\\				
				\includegraphics[trim= 5 5 5 5, clip, width=0.9\textwidth]{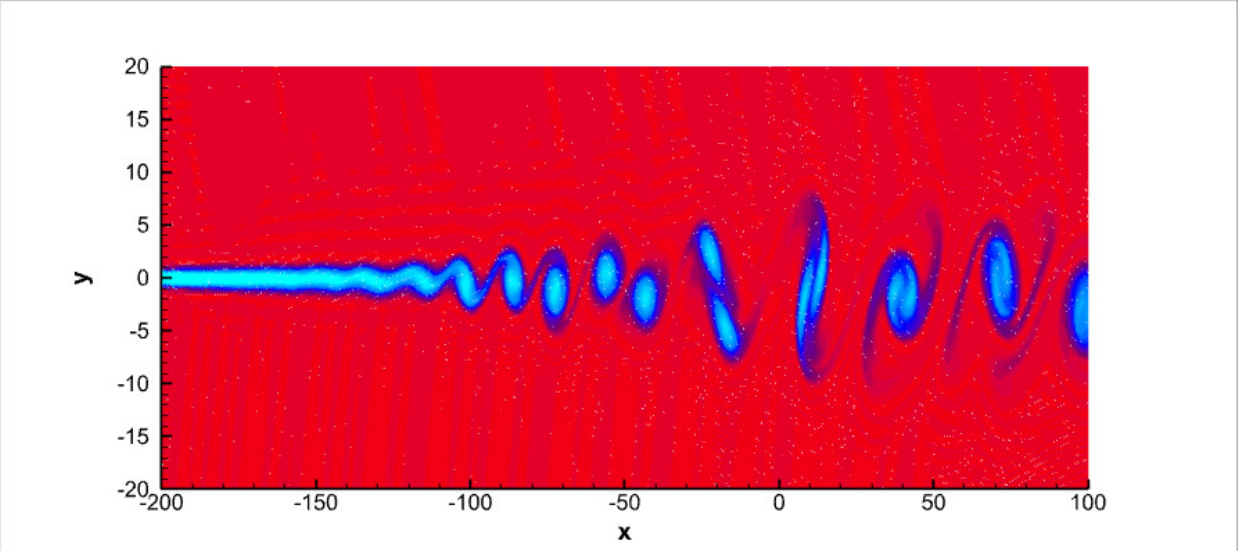}\\
				\includegraphics[trim= 5 5 5 5, clip, width=0.9\textwidth]{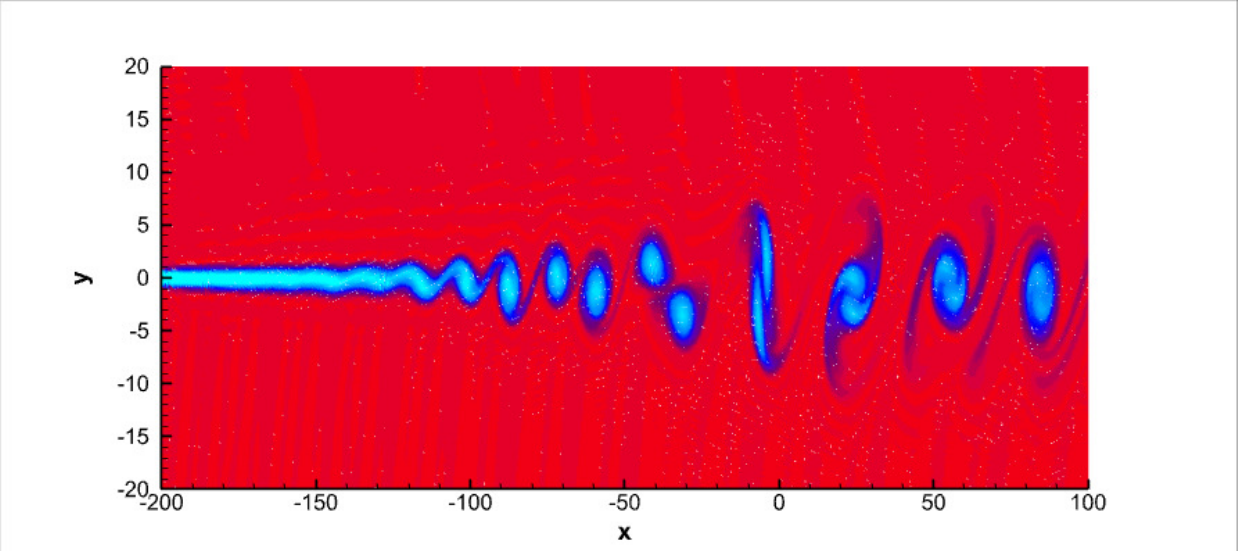}\\				
			\end{tabular} 
			\caption{Compressible mixing layer at time $t=500$, $t=1000$ and $t=1596.8$ (from top to bottom row). Third order numerical results with Runge-Kutta VEM-DG for $z-$vorticity. 51 contour levels in the range $[-0.12,0.12]$ have been used for plotting the vorticity distribution on the sub-domain $[-200,100]\times[-20,20]$.}
			\label{fig.MixLayer}
		\end{center}
	\end{figure}

On the right boundary we set outflow conditions, while the free stream velocities are imposed in the $y-$direction, thus we set $u_{+\infty}=0.5$ and $u_{-\infty}=0.25$ for $y \to + \infty$ and $y \to -\infty$, respectively. Heat conduction is neglected, and the viscosity coefficient is set to be $\mu=10^{-3}$. For this final test, the simulation is run with the third order VEM-DG scheme with Runge-Kutta time stepping up to the final time $t_f=1596.8$. We show the vorticity of the flow field at three different output times in Figure \ref{fig.MixLayer}. The resolution of the numerical results permits to appreciate the vortical structures generated by the perturbation at the inflow boundary, demonstrating the capability of the novel methods to capture complex vortical patterns in the flow field. We also notice that these results are qualitatively in good agreement with those obtained with the method presented in \cite{ADERAFEDG}.
	
\section{Conclusions} \label{sec.concl}
In this article, a new family of discontinuous Galerkin finite element
methods has been introduced for the solution of nonlinear systems of hyperbolic PDEs on unstructured meshes composed of Voronoi cells. The numerical solution is approximated by means of novel Virtual Element basis functions, borrowing the $L_2$ projection operators defined in the Virtual Element Method \cite{vem2}. The resulting basis leads to a mixed nodal/modal approach that permits to represent the discrete solution with arbitrary order of accuracy within each control volume. The time marching is carried out relying on the ADER strategy, which yields a fully discrete one-step DG scheme in space and time based on a variational formulation of the governing PDEs. The global solution space results to be discontinuous, thus the novel Virtual Element basis functions are said to be nonconforming, in accordance with the definition of nonconforming VEM spaces in \cite{VEM_nc_org}. The space-time basis functions are used to construct mass and stiffness matrices that are stabilized relying on an extension of the dof--dof VEM stabilization proposed here for the first time to handle space-time polynomial spaces. An orthogonalization of the basis functions is also forwarded, following the ideas outlined in \cite{Berrone2017}, which allows the condition number of the mass and stiffness matrices to be improved. Applications to compressible inviscid and viscous flows demonstrate the accuracy, robustness and capabilities of the novel VEM-DG schemes.

By exploiting the nodal interpolation property of the VEM basis, future work will concern the investigation of quadrature-free VEM-DG schemes, aiming at improving the computational efficiency of the algorithm. Furthermore, we remark that an implicit treatment of the viscous fluxes \cite{BDT_cns} would lead to a remarkable relaxation of the $\CFL$ stability condition on the maximum admissible time step, thus contributing to make the scheme more efficient as well. Finally, we foresee also the employment of the novel nonconforming Virtual Element basis for the construction of structure-preserving differential operators \cite{SPDGdivcurl} to be used in magnetohydrodynamics or elasticity in solid mechanics. Stabilization-free techniques devised in the VEM framework \cite{BERRONE2023108641} will also be part of further developments in the definition of the space-time matrices.

\section*{Acknowledgments}

WB received financial support by Fondazione Cariplo and Fondazione CDP (Italy) under the project No. 2022-1895 and by the Italian Ministry of University and Research (MUR) with the PRIN Project 2022 No. 2022N9BM3N. 
GB has been partially funded by the University of Ferrara under the call “Bando Giovani anno 2022”.
Both the authors are members of GNCS--INdAM (\textit{Gruppo Nazionale per il Calcolo Scientifico} of the Italian \textit{Istituto Nazionale di Alta Matematica}). 
This work was partially carried out at the Institute des Mathematiques de Bordeaux (IMB, Bordeaux-France) during the visiting program of WB. 
Finally, the authors also wish to acknowledge the great opportunity for the exchange of ideas provided by the SHARK-FV 2023 workshop, where the main lines of research presented in this paper were initially conceptualized.

\bibliographystyle{abbrv}
\bibliography{paper_vem-dg}

\end{document}